\documentclass[amsmath,article,aps,fleqn]{revtex4-1}

\begin{document}
\title{On the Matrix-Element Expansion of a Circulant Determinant}
\author{Jerome Malenfant}
\date{\today}
\begin{abstract} The determinant of an $N \times N$ circulant matrix $M = {\rm CIRC}[x_0,  x_1, \ldots , x_{N-1}$] can be expanded in the form det$ ~M= \sum C_{a_0 a_1 \cdots a_{N-1}} x_{a_0} x_{a_1}\cdots x_{a_{N-1}}$.    By using the generating function of a restricted, 
mod-$N$ partition function, we derive a formula for the coefficients in this expansion as finite sums over products of binomial coefficients
with integer variables.
%SourceDoc 
\end{abstract}

\maketitle
\newtheorem{theorem} {Theorem}
\newtheorem{prop}{Proposition }
\newtheorem {lemma} {Lemma}
\newtheorem {corollary} {Corollary}
\newtheorem{proof} {Proof}

\section{Introduction}
A circulant matrix is a matrix of the form  
       \begin{eqnarray*}
      ~~~~~~~~~~~~~~~~~~~~~~~ \left(    \begin{array} {ccccc}
       x_0 & x_{N-1}  & \cdots & x_2 &x_1 \\ 
       x_1 & x_0  &\cdots & x_3 & x_2\\
       x_2 & x_1 & \cdots & x_4& x_3 \\
       \vdots & \vdots & ~ & \vdots & \vdots \\
        x_{N-2} & x_{N-3}  &\cdots &x_0 & x_{N-1} \\
        x_{N-1} & x_{N-2}  &\cdots &x_1 & x_0  \end{array} \right)  \equiv {\rm CIRC} \,[x_0, x_1, \ldots , x_{N-1}]     
         \end{eqnarray*} 
where successive columns are circular permutations of the first column.  (An alternative definition used by some authors is to take successive {\it rows} as circular permutations of the first row.)  It is a particular form of  a Toeplitz matrix, a matrix whose elements are constant along the diagonals.     The $N$ eigenvalues of the $ N \times N$ circulant matrix  above are
\begin{eqnarray}
x_0 + \omega^p x_1 + \omega^{2p} x_2 + \cdots + \omega^{p(N-1)} x_{N-1}, ~~p = 0, \ldots , N-1
\end{eqnarray}
where $\omega = e^{2 \pi i /N} $ is an $N$th root of unity.  The determinant of the matrix is therefore the product of these eigenvalues:
 \begin{eqnarray} 
       \left|    \begin{array} {ccccc}
       x_0 & x_{N-1}  & \cdots & x_2 &x_1 \\ 
       x_1 & x_0  &\cdots & x_3 & x_2\\
       x_2 & x_1 & \cdots & x_4& x_3 \\
       \vdots & \vdots & ~ & \vdots & \vdots \\
        x_{N-2} & x_{N-3}  &\cdots &x_0 & x_{N-1} \\
        x_{N-1} & x_{N-2}  &\cdots &x_1 & x_0  \end{array} \right|
        ~\equiv {\rm det} [x_0, x_1, \ldots , x_{N-1}] 
        = \prod_{p=0}^{N-1} \left( x_0 + \omega^p x_1  +\cdots + \omega^{(N-1)p}x_{N-1} \right).
       \end{eqnarray} 

Some examples of these determinants for small values of $N$ are:                          
  \begin{eqnarray*}
  {\rm det} [A,B,C]  &=& A^3 +B^3 + C^3 -3ABC;  \\~\\
{\rm det} [A,B,C,D] &=& A^4 - B^4  + C^4 -D^4  - 2A^2C^2 + 2B^2D^2- 4A^2BD + 4AB^2C -4BC^2D+ 4ACD^2; \\~\\
{\rm det} [A,B,C,D,E] &=&   ~~ A^5 + B^5  + C^5  + D^5 + E^5   \\ &&
                 -5A^3BE  -5 A^3CD  -  5AB^3C  - 5B^3DE  - 5AC^3E\\ 
                 &&     - 5BC^3D   - 5ABD^3 - 5CD^3E - 5ADE^3 - 5BCE^3\\
                 &&  +  5A^2B^2D  + 5A^2BC^2  + 5A^2CE^2  + 5A^2 D^2E  + 5 AB^2E^2  \\ 
                 &&     + 5AC^2D^2   + 5B^2C^2E + 5B^2CD^2  + 5BD^2E^2  + 5C^2DE^2 \\ 
                 && -5ABCDE .              
\end{eqnarray*}

The product on the right in eq.(2), when expanded out, contains $N^N$ terms.  There are considerably fewer terms in these examples, and show some apparent patterns in the coefficients.    There appears then to be considerable cancelations occurring in the expansion of eq. (2), as there must be, since all the complex terms cancel.  This  suggests that there is a simpler way to calculate these determinants.  In addition, expression (2) has the disadvantage of using a complex number, $\omega$, to find integer coefficients.\\  

This work was initiated then by the desire  to derive an explicit formula for the coefficients in the expansion of the determinant,  and to see if the apparent patterns in the coefficients for small $N$ appear in the general case as well.\\

Near the completion of these results I discovered the monograph on circulant matrices by Wyn-jones  \cite{Wyn}, where chapters 10 and 11 are devoted to the same problem that inspired this effort, that of finding a formula for the coefficients.   The formula in ref. [1], which was first stated in a 1951 article by Ore \cite{Ore},  bears some similarities to the one derived in this article; both formulas involve sums over the partitions of sets of integers, for example.  But there are significant differences between the two, and Wyn-jones' proof uses a completely different method than the one presented here.   In addition, it is likely that one formula is more efficient than the other for some coefficients and less efficient for others.  So it seems reasonable to present this alternative one as well.  Some of the other results in this article were proved in ref. [1], but again the proofs here are by different methods.   \\

\section{Preliminary results:  symmetries of the $C_{[a]}$ coefficients}

As shown in the examples in Section I, the determinants can be expanded as sums over products of the matrix elements:
\begin{subequations}
\begin{eqnarray}
    {\rm det ~} [  x_0, x_1, \ldots , x_{N-1}] &=& \sum_{0 \leq a_0 \leq  \cdots \leq a_{N-1} \leq N-1} 
     C_{a_0\cdots  a_{N-1}} x_{a_0}\cdots x_{a_{N-1}} .
\end{eqnarray}
The sum can also be written in the form
\begin{eqnarray}
    {\rm det ~} [  x_0, x_1, \ldots , x_{N-1}] &=& \sum_{0 \leq M_0 ,  \ldots , M_{N-1} \leq N} 
     C^*_{M_0 \cdots  M_{N-1}} x_0^{M_0 }\cdots x_{N-1}^{M_{N-1} } ,
\end{eqnarray}
\end{subequations}
where $ C^*_{M_0 \cdots M_{N-1}} = C_{a_0 \cdots a_{N-1}}, ~ M_0+ \cdots + M_{N-1} =N,$  and $M_m$ is the multiplicity of $m$, the number of times the integer $m$ occurs in the index set  $[a_0, \ldots, a_{N-1}]$.
 Let ${\cal S}_{[a_0, \cdots , a_{N-1}]}$ be the set of permutations of the index set $[a_0,\ldots , a_{N-1}]:  ~\sigma[a_0, \ldots, a_{N-1} ] \rightarrow   [\sigma_0, \ldots, \sigma_{N-1}]$.   From eq.(2), the coefficients in eq.(3a) are 
\begin{eqnarray}
     C_{a_0 \cdots a_{N-1}}  &=& \sum_{\sigma \in {\cal S}_{[a_0, \ldots, a_{N-1}]}} ~ \omega^{\sigma_1+2 \sigma_2+ \cdots + (N-1)\sigma_{N-1} }
     \end{eqnarray}
where the sum is over all elements of ${\cal S}_{[a_0, \cdots, a_{N-1}]}$.   One could of course use this expression to calculate these coefficients, but it is not very convenient, involving as it does a sum over a large number of complex numbers.  \\

We will use the notations $[a] \equiv [a_0, \cdots, a_{N-1}],~(a_0 \leq a_1 \leq \cdots \leq a_{N-1}), $ and $[\sigma] \equiv [\sigma_1, \ldots , \sigma_{N-1}]$, where $\sigma$ is an arbitrary permutation of $[a]$.  Let $\Omega^{(k)}_{[a]}$ be the subset of ${\cal S}_{[a]}$ defined by
\begin{eqnarray*}
\Omega^{(k)}_{[a]} = \left\{ ~\sigma \in {\cal S}_{[a]}  \left|  \right. \sigma_1 + 2 \sigma_2 + \cdots + (N-1) \sigma_{N-1} \equiv k {\rm ~mod~} N ~\right\}
\end{eqnarray*}
and let $A_{[a]}^{(k)}$ be the number of elements in $\Omega^{(k)}_{[a]}$, which we will denote as $\# \Omega^{(k)}_{[a]}$.  More generally, if $S$ is a finite set, then $\#S$ is the number of elements in $S$.  The elements of $\Omega^{(k)}_{[a]}$ will be referred to as ``$k$-mod permutations'' (of $[a]$).  Since $\omega^{mN+k} =\omega^k$, the sum in eq.(4) reduces to the form
\begin{eqnarray}
C_{[a]} &=& A_{[a]}^{(0)} + A_{[a]}^{(1)} \omega + A_{[a]}^{(2)} \omega^2 + \cdots + A_{[a]}^{(N-1)} \omega^{N-1}
\end{eqnarray}
with
\begin{eqnarray}
A^{(0)}_{[a]} + A^{(1)}_{[a]} + \cdots + A^{(N-1)}_{[a]} = \frac{N!}{M_0! \cdots M_{N-1}!} .
\end{eqnarray}
When all the $a$'s are equal, ($a_0  = a_1 = \cdots = a_{N-1} =a; ~M_a=N)$,  ${\cal S}_{[aa \cdots a]}$ has just one element, so there is only one term in the sum in eq.(4).  The corresponding coefficient is
    \begin{eqnarray*}
    C_{a \cdots a} = e^{2 \pi i a(1 +2+ \cdots + N-1)/N} = e^{\pi ia(N-1)} = (-1)^{a(N-1)}.
    \end{eqnarray*}
In the rest of this article  we will consider coefficients with more than one multiplicity. \\

We have $\omega^N=1$; more generally, $(\omega^d)^N=1$ for any integer $d$.  If $d \neq 0 $ mod $N$,  the set 
$\{\omega^d, \omega^{2d}, \ldots , \omega^{(N-1)d}, \omega^{Nd}=1\}$ consists of the zeros of the polynomial $z^N -1$:
\begin{eqnarray*}
z^N-1&=& (z-\omega^d)(z-\omega^{2d}) \cdots (z-\omega^{(N-1)d})  (z-1 ). 
\end{eqnarray*}
Expanding on the right, we have
\begin{eqnarray*}
    z^N-1 = z^N -(\omega^d + \omega^{2d} + \cdots + \omega^{(N-1)d}+1) z^{N-1} + \cdots ,
\end{eqnarray*}
and since the coefficient of $z^{N-1}$ on the left is zero we  get the identity,
\begin{eqnarray}
\sum_{q=1}^{N-1} \omega^{qd}  = -1 +N \delta^N_{d,0} , 
\end{eqnarray}
where the mod-$N$ Kronecker delta is defined as
\begin{eqnarray*}
   \delta^N_{d,0} = \left\{ \begin{array} {l} 1~{\rm ~if~} d \equiv 0 {\rm ~mod~} N, \\~\\ 0 ~~~{\rm otherwise}. \end{array}\right.
 \end{eqnarray*}

\begin{theorem}
 ~\\~\\
  {\it $C_{[a]} =0$ if $a_0 +a_1 +  \cdots + a_{N-1} \neq 0 {\rm ~mod~}N.$ }
\end{theorem}
~\\
(This is Proposition (10.4.3) in ref.[1].)\\
~\\
{\bf Proof}:   Let $[a]$ be such that 
\begin{eqnarray*}
a_0 + a_1 + \cdots +a_{N-1} \equiv X {\rm ~mod~} N,~~X \neq 0,
\end{eqnarray*}
and therefore, for any permutation $[\sigma]$ of $[a]$,
\begin{eqnarray*}
\sigma_0 + \sigma_1 + \cdots +\sigma_{N-1} \equiv X {\rm ~mod~} N.
\end{eqnarray*}
We can assume $X$ to be positive; (if not, we can make the replacement $X \rightarrow X' = N-|X|$).    If we add the second equation  repeatedly to the condition 
\begin{eqnarray*}
\sigma_1 + 2  \sigma_2 + \cdots + (N-1)  \sigma_{N-1}\equiv k {\rm ~mod ~}N,
\end{eqnarray*}
 for permutations in $\Omega^{(k)}_{[a]}$, we get
\begin{eqnarray*}
 \sigma_0 + 2  \sigma_1 + \cdots + (N-1)  \sigma_{N-2}&\equiv&( k+X) {\rm ~mod ~}N ,\\
 \sigma_{N-1} + 2  \sigma_0 + \cdots + (N-1)  \sigma_{N-3} &\equiv& (k+2X) {\rm ~mod ~}N ,\\
 & \vdots&\\
  \sigma_2 + 2  \sigma_3 + \cdots + (N-1)  \sigma_{1} &\equiv &(k+N-X) {\rm ~mod ~}N.
\end{eqnarray*}
There is therefore, for each $k$,  a sequence of one-to-one, onto mappings:
\begin{eqnarray*}
 \Omega^{(k)}_{[a]} \rightarrow \Omega^{(k+X)}_{[a]} \rightarrow \Omega^{(k+2X)}_{[a]} 
     \rightarrow \cdots \rightarrow \Omega^{(k+N-X)}_{[a]} \rightarrow \Omega^{(k)}_{[a]} 
\end{eqnarray*}
induced by circular permutations of $\sigma_0, \ldots , \sigma_{N-1}$.  Accordingly,
\begin{eqnarray*}
A^{(k)}_{[a]} = A^{(k+X)}_{[a]} = A^{(k+2X)}_{[a]} = \cdots = A^{(k+N-X)}_{[a]},~~{\rm for~}k=0,1, \ldots , X-1.
\end{eqnarray*}
If $X=1$, or if $X$ and $N$ are relatively prime, then $A^{(0)}_{[a]} =A^{(1)}_{[a]} = \cdots = A^{(N-1)}_{[a]}$, and the summation in eq.(5)
for $C_{[a]}$ equals zero as a result of the identity (7).  If on the other hand gcd$(X,N) >1$, then $MX \equiv 0 {\rm ~mod~}N$ for some $M<N$, and we have instead $X$ sets of equalities:
\begin{eqnarray*}
A^{(0)}_{[a]} &=& A^{(X)}_{[a]} = \cdots = A^{((M-1)X)}_{[a]},\\
A^{(1)}_{[a]} &=& A^{(1+X)}_{[a]} = \cdots = A^{(1 +(M-1)X)}_{[a]},\\
&\vdots& \\
A^{(X-1)}_{[a]} &=& A^{(2X-1)}_{[a]} = \cdots = A^{(N-1)}_{[a]}.
\end{eqnarray*}
In this case, eq.(5) becomes
\begin{eqnarray*}
C_{[a]} &=& ( ~ A^{(0)}_{[a]} +A^{(1)}_{[a]} \omega  + \cdots + A^{(X-1)}_{[a]} \omega^{X-1} ~  ) (~ 1 + \omega^X + \omega^{2X} + \cdots +\omega^{(M-1)X} ~) .
\end{eqnarray*}
But since $\omega^X= e^{2 \pi i X/N} = (e^{2 \pi i /M})^x$ for some integer $x \neq 0 \mod N$, the sum $  1 + \omega^X + \omega^{2X} + \cdots +\omega^{(M-1)X}$ equals zero as a result of the identity corresponding to eq. (7) with the substitutions $N \rightarrow M$ and $\omega \rightarrow \omega_M = e^{2 \pi i /M}$. $ \triangle$\\~\\

A number of corollaries follow from Theorem 1.

\begin{corollary}
 For even $N$, $C_{0123 \cdots N-1} = 0$.
\end{corollary}

\begin{corollary}
 For odd $N$, $C_{a_0 \cdots a_{N-1}} =0$ if the multiplicities are 2, 1, \ldots , 1, 0.
 \end{corollary}

\begin{corollary} 
$C_{a\cdots ab\cdots b} = 0$ if the multiplicities  $M_a$ and $M_b$ are relatively prime to $N$.
\end{corollary}
The proofs are straightforward and we omit them.\\

As noted, if $[a_0, \ldots , a_{N-1}]$ satisfies the condition, 
\begin{subequations}
\begin{eqnarray}
a_0 + a_1 + \cdots + a_{N-1} \equiv 0 {\rm ~mod~} N,
\end{eqnarray}
then any permutation $[\sigma]$ of $[a]$ also satisfies it.   Expressed in terms of the set of  multiplicities, 
$[M] = [M_0,M_1, \ldots , M_{N-1}]$, $ 0 \leq M_0, M_1, \ldots , M_{N-1}\leq N$, this condition becomes
\begin{eqnarray}
&&  0 \cdot M_0 + 1 \cdot M_1 + \cdots + (N-1) \cdot M_{N-1} \equiv 0 ~{\rm ~mod~} N,\\ 
&& {\rm subject~to~the~restriction} \nonumber \\ 
&&M_0 + M_1 + \cdots + M_{N-1} =N.
\end{eqnarray}
\end{subequations}
 Since we are only interested in non-zero coefficients, unless otherwise stated we will assume that all $[a]$'s and  $[\sigma]$'s satisfy (8a), and all $[M]$'s satisfy (8b,c).   \\
 
We define $P$ to be the circular permutation operator,  $P( \sigma_q) =\sigma_{q+1}$:
\begin{eqnarray}
P[\sigma] = P [\sigma_0, \sigma_1, \ldots, \sigma_{N-2}, \sigma_{N-1}]  &=& [ \sigma_1, \sigma_{2}, \ldots , \sigma_{N-2} ,\sigma_{N-1}, \sigma_0 ] .
\end{eqnarray}
Let $ \{ \sigma \} = \{\sigma_0, \ldots , \sigma_{N-1} \}$ denote the  equivalence class under $P$  of a $k$-mod permutation $[\sigma]$:
\begin{eqnarray*}
\{ \sigma \} = {\rm distinct~elements~ of ~the~set~} \{ ~[\sigma], P [\sigma], \ldots ,P^{N-1} [\sigma] ~\}.
\end{eqnarray*} 
If $[a]$, with multiplicities $M_0, M_1, \ldots , M_{N-1},$ is such that gcd$(M_0,M_1 \ldots , M_{N-1}) =1$, each equivalence class of the permutations of $[a]$ has $N$ elements  and each $A^{(k)}_{[a]}$ is $N$ times the total number of such equivalence classes.   
Therefore $C_{[a]} \equiv 0 $ mod $N$.   On the other hand, if the multiplicities have a common divisor $d$ and if a permutation in $\Omega^{(k)}_{[a]}$ takes the form  $[a_p a_q a_r \cdots a_s~ a_p \cdots a_s ~ a_pa_qa_r \ldots a_s]$ where the sequence $a_pa_qa_r\cdots a_s$ is repeated $d$ times, then the equivalence class containing this permutation has only $N/d$ distinct elements.  So, from eq.(5), we have

\begin{prop}
If gcd$(M_0,M_1, \ldots, M_{N-1} ) =d$, then $C_{a_0\cdots a_{N-1} } \equiv 0$ mod (N/d).
\end{prop}
~\\
  For an integer $n$ and a permutation $[~\sigma_0, \ldots, \sigma_{N-1} ~]$ we define the operations 
\begin{eqnarray*}
&&n + [~\sigma_0, \ldots , \sigma_{N-1}~] \equiv [~n+\sigma_0 , \ldots , n +\sigma_{N-1}~], \\
&&n \times [~ \sigma_0, \ldots, \sigma_{N-1} ~] \equiv [ ~ n \sigma_0, \ldots, n \sigma_{N-1} ~],
\end{eqnarray*}
where addition and multiplication is modulo $N$.  Further, we define $|\sigma|$ to be the permutation $[\sigma']$ of  $ [\sigma_0, \sigma_1,  \ldots , \sigma_{N-1} ]$ that puts it into ascending order:   $\sigma'_0 \leq \sigma'_1 \leq \cdots \leq \sigma'_{N-1}$.\\

\begin{prop}
~\\~\\
 {\rm (A)} $C_{a_0 \cdots a_{N-1}}= (-1)^{n(N-1)}C_{| n+[a_0, \cdots, a_{N-1}] | }.$ \\~\\
{\rm (B)}  $C_{a_0\cdots a_{N-1}} = C_{| n \times [a_0 \cdots a_{N-1}]|}$ {\it for n and N relatively prime.}.\\
\end{prop}
(This is Proposition 10.5.5 in ref [1].)\\~\\
{\bf Proof}:  \\
~\\
 (A) Condition (8a) is invariant under the transformation $\sigma_k \rightarrow (n+\sigma_k )~{\rm mod~} N$ for $n=1, \ldots , N-1$, while  
  $ \sigma_1+2 \sigma_2+ \cdots + (N-1)\sigma_{N-1} $  changes by
 \begin{eqnarray*}
\Delta \left( \sigma_1+2 \sigma_2+ \cdots + (N-1)\sigma_{N-1}  \right) = n~ \frac{(N-1)N}{2}.
\end{eqnarray*}
This transformation will then give a factor of $e^{2 \pi i \Delta/N} $ to each term on the right side of eq.(4), so
\begin{eqnarray*}
C_{|n+a_0,  \cdots ,n+ a_{N-1}|} = (-1)^{n(N-1)}   C_{a_0  \cdots   a_{N-1}} .
\end{eqnarray*} 
~\\
(B) The set $\{ 1, \omega^n, \omega^{2n}, \ldots , \omega^{(N-1)n} \}$ contains the same elements, up to order, as the set $\{ 1, \omega, \omega^2, \ldots , \omega^{N-1} \}$  when $N$ and $n$ are coprime.  We can therefore make the replacement $\omega \rightarrow \omega^n$ in all the terms and factors in the product on the right in eq.(2), since this amounts to just a re-ordering of the factors.   Eq.(4) then becomes
 \begin{eqnarray*}
     C_{[a]}  &=& \sum_{\sigma \in {\cal S}_{[a_0, \ldots, a_{N-1}]}} ~ (\omega^n)^{\sigma_1+2 \sigma_2+ \cdots + (N-1)\sigma_{N-1} }
     = \sum_{\sigma \in {\cal S}_{[a_0, \ldots, a_{N-1}]}} ~ \omega^{n\sigma_1+2n \sigma_2+ \cdots + (N-1)n\sigma_{N-1} }
      = C_{| (n \times [a]  |}. ~~~~ \triangle\\
     \end{eqnarray*}

In section IV below, we will use the two symmetry operations above to classify the coefficients, (or, equivalently, their multiplicative index sets), into ``additive multiplets'', in which the coefficients are related by the additive symmetry, and into ``super-multiplets'', in which they are related by the additive symmetry, the multiplicative symmetry, or a combination of both.   We will use the notation $\{ C_{[a]} \}$ or
 $\{C^*_{[M]}\} $ to represent the additive multiplet containing $C_{[a]}$, and $( C_{[a]})$ or $( C^*_{[M]}) $ to represent the corresponding super-multiplet.  As a result of Proposition 2, coefficients in an additive  multiplet or in a super-multiplet are equal, up to sign, and have the same multiplicities, up to order.   The number of super-multiplets is thus equal to the number of coefficients that need to be calculated to solve for the determinant. \\

As pointed out by Wyn-jones, condition (8a) is a necessary but not sufficient condition for a coefficient to be non-zero;   he gives two examples,  $C_{001335}$ and $C_{0000111368}$, which are zero despite satisfying (8a).   A sufficient (but not necessary) criteria for ``condition-(8)'' coefficients to be zero is given in Section IIID, as a corollary to Theorem 3.  \\

Wyn-jones  did prove that, for $N$ a prime number, condition (8a) {\it is} a sufficient condition for $C_{[a]}$ to be non-zero.  An alternate proof of this is given as a corollary of Theorem 2 below.\\

The  M\"{o}bius function  \cite{Moebius} $\mu(n)$ is defined as
\begin{eqnarray*}
 \mu (n) &=& \left\{ \begin{array} {l} ~~~1~~~~{\rm ~if~}n=1;\\~\\~~~0~~~~{\rm ~if~}n{\rm~has~one~or~more~repeated~prime~factors};\\~\\
                                                 (-1)^k~{\rm ~if~}n {\rm ~is~the~product~of~}k{\rm~distinct~prime~factors}. \end{array} \right.                                                                                       
 \end{eqnarray*}
Then
\begin{theorem}
 \begin{eqnarray*}
C_{[a]} =  \sum_{d | N}~ \mu\left( \frac{N}{d} \right) A^{(d)}_{[a]} 
\end{eqnarray*}
where, in the sum over the divisors of N, $A^{(N)}_{[a]}$ is to be identified with $ A^{(0)}_{[a]}$.
\end{theorem}
{\bf Proof}:     We define the mapping $n \cdot [\sigma] \rightarrow [\sigma']$, where  $n$ is an integer that is coprime to $N$, by $\sigma_q \rightarrow \sigma'_q  = \sigma_{nq {\rm ~mod~}N}$, (or just $\sigma_{nq}$, with mod $N$ multiplication understood).   Under this map, 
\begin{eqnarray*}
\sum_{q=0}^{N-1} q \sigma_q \rightarrow \sum_{q=0}^{N-1} q \sigma_{nq} = \sum_{q'=0}^{N-1} (n^{-1}q')\sigma_{q'}
  = n^{-1} \sum_{q'=0}^{N-1}q' \sigma_{q'}
\end{eqnarray*}
where $n^{-1}$ is the multiplicative inverse of $n$.  $[\sigma']$ is therefore an element of $\Omega^{(n^{-1}k )}_{[a]}$ and $n$ induces a mapping:
\begin{eqnarray*}
 \Omega^{(k)}_{[a]} \longrightarrow_n \Omega^{(n^{-1}k)}_{[a]}.
 \end{eqnarray*}
If $n$ is coprime to $N$ then so is $n^{-1}$, which induces the inverse map:
\begin{eqnarray*}
\Omega^{(k)}_{[a]} ~_{n^{-1}}\longleftrightarrow_n~ \Omega^{(n^{-1}k )}_{[a]} .
\end{eqnarray*}
This is therefore a one-to-one onto mapping between $\Omega^{(k)}_{[a]}$ and $\Omega^{(nk )}_{[a]}$.  These two sets have the same number of elements then, and $A^{(k)}_{[a]} =A^{(nk )}_{[a]}$ for all $n$ coprime to $N$.  In particular,  $A^{(1)}_{[a]} =A^{(n)}_{[a]}$ for all $n$ coprime to $N$ and, if  $N$ is a prime number,  $A^{(1)}_{[a]} = A^{(2)}_{[a]} = \cdots = A^{(N-1)}_{[a]}$.\\

Let  $\{ d_0, d_1, \ldots , d_K \}$ be the set of divisors of $N$:  $1=d_0 <d_1 < \cdots < d_K=N$.  With the replacement 
$A^{(0)}_{[a]} \rightarrow A^{(N)}_{[a]}$, we separate the sum in eq.(5)  in the manner below: 
\begin{eqnarray*}
C_{[a]} = A^{(N)}_{[a]} + \sum_{{\rm gcd}(N,q)= 1}^{N-1} A^{(q)}_{[a]} \omega^q +  \sum_{{\rm gcd}(N,q)= d_1}^{N-d_1} A^{(q)}_{[a]} \omega^q
+ \cdots +  \sum_{{\rm gcd}(N,q)= d_{K-1}}^{N-d_{K-1}} A^{(q)}_{[a]} \omega^q.
\end{eqnarray*}
In the sum corresponding to gcd$(N,q) = d_p$, the summation variable $q$ equals  $kd_p$ for some $k$ coprime to $N$.  By the relation above,  $A_{[a]}^{(kd_p)} = A_{[a]}^{(d_p)}$ for this sum.   The $A$ coefficients can therefore be brought outside:
\begin{eqnarray*}
C_{[a]} = A^{(N)}_{[a]} + \sum_{p=0}^{K-1} A^{(d_p)}_{[a]}  \sum_{{\rm gcd}(N,k)= 1 }^{N/d_p-1} \omega^{kd_p}
              =   \sum_{p=0}^K A^{(d_p)}_{[a]}  \sum_{{\rm gcd}(N,k)= 1 }^{N/d_p-1} \omega^{kd_p}.
\end{eqnarray*}
  The result then follows by noting that the sum over $k$
\begin{eqnarray*}        
 \sum_{{\rm gcd}(N,k) =1}^{N/d_p-1} \omega^{kd_p}
  \end{eqnarray*}
 is the M{\"o}bius function $\mu(N/d_p)$\cite{Moebius}.    $\triangle$\\
~\\

The coefficient $C_{[a]}$ can thus be expressed as a sum of integers, without the use of complex numbers.  This expression is 
not very useful in calculating $C_{[a]}$ since it still leaves us with the need to count the number of elements in the $\Omega^{(k)}_{[a]}$ sets, although for fewer of these sets than appear in eq.(5).  \\

Using this theorem though we can prove that $C_{[a]}$ is non-zero when $N$ is prime, provided of course that the index set $[a]$ satisfies condition (8a).  Changing back to the previous notation $A^{(0)}_{[a]}$, we have
 \begin{corollary}
 For $N$ prime,
\begin{eqnarray*}
  C_{[a]} =  \frac{N}{N-1} ~ \left( ~A^{(0)}_{[a]} - ~ \frac{(N-1)!}{M_0! \cdots M_{N-1}!} ~\right) \neq 0.
\end{eqnarray*}
\end{corollary}
{\bf Proof}: The divisors of $N$ consist of just 1 and $N$, so from Theorem 2
\begin{eqnarray*}
  C_{[a]} =  A^{(0)}_{[a]} -A^{(1)}_{[a]} .
\end{eqnarray*}
We also have in this case that 
\begin{eqnarray*}
 \frac{N!}{M_0! \cdots M_{N-1}!}=  A^{(0)}_{[a]} + A^{(1)}_{[a]} + \cdots + A^{(N-1)}_{[a]} =A^{(0)}_{[a]} + (N-1)A^{(1)}_{[a]}
\end{eqnarray*}
from the proof of Theorem 2 and from eq.(6).  The first equality  follows from these relations.   To show that $C_{[a]} \neq 0$,  $C_{a\cdots a } = (-1)^{(N-1)a} \neq 0$, so we need consider  just the cases where all the multiplicities are less than $N$.   For these, $A^{(0)}_{[a]} \propto N$, so $A^{(0)}_{[a]} = N B_{[a]}$ for some integer $B_{[a]}$.   If $C_{[a]}=0$, then
\begin{eqnarray*}
 B_{[a]} =  \frac{(N-1)!}{N~M_0!\cdots M_{N-1}!} .
\end{eqnarray*}
$N$ must therefore divide $(N-1)!$.  But this is a contradiction since $N$ is prime and all the factors in $(N-1)!$ are less than $N$.  $\triangle$.\\
~\\

Corollary 5 in Section IIIB below, also a corollary to Theorem 3, extends the criteria for $C_{[a]} \neq 0$, i.e., $N$ being a prime number, to $N$ not dividing $\frac{(N-M_0-1)!}{M_1!}$.\\

~\\
\section{ Derivation of a Formula for $C_{a_0 \cdots a_{N-1}} $}

\subsection{ Definitions and statement of the theorem}

We will denote the set of partitions of the integer $p$ by ${\cal P}_p$, with a particular partition consisting of $j$ parts denoted by $Z= (z_1, \ldots , z_j)$ or, equivalently, by the set of multiplicities $[k_1, \ldots , k_p]$: 
\begin{eqnarray*}
 p &=& z_1 + \cdots + z_j \\
 &=& k_1 + 2 k_2 + \cdots + p k_p;~~k_1+ \cdots +k_p =j.
\end{eqnarray*}

The set of partitions of a set of not-necessarily-distinct integers $\{a,b, \ldots , c\} $ will be denoted by ${\cal P}\{ ab \cdots c\}$.  For a particular partition $\Theta \in {\cal P}\{ab \cdots c\}$, $\Theta = ( \theta_1) \cdots ( \theta_j)$, where each $\theta_k$ is a subset of 
$\{a,b,\ldots , c\}$.  We will express these subsets, for example, as $\theta = abc$, in contexts where there should be no confusion. 
  As an illustration, the set of partitions of the set  $\{a,a,b,c\}$, where $a \neq b \neq c$, is

 \begin{eqnarray*}
{\cal P}\{aabc\} &=& \left\{ \frac{}{}(aabc),~(aab)(c),~(aac)(b),~(abc)(a),~(aa)(bc),~(ab)(ac),\right. \\
&& \left. ~~ (aa)(b)(c),~  (ab)(a)(c),~(ac)(a)(b),~(bc)(a)(a),~    (a)(a)(b)(c) \frac{}{} \right\}.
 \end{eqnarray*}

For a partition $\Theta$  and for any subset $\theta$ of $\{a,b, \ldots ,c\}$, $\kappa_{\theta}$ will denote the multiplicity 
of $\theta$ in the partition; for $\theta \in \Theta$, $m_{a}^{(\theta)}$ will denote the multiplicity of  the element $a$ in $\theta$.  Finally, the sum of the elements of $\theta$ will be denoted as the trace of $\theta$, Tr$(\theta)$:
 \begin{eqnarray*}
 {\rm Tr} (\theta) = \sum_{a\in \theta}~a.
 \end{eqnarray*}

    The map given by $z_k = \# \theta_k$ sends the set partition $\Theta= (\theta_1)\cdots ( \theta_j) $ to the integer partition $Z= (z_1, \ldots , z_j)$.  Two partitions $\Theta$ and $\Theta'$ in ${\cal P} \{ ab \cdots c\} $ are equivalent iff they map to the same partition $Z$.   In this case, each $\theta_k$ and $\theta'_k$ have the same number, $z_k$, of elements.  The $\Theta$ equivalence classes can then be labeled by the partition $Z$, and each $\Theta \in \{Z\}$ for some $Z$. \\
 
We will also use the multinomial coefficient $(p; k_1, \ldots, k_p)^*$, defined as \cite{Abram} 
\begin{eqnarray*}
   (p;~k_1,k_2, \ldots, k_p)^* &=& \frac{p!}{1^{k_1}k_1! \, 2^{k_2}k_2! \cdots p^{k_p}k_p!},~~
    p  = k_1+2k_2 + \cdots + pk_p, \\
\end{eqnarray*}
and  the Heaviside step function $H[n]$:
 \begin{eqnarray*}
H[n] = \left\{  \begin{array} {l} 1 {\rm ~~if~} n\geq 0, \\ ~\\ 0 {\rm ~~if~} n <0 .\end{array} \right.\\
\end{eqnarray*}  

For convenience, we will write the coefficient $C_{a_0 \cdots a_{N-1}}$ as $C_{0\cdots 0\,1\cdots 1 \, A_1 \cdots A_p A_{p+1}}$
with the 0's and 1's separated off, with $A_1 = a_{M_0+M_1},~A_{p+1} = a_{N-1}$, and $p= N-M_0-M_1-1 \geq 0 $.  The $A$ indices lie between 2 and $N-1$, with multiplicities $M_2, M_3, \ldots , M_{N-1}$. \\

With these definitions,  Theorem 3 is the statement below.

\begin{theorem}
~\\~\\
The coefficients  with multiplicities $0 \leq M_0, M_1, \ldots M_{N-1} \leq N-1$ in the expansion (3a) are given by 

  \begin{eqnarray*}
C_{0\cdots 0 \, 1 \cdots 1\, A_1 \cdots A_p A_{p+1} } &=& \delta^N_{ 0,M_1 +2M_2 + \cdots + (N-1)M_{N-1} }~  (-1)^{N -M_0 -1}  N
 \\
  && \times  \left\{~ \frac{(N-M_0-1)!}{M_1!M_2! \cdots M_{N-1}!} +   \frac{1}{M_{A_{p+1}}} \sum_{\Theta \in {\cal P} \{A_1  \cdots A_p  \} }  
                       \prod_{\theta }     \frac{1}{\kappa_{\theta}!}  ~  \prod_{k=1}^j 
                       ~\frac{   (z_k-1)! } {m_{2}^{(\theta_k)}! \cdots m_{N-1}^{(\theta_k)}!  }  \right. \\
   && ~~~~  \times  ~ \sum_{\mu=1}^j  (-N)^{\mu}  \sum_{\lambda_1, \ldots , \lambda_j =0,1 } \delta_{\mu, \lambda_1 + \cdots + \lambda_j } 
                     ~ H \left[ M_1 - \sum_{t=1}^j \lambda_t X(\theta_t) \right]     \\ 
    &&\left.  ~~~    \times   \left( \begin{array} {c} N-M_0-1 - \sum_{t=1}^j \lambda_t (X(\theta_t) +z_t ) \\
            ~\\   M_1 - \sum_{t=1}^j \lambda_t X(\theta_t)  \end{array} \right) 
                \prod_{s=1}^j    \left( \begin{array} {c} X(\theta_s) + z_s -1\\~\\ z_s-1 \end{array} \right)^{\lambda_s} \right\}  \\
        \end{eqnarray*}
where $j$ is the number of parts in the partition $ \Theta$, $z_s= \# \theta_s$ is the number of elements in the part $\theta_s$, and   
  \begin{eqnarray*}
     && X(\theta_s) \equiv - Tr (\theta_s) {\rm ~mod~} N
       \end{eqnarray*}
such that $ 0 \leq X(\theta_s)  \leq N-1$.\\
\end{theorem}
   
For each partition, the sum over the $\lambda$'s has $2^j -1$ terms.  However, the presence of the step function in the expression means that, for a term to be non-zero, it must satisfy the condition
  \begin{eqnarray*}
  M_1 -  \sum_{t=1}^j \lambda_t X(\theta_t)  \geq 0.
  \end{eqnarray*}         
 As a consequence, for many coefficients, many if not most of the terms in the sum over the $\lambda$'s are zero.  \\
 
 As an example, and as a comparison to the method used in \cite{Wyn}, we will calculate the coefficient  $C_{0011113788}$ using  this theorem.   In this case, the  $\kappa,~m$ and the $M_{k>1}$ multiplicities  all equal 1 except for $M_{A_{p+1}}$.   We have then
   \begin{eqnarray*}
C_{0011113788}
         &=&-10
        \left\{~ \frac{7!}{4!1!1!2!}  + \frac{1}{2}  \sum_{\Theta \in {\cal P} \{378 \} }    (z_1-1)! \cdots (z_j-1)!   
         ~ \sum_{\mu=1}^j  (-10)^{\mu}  \sum_{\lambda_1, \ldots , \lambda_j =0,1 } \delta_{\mu, \lambda_1 + \cdots + \lambda_j }\right. \\ 
                &&\left.~~~~~~~~  \times       H \left[ 4 - \sum_{t=1}^j \lambda_t X(\theta_t) \right ]     \left( \begin{array} {c} 7 - \sum_{t=1}^j \lambda_t (X(\theta_t)  + z_t ) \\
            ~\\   4 - \sum_{t=1}^j \lambda_t X(\theta_t)  \end{array} \right)
                \prod_{s=1}^j    \left( \begin{array} {c} X(\theta_s) + z_s -1\\~\\ z_s-1 \end{array} \right)^{\lambda_s} \right\}  .
   \end{eqnarray*}              
 The partitions of the set $\{3,7,8\}$  and the values of $(z_1, \ldots , z_j) $ and $(\lambda_1, \ldots , \lambda_j)$ for the non-zero terms are:
 \begin{eqnarray*}
~  (378): &&~z_1=3;~X(378)=2;~ \lambda_1 =1;\\
~\\
~(37)(8) :&&~(z_1,z_2)= (2,1);~X(37)= 0,~X(8)=2; ~(\lambda_1,\lambda_2) = (1,0),(0,1),(1,1);\\
 ~\\
~(38)(7): &&~(z_1,z_2) = (2,1);~X(38) = 9,~ X(7)=3; ~(\lambda_1, \lambda_2)= (0,1) ;\\
   ~\\     
~(3)(7)(8) :&& ~(z_1,z_2,z_3) = (1,1,1);~X(3)=7,~X(7)=3,~X(8)=2;~ (\lambda_1, \lambda_2,\lambda_3) = (0,1,0),(0,0,1).
\end{eqnarray*}
Then
         \begin{eqnarray*}
 C_{0011113788}
         &=& -10  \left\{~ 105  -5 \left[     2  \left( \begin{array} {c} 2 \\~\\   2 \end{array} \right) \left( \begin{array} {c} 4\\~\\ 2 \end{array} \right)  
        +   \left( \begin{array} {c} 5 \\ ~\\   4   \end{array} \right)   \left( \begin{array} {c} 1\\~\\ 1 \end{array} \right)
                   + \left( \begin{array} {c} 4 \\   ~\\   2  \end{array} \right)  \left( \begin{array} {c} 2 \\~\\   0 \end{array} \right)
                    +       \left( \begin{array} {c} 3\\   ~\\   1  \end{array} \right)  \left( \begin{array} {c} 3 \\~\\   0 \end{array} \right)
                    \right. \right.\\
                 && \left.  \left. ~~~~~~~~~~~~~~~~~~~~   + \left( \begin{array} {c} 3\\ ~\\   1 \end{array} \right)   \left( \begin{array} {c} 3\\~\\ 0 \end{array} \right)
                   + \left( \begin{array} {c} 4  \\   ~\\   2  \end{array} \right)   \left( \begin{array} {c} 2 \\~\\   0  \end{array} \right)     \right]                   
                     +   (-10) \left( \begin{array} {c} 2  \\~\\   2 \end{array} \right)  \left( \begin{array} {c} 1\\~\\ 1 \end{array} \right)
                      \left( \begin{array} {c} 2 \\~\\   0 \end{array} \right)\right\}\\
          &=&  -10 \{ ~105 - 175 + 50~ \} = 200 ,\\    
   \end{eqnarray*}              
in agreement with the value obtained in ref. [1]. \\

By way of comparison, the method used by Wyn-jones involves partitions of the partial index set $\{1,1,1,1,3,7,8,8\} $ that are make up entirely of ``null multisets'';   i.e. (in my notation) subsets $\theta$ of $\{1,1,1,1,3,7,8,8\} $ such that Tr($\theta) \equiv 0 $ mod $N$.   Such a partition consisting of $j$ subsets contributes to the $N^{j}$ term in the expansion, unlike the calculation above, where each partition can contribute to several $N^{j}$ terms.    Although ${\cal P}(11113788)$ contains many more partitions than ${\cal P}(378)$, (4,140 compared to 5), the calculation of this coefficient in ref. [1] involves only the 5 partitions,
\begin{eqnarray*}
(37)(118)(118),~(1117)(1388),~(37)(111188),~(118)(11378),~(11113788).\\
\end{eqnarray*}

 \subsection{Proof of Theorem 3}

  We will assume in the following that the indices of all coefficients considered satisfy conditions (8) and will omit writing out explicitly the mod-$N$ Kronecker delta $\delta^N_{0, M_1 + \cdots + (N-1)M_{N-1}}$ in the formulas for the coefficient. \\

Consider first the case $M_A=1$ and the coefficient $C_{0 \cdots 0\,1 \cdots 1A}$.   In this case we have $M_0+M_1+1 =N$ and $0 \cdot M_0 + 1\cdot M_1 + A \cdot 1 \equiv 0$, so  $A= N-M_1$.\\

With the freedom we have from the circular-permutation symmetry we can choose, when counting the number of elements in each set 
$\Omega^{(k)}_{0\cdots 0 1 \cdots 1 (N-M_1)}$,  to consider only permutations of $[0 \cdots 0 \, 1 \cdots 1\,(N-M_1)]$  for which 
\begin{eqnarray*}
&&\sigma_0 = A= N-M_1, \\
&&\sigma_{p_1}=\sigma_{p_2} = \cdots = \sigma_{p_{M_1}}=1,~1 \leq p_1<p_2< \cdots < p_{M_1} \leq N-1,\\
&&{\rm all~ other~} \sigma_q =0. 
\end{eqnarray*}
~\\
The number of elements in $\Omega^{(k)}_{0\cdots 0\,1 \cdots 1 \,(N-M_1)}$ is $N$ times the number of such permutations that satisfy the condition\\
\begin{eqnarray*}
k \equiv 0\cdot \sigma_0 + 1 \cdot \sigma_1 + \cdots + (N-1) \cdot \sigma_{N-1} \equiv p_1 + p_2 + \cdots + p_{M_1} .
\end{eqnarray*}
~\\
$A^{(k)}_{0 \cdots 0\, 1 \cdots 1 \,(N-M_1)}$ therefore equals  $N$ times the number of partitions (mod $N$) of $k$ into exactly $M_1$ distinct parts less than or equal to $N-1$.   We define $Q(n;X,b)$ to be the restricted, mod-$N$ partition function \\
\begin{eqnarray*}
Q(n;X,b) = \#\{~ (p_1,p_2,\ldots ,p_X)~|~ p_1 + \cdots + p_X \equiv  n {\rm ~mod~} N;
                                  ~1 \leq p_1 < \cdots < p_X \leq b<N\};
\end{eqnarray*}
i.e., the number of partitions mod $N$ of $n$ into $X$, distinct parts, with each part less than or equal to $b$.
For $n=k,~X= M_1,$ and $b=N-1$, each such partition corresponds to a circular-permutation equivalence 
class $\{N-M_1, \sigma_1, \ldots , \sigma_{N-1} \}$.  $C_{0\cdots 0 \, 1 \cdots 1\, (N-M_1)}$ is then, from eq.(5)\\
\begin{eqnarray*}
C_{0 \cdots 0 1\cdots 1 (N-M_1)}= \sum_{k=0}^{N-1} A^{(k)}_{0 \cdots 0 1 \cdots 1(N-M_1)}  \omega^{k} 
                                                      = N \cdot \sum_{k=0}^{N-1} Q(k;M_1,N-1) \omega^k .
\end{eqnarray*}
By inspection, it can be seen that the partition function $Q(n;X,N-1)$ has the generating function
\begin{eqnarray*}
\sum_{X=0}^{N-1} \sum_{n=0}^{N-1} Q(n;X,N-1) y^X \omega^n &=& \prod_{q=1}^{N-1}(1+y \omega^q )=1 + \cdots 
    + y^s( \omega^{p_1 + \cdots +p_s} + \cdots ) + \cdots  .
 \end{eqnarray*}
Now referring back to eq.(2), the product over $q$ in the equation above is, except for the missing $q=0$ factor, the determinant of an $N \times N$ circulant matrix with elements $1,k,0,\ldots ,0$.  Using this fact and evaluating the determinant, we find
\begin{eqnarray*}
 \prod_{q=1}^{N-1} (1+y\omega^q)  
      = \frac{1}{1+y}~
      \left|    \begin{array} {cccccc}
       1 & 0 &0&  \cdots &0  &y \\ 
       y & 1  & 0& \cdots  & 0 & 0\\
       0 & y & 1 & \cdots & 0& 0 \\
       \vdots & ~ & \ddots  &  \ddots & ~  & \vdots \\
        0 & 0  &~ & y  &1 & 0 \\
        0 & 0  &0 & \cdots  &y & 1  \end{array} \right|_{N \times N}
                = \sum_{m=0}^{N-1} (-y)^m.
\end{eqnarray*}
Then
\begin{eqnarray*}
\sum_{n=0}^{N-1} Q(n;X,N-1) \omega^n =(-1)^X.
\end{eqnarray*}
and our result for this coefficient is 
\begin{eqnarray*}
C_{0 \cdots 0 1\cdots 1 (N-M_1)} = (-1)^{M_1} N.
\end{eqnarray*}

Next we consider the coefficient $C_{0 \cdots 0 \,1 \cdots 1 \,A_1A_2 \cdots A_pA_{p+1}},$ with $p =N-M_0-M_1-1$.  
For given sets $\{p_1, \ldots , p_{M_1} \}$ and $\{q_1, \ldots, q_{p}\}$ of, respectively,  $M_1$ and $p$ distinct integers between 1 and $N-1$, we assign the permutation sending $[0 \cdots 0 \,1 \cdots 1 \, A_1 \cdots A_{p+1}]$ to $ [\sigma_0 \cdots \sigma_{N-1}]$ as follows:
\begin{eqnarray*}
&&  \sigma_0 =A_{p+1};\\
&&\sigma_{p_1}= \cdots =\sigma_{p_{M_1}} =1,~ 1 \leq p_1 < \cdots < p_{M_1} \leq N-1; \\
&& \sigma_{q_1} =A_1, ~\sigma_{q_2}=A_2,~ \ldots , ~\sigma_{q_{p}} =A_p;
~1 \leq q_1, ~q_2, \ldots ,~ q_p \leq N -1, ~q_i \neq q_j {\rm ~for~all~} i \neq j.
\end{eqnarray*}
The requirement that a particular permutation be an element of  $\Omega^{(k)}_{0\cdots 0\,1 \cdots 1 \, A_1 \cdots A_{p+1}}$ is now
\begin{eqnarray*}
  p_1 + \cdots + p_{M_1}\equiv   k - A_1 q_1 - A_2 q_2 -\cdots -A_p q_p .
 \end{eqnarray*}

For the above set of distinct integers $\{q_j\},$  $1 \leq q_1, \ldots , q_p \leq N-1$, we  define $Q(n;X,b~|~q_1, \ldots , q_p)$ to be the mod-$N$ restricted partition function, 

\begin{eqnarray*}
&&Q (n; X,b~|~ q_1, \ldots , q_p)\\~\\
&&~~~~~~~~ = \#\left\{~\frac{}{}  (p_1,\ldots ,p_X)~|~ p_1 + \cdots + p_X \equiv  n {\rm ~mod~} N;~1 \leq p_1 < \cdots < p_X \leq b;
                                        ~p_i \neq q_j ~{\rm for~all~} i {\rm ~and~}j~ \right\} .
\end{eqnarray*}
 As before, $Q$ is the number of partitions mod $N$ of $n$ into $X$, distinct parts less than or equal to $b$, but now with the additional restrictions that none of the parts are equal to any $q_j$.  The set of all permutations in $\Omega^{(k)}_{0\cdots 0 \,1\cdots 1 \, A_1 \cdots A_{N-1}}$ is then the sum of $N \cdot Q (k-A_1q_1 -\cdots - A_pq_p; M_1,N-1~|~ q_1, \ldots , q_p)$ over all allowed values of the $q$'s,  divided by the factorials of the multiplicities to avoid overcounting:

  \begin{eqnarray*}
A^{(k)}_{0\cdots 0\,1 \cdots 1\, A_1 \cdots  A_{p+1}} &=&\left.  \frac{ N}{ M_2!  \cdots M_{N-1}! }   ~ \sum_{1 \leq  q_1, \ldots,  q_p \leq N-1}
                                                 Q (k-A_1q_1 -\cdots - A_pq_p; M_1,N-1~|~ q_1, \ldots , q_p) \right|_{q_i \neq q_l} .
 \end{eqnarray*}
 I.e. the sum is over all lattice points $(q_1, \ldots , q_p),~ q_i \in [1, N-1]$ such that no coordinates are equal.  In the following, we will omit explicitly writing out the limits of the summation over the $q$ variables, with the understanding that the limits will, unless otherwise indicated, always be from 1 to $N-1.$ \\
 ~\\
 ~[Note:  As pointed out in a previous section, in the case where the multiplicities of an index set all have a common factor $d$, there are some equivalence classes that have only $N/d$ distinct elements.  The expression above is however valid in those cases as well; see Appendix B for a discussion of this issue. ]\\ 
 
 Therefore, again using eq.(5),
\begin{subequations} 
  \begin{eqnarray}
C_{0\cdots 0\,1 \cdots 1\, A_1 \cdots  A_{p+1}} 
        &=&  \frac{ N}{ M_2!  \cdots M_{N-1}! } \sum_{k=0}^{N-1} ~ \omega^k  \sum_{q_i \neq q_j}
                                                 Q (k-A_1q_1 -\cdots - A_pq_p; M_1,N-1~|~ q_1, \ldots , q_p)  \nonumber \\
        &=&  \frac{ N}{ M_2! \cdots M_{N-1}! } ~ \sum_{ q_i \neq q_l} \omega^{ A_1q_1  +\cdots + A_pq_p } 
                ~\sum_{k'=0}^{N-1}~ \omega^{k'} ~ Q (k'; M_1,N-1~|~ q_1, \ldots , q_p)                                  
 \end{eqnarray}
 where we've changed the summation index to $k' = k -Aq_1 - \cdots - Aq_{M_A-1}$; since $Q$ uses modular arithmetic, this does not change the limits of the sum over $k$. \\

  The generating function for this restricted partition function can be obtained by omitting the factors $(1+y\omega^{q_j})$ for $j= 1, \ldots ,p$  in the product in the previous section for $Q(n;X,b)$'s generating function or, equivalently, by dividing that product by these factors:  
  
\begin{eqnarray*}
\sum_{X=0}^{N-1} y^X  \sum_{n=0}^{N-1} Q (n;X,N-1~|~    q_1, \ldots , q_p) \omega^n 
                           &=& \sum_{m=0}^{N-1} (-y)^m ~\frac{1}{(1+y\omega^{q_1})\cdots (1+y\omega^{q_p})}.\\
\end{eqnarray*}
\begin{lemma}
~\\~\\
Let  $\omega = e^{2 \pi i/N}$ and let $\{q_1, \ldots , q_p~| ~1 \leq q_i \leq N-1\}$ be a set of p distinct integers. Then 
  \begin{eqnarray*}
\sum_{m=0}^{N-1} (-y)^m  ~  \prod_{s=1}^{p} \frac{1} { 1+y\omega^{q_s} }
                 &=&  \sum_{m=0}^{N-1-p} (-y)^m   ~ 
                 \sum_{0 \leq \kappa_1, \ldots , \kappa_p \leq m}  \omega^{\kappa_{1} q_1  +\kappa_2 q_2 + \cdots +
                 \kappa_p q_p}~H[ m- ( \kappa_1 + \cdots + \kappa_p)] .
 \end{eqnarray*}

\end{lemma}
~\\
 For the proof of this lemma, and that of Lemma 2 below, see Appendices C and D, respectively.\\
 
  Lemma 1 allows us to write eq.(10a) as, 
  \begin{eqnarray}  
          C_{0\cdots 0 \,1\cdots 1\, A_1\cdots A_{p+1}} 
                     &=&  \frac{ (-1)^{M_1}  N}{M_2! \cdots M_{N-1}!} ~ \sum_{ q_i \neq  q_l }~  \sum_{0 \leq \kappa_1, \ldots , \kappa_p \leq M_1 } 
                                  ~ \omega^{(\kappa_1+A_1) q_1  + \cdots + (\kappa_p+A_p) q_p} H[M_1- ( \kappa_1 + \cdots + \kappa_p) ] . \nonumber 
    ~\\
 \end{eqnarray}
We plan to use eq.(7) to perform the summations over the $q$ variables in eq.(10b).  This is complicated by the fact that the sum in (7) is over the complete range of $q$, from 1 to $N-1$, but each $q$-sum above, while having this range, is restricted by the condition that $q_i \neq q_l$ for all $i \neq l$.   We will therefore express the restricted, $p$-dimensional sum in eq.(10b) by a sum over unrestricted, lower-dimensional $q$-sums.  An elementary example of this is
\begin{eqnarray*}
\sum_{x<y}~ h(x,y) = \frac{1}{2} \left\{ \sum_{x,y}~ h(x,y) -  \sum_{x} ~h(x,x)  \right\}
\end{eqnarray*}
for a symmetric function $h(x,y)$.  Lemma 2 below generalizes this example to higher dimensions.\\

We first define the $p$-dimensional vector ${\bf q}_{_Z}$:  For each part $z_k$ of a partition $Z$ of $p$, with $j$ parts, 
 the $p$-dimensional vector ${\bf u}_k$ has components 
\begin{eqnarray*}
  u^i_k = \left\{ \begin{array} {l} 1 {\rm ~~if~~} z_1+\cdots + z_k \leq i < z_1 + \cdots + z_{k+1} ; \\~\\
                                                        0 {\rm ~~otherwise.} \end{array} \right.
 \end{eqnarray*}
Then for variables $(q_1,q_2, \ldots , q_j)$,  ${\bf q}_{_Z}$ is defined as
\begin{eqnarray*}
{\bf q}_{_Z} = q_1 {\bf u}_1 + \cdots q_j {\bf u}_j  &=& (q_1, \ldots , q_1, q_2, \ldots , q_2, ~\ldots \ldots ~, q_j, \ldots , q_j) .\\
 && ~~~\,   \longleftrightarrow ~~~~~ \longleftrightarrow~~~~~~~~~~~~~~~~~~~\longleftrightarrow  \\
 &&~~~~~~z_1~~~~~\,~~~~z_2~~~~~~~~~\cdots ~~~~~~~~~~z_j
\end{eqnarray*}
${\bf q}_{_Z} $ thus lies on a $j$-dimensional sublattice determined by the partition $Z$. 
 
\begin{lemma}
~\\~\\
  Let $h(q_1,\ldots ,q_p)=h({\bf q})$ be a completely symmetric function of the variables $q_1, \ldots, q_p$.  Then
\begin{eqnarray*}
     \sum_{1 \leq q_1< \cdots <  q_p \leq M}   h({\bf q})   
     &=& \frac{1}{p!} \sum_{Z \in {\cal P}_p}  (-1)^{p+j } ~  (p; k_1,k_2, \ldots , k_p)^*
      ~ \sum_{1 \leq q_1, \ldots , q_j \leq M} h( {\bf q}_{_Z}). \\
\end{eqnarray*}
\end{lemma}

The function $ f(q_1,\ldots, q_p)$ that appears in eq.(10b) above,
\begin{eqnarray*} 
        f(q_1,\ldots, q_p) 
                     =\omega^{A_1 q_1 + \cdots  +  A_pq_p} ~ \sum_{0 \leq \kappa_1, \ldots , \kappa_p \leq M_1 } 
                                  ~ \omega^{\kappa_1 q_1  + \cdots + \kappa_p q_p} H[M_1-( \kappa_1 + \cdots + \kappa_p)],
 \end{eqnarray*}
 is however not a symmetric function on account of the factor in front.  To correct this we divide the sum over $q_1, \ldots , q_p$ in eq.(10b) into $p!$ subregions:  
  \begin{eqnarray*}
 \sum_{q_i \neq q_l} f(q_1, \ldots , q_p) &=&  \left\{  \sum_{ q_1<q_2<  \cdots <  q_p } + \sum_{q_2<q_1 < \cdots < q_p} +
                \cdots + \sum_{q_p < \cdots < q_1} \right\}   f(q_1, \ldots ,q_p) ,  \\                                                  
   \end{eqnarray*}
where the sums inside the brackets correspond to the $p!$ orderings of $(q_1, \ldots , q_p)$.  Relabeling the summation indices, this is
  \begin{eqnarray*}
 \sum_{q_i \neq q_l} f(q_1, \ldots , q_p) &=& \sum_{ q_1< \cdots <  q_p }  \left\{ ~ f(q_1, \ldots ,q_p)+ \cdots + f(q_p, \dots , q_1) ~\right\} 
 =  p!  \sum_{ q_1< \cdots <  q_p }   f_S(q_1, \ldots ,q_p)                                                   
   \end{eqnarray*}   
where $f_S(q_1,\ldots , q_p)$ is the symmetrized function,
\begin{eqnarray*}
f_S(q_1, \ldots ,q_p) = \frac{1}{p!}~\sum_{\sigma \in {\cal S}_p} f(q_{\sigma_{(1)}}, \ldots , q_{\sigma_{(p)}}).
\end{eqnarray*}
Now, for $\sigma \in {\cal S}_p$,  define the $p \times p$ matrix
  \begin{eqnarray*}
({\cal O}(\sigma))_{ik} &=& \delta_{\sigma(i),k} .       
\end{eqnarray*}
Then
   \begin{eqnarray*}
  f_S ({\bf q}_{_Z}) &=&   \frac{1}{p!} ~ \sum_{\sigma \in S_p} ~\omega^{{\bf A} \cdot {\cal O}(\sigma)   {\bf q}_{_Z}}
                                    ~ \sum_{0 \leq \kappa_1, \ldots , \kappa_p \leq M_1 } 
                                  ~ \omega^{(\kappa_1 + \cdots +\kappa_{z_1})  q_1 + \cdots + (\kappa_{p-z_j} + \cdots + \kappa_p) q_j} 
                                  H[ M_1-( \kappa_1 + \cdots + \kappa_p )] \\        
\end{eqnarray*}
where ${ \bf A } = (A_1, A_2, \ldots, A_p)$.  Now defining $X^{(\sigma)}_s$ as the integer in $[0,N-1]$ such that 
\begin{eqnarray*}
X_s^{(\sigma)} &\equiv& - ( {\bf A} \cdot {\cal O}(\sigma) {\bf u}_s)  {\rm ~mod~} N\\
&\equiv& - \left(  A_{\sigma(d_s+1)} + \cdots + A_{\sigma(d_s+z_s)}  \right) {\rm ~mod~}N, ~~{\rm with ~} d_s = z_1 + \cdots + z_{s-1},
\end{eqnarray*}
we have, applying Lemma 2,
 \begin{eqnarray*}  
      \sum_{ q_i \neq  q_l } f(q_1, \ldots , q_{p}),    
                         &=&  \sum_{Z \in {\cal P}_p}  (-1)^{p+j } ~ \frac{ (p; k_1,k_2, \ldots , k_p)^*}{p!}
                    ~ \sum_{q_1, \ldots , q_j }  \sum_{\sigma \in S_p} ~ \omega^{ - (q_1 X_1^{(\sigma)} + \cdots + q_j X_j^{(\sigma)} )}\\
                    && \times   \sum_{0 \leq \kappa_1, \ldots , \kappa_p \leq M_1 } 
                                  ~ \omega^{(\kappa_1 + \cdots +\kappa_{z_1})  q_1 + \cdots + (\kappa_{p-z_j} + \cdots + \kappa_p) q_j} 
                                 H[~M_1-(  \kappa_1 + \cdots + \kappa_p)~]. \\ 
\end{eqnarray*}
The $p$ sums over the $\kappa$ variables can be reduced to $j$ sums over  $\alpha$ variables, defined as 
\begin{eqnarray*}
\alpha_1 &=& \kappa_{_1} + \cdots + \kappa_{_{z_1}},\\
& \vdots & \\
\alpha_k &=& \kappa_{_{z_1 + \cdots + z_{k-1} +1}} + \cdots + \kappa_{_{ z_1 + \cdots + z_{k}   }},   \\
& \vdots &\\
 \alpha_j &=&  \kappa_{_{p-z_j }}+ \cdots + \kappa_{_p;}\\
  \alpha_1 &+& \cdots + \alpha_j = \kappa_1 + \cdots + \kappa_p.
 \end{eqnarray*}
Then 
\begin{eqnarray*}  
     &&   \sum_{0 \leq \kappa_1, \ldots , \kappa_p \leq M_1 }  
             \omega^{(\kappa_1 + \cdots +\kappa_{z_1})  q_1 + \cdots + (\kappa_{p-z_j} + \cdots + \kappa_p) q_j} 
             H[~M_1-( \kappa_1 + \cdots + \kappa_p)~]\\
             &&~~~~~~~~~~~~   = \sum_{0 \leq \alpha_1, \ldots , \alpha_j \leq M_1 } 
            ~ \omega^{\alpha_1  q_1 + \cdots + \alpha_j q_j}  H[~M_1 -(  \alpha_1 + \cdots + \alpha_j ) ~]
             \prod_{s=1}^j  \left(  \begin{array} {c} \alpha_s+z_s -1 \\~\\z_s-1 \end{array} \right) ,
\end{eqnarray*}
where the multiplicative factor $ \left(  \begin{array} {c} \alpha_s+z_s -1 \\~\\z_s-1 \end{array} \right)$ is the number of compositions of  
$\alpha_s$ consisting of $z_s$ or fewer parts.  Making this substitution, we can finally use eq.(7) to evaluate the sums over the $q$'s:
   \begin{eqnarray*}  
      \sum_{ q_i \neq  q_j } f(q_1, \ldots , q_{p}),    
               &=&  (-1)^{p }\sum_{\sigma \in S_p}\sum_{Z \in {\cal P}_p}   ~ \frac{ (p; k_1,k_2, \ldots , k_p)^*}{p!}
              ~ \sum_{0 \leq \alpha_1, \ldots , \alpha_j \leq M_1 } ~H[~M_1 - ( \alpha_1 + \cdots + \alpha_j)~]  \\
               && \times \prod_{s=1}^j  \left(  \begin{array} {c} \alpha_s+z_s -1 \\~\\z_s-1 \end{array} \right)
                  ~\prod_{s=1}^j  (~1 -N \delta_{\alpha_s, X_s^{(\sigma)} }~) .
 \end{eqnarray*}
We expand the product $ \prod_{s=1}^j  (~1 -N \delta_{\alpha_s, X_s^{(\sigma)} }~)$ in powers of $(-N)$ as 
  \begin{eqnarray*}
  \prod_{s=1}^j  (~1 -N \delta_{\alpha_s, X_s^{(\sigma)} }~)
      &=& ~ \sum_{ \lambda_1, \ldots , \lambda_j = 0,1}   (-N)^{ \lambda_1 + \cdots + \lambda_j}  ( \delta_{\alpha_1, X_1^{(\sigma)} }~)^{\lambda_1} 
      \cdots  ( \delta_{\alpha_j, X_j^{(\sigma)} }~)^{\lambda_j}  
  \end{eqnarray*}   
and get,  after rearranging the order of some of the sums,
  \begin{eqnarray*}  
      \sum_{ q_i \neq  q_j } f(q_1, \ldots , q_{p}),    
               &=& (-1)^p  \sum_{\sigma \in S_p}~ \sum_{Z \in {\cal P}_p}  ~ \frac{ (p; k_1,k_2, \ldots , k_p)^*}{p!}
                 \sum_{ \lambda_1, \ldots , \lambda_j = 0,1} (-N)^{ \lambda_1 + \cdots + \lambda_j}  \\
               && \times ~ \sum_{0 \leq \alpha_1, \ldots , \alpha_j \leq M_1 } ~ H[~M_1 - ( \alpha_1 + \cdots + \alpha_j)~] 
                \prod_{s=1}^j  \left(  \begin{array} {c} \alpha_s+z_s -1 \\~\\z_s-1 \end{array} \right)   
                \prod_{s=1}^j ( \delta_{\alpha_s, X_s^{(\sigma)} })^{\lambda_s} .
 \end{eqnarray*}
  
Now if, for fixed $(\lambda_1, \ldots , \lambda_j)$,  a particular $\lambda_t$ is equal to 1, then that term will contain a factor of  $\delta_{\alpha_t, X_t^{(\sigma)}}$, and we can replace the $\alpha_t$ in the binomial coefficient by $X_t^{(\sigma)}$. 
 The above equation can therefore be written as 
   \begin{eqnarray*}  
      \sum_{ q_i \neq  q_j } f(q_1, \ldots , q_{p}),    
               &=&   (-1)^{p \sum_{\sigma \in S_p}  ~ }\sum_{Z \in {\cal P}_p} ~ \frac{ (p; k_1,k_2, \ldots , k_p)^*}{p!} 
               ~ \sum_{ \lambda_1, \ldots , \lambda_j = 0,1} (-N)^{ \lambda_1 + \cdots + \lambda_j}\\
           && \times   \sum_{0 \leq \alpha_1, \ldots , \alpha_j \leq M_1 } ~ H[M_1 - ( \alpha_1 + \cdots + \alpha_j)]  \\
             && \times  \prod_{s=1}^j \left[  \left(  \begin{array} {c} \alpha_s+z_s -1 \\~\\z_s-1 \end{array} \right)^{1-\lambda_s}
                \left(  \begin{array} {c} X_s^{(\sigma)} +z_s -1 \\~\\z_s-1 \end{array} \right)^{\lambda_s}   (\delta_{\alpha_s, X_s^{(\sigma)} })^{\lambda_s} \right] .
 \end{eqnarray*}
   
Again take $(\lambda_1, \ldots , \lambda_j)$ as fixed and assume that a particular $\lambda_t$ is equal to  1.  The summand in the sum over the corresponding $\alpha_t$ is then  $  H[M - ( \alpha_1 + \cdots + \alpha_j)]  \delta_{\alpha_t, X_t^{(\sigma)} }$ times factors independent of $\alpha_t$.  We can thus replace $\alpha_t$ in the step function by $X_t^{(\sigma)}$.  The sum over $\alpha_t$ is from 0 to $M_1$, and the sum over $\delta_{\alpha_t, X_{z_t}}$ will be zero if $X_t^{(\sigma)}  > M_1$.   But, with the above replacement, the step function is also zero in this case.  So the ``zero condition'' is enforced by $H$ and we can simply make the replacement:  $\sum_{\alpha_t}~ H[M - ( \alpha_1 + \cdots + \alpha_j)]   \delta_{\alpha_t, X_t^{(\sigma)}}  \rightarrow \left.   H[M - ( \alpha_1 + \cdots + \alpha_j)]  \right|_{\alpha_t = X_t^{(\sigma)} }$.\\
 
After performing the summations over all the $\alpha_s$'s for which $\lambda_s =1$,  we are left with the sums over the $\alpha$'s for which  $\lambda_s=0$.  For $\lambda_t=0$,

 \begin{eqnarray*}
 \sum_{\alpha_t=0}^{M_1} \left( \begin{array} {c} \alpha_t + z_t-1\\ ~ \\ z_t-1 \end{array} \right)
               H \left[  M_1-\alpha_t   - {\sum_{q=1}^j}' ~ (1-\lambda_q) \alpha_q   - \sum_{q=1}^j \lambda_q X_q^{(\sigma)}  \right]
               &=&  \sum_{\alpha_t=0}^{M'_1} \left( \begin{array} {c} \alpha_t + z_t-1\\ ~ \\ z_t-1 \end{array} \right)\\
                &=&   \left( \begin{array} {c} M'_1  + (1-\lambda_t) z_t\\ ~ \\ (1-\lambda_t) z_t \end{array} \right)H[M'_1],
 \end{eqnarray*}
 where the prime on the first sum over $q$ means that the sum is over all $q \neq t$,  and where
 \begin{eqnarray*}
 M'_1 =M_1 - {\sum_{q=1}^j}' ~ (1-\lambda_q) \alpha_q   - \sum_{q=1}^j \lambda_q X_q^{(\sigma)}.
 \end{eqnarray*}
Summing again for $\lambda_s =0$, with $M''_1 = M'_1 -\alpha_s$, we have
  \begin{eqnarray*}
&& \sum_{\alpha_s=0}^{M_1} \left(\begin{array} {c} \alpha_s + z_s-1\\ ~ \\ z_s-1 \end{array} \right)
                \left( \begin{array} {c} M_1''-\alpha_s +(1-\lambda_t)z_t \\~\\ (1-\lambda_t)z_t \end{array} \right)   H[ M''_1 -\alpha_s] \\
            && ~~~~~~~~~~ =  \sum_{\alpha_s=0}^{M''_1} \left(\begin{array} {c} \alpha_s + z_s-1\\ ~ \\ z_s-1 \end{array} \right)
                \left( \begin{array} {c} M_1''-\alpha_s +(1-\lambda_t)z_t \\~\\ (1-\lambda_t)z_t \end{array} \right) \\
               &&~~~~~~~~~~ =  \left( \begin{array} {c} M_1'' +(1-\lambda_t)z_t+(1-\lambda_s)z_s \\
                          ~\\ (1-\lambda_t)z_t +(1-\lambda_s)z_s \end{array} \right)   H[M''_1 ] ,
  \end{eqnarray*}
using the identity
   \begin{eqnarray*}
\sum_{m=0}^n \left( \begin{array} {c} m\\ ~ \\ l \end{array} \right) \left( \begin{array} {c} n-m \\~\\ k-l \end{array} \right) 
= \left(  \begin{array} {c} n+1 \\ ~ \\ k+1 \end{array} \right),~~ 0 \leq l \leq k \leq n.
\end{eqnarray*}
Continuing, the effect of summing over these $ \alpha_t$'s is to remove a $(-\alpha)$ from the step function and the binomial coefficient, and to  add $(1-\lambda_t) z_t$ to both the upper and lower elements of the binomial coefficient.  So, after summing over all $\alpha$ variables, we are left with the factor
\begin{eqnarray*}  
    H \left[~M_1 - \sum_{t=1}^j \lambda_t X_t^{(\sigma)}~\right] 
     \left( \begin{array} {c}  M_1  + \sum_{t=1}^j (1-\lambda_t) z_t - \sum_{t=1}^j \lambda_t X_t^{(\sigma)} \\
              ~\\ \sum_{t=1}^j (1-\lambda_t) z_t \end{array}  \right)  ~ 
               \prod_{s=1}^j \left(  \begin{array} {c} X_s^{(\sigma)} +z_s -1 \\~\\z_s-1 \end{array} \right)^{\lambda_s} .                    
  \end{eqnarray*}
     Putting this into the equation above for the sum over $f$ and then that into (10b), and then making the replacement  $\sum(1-\lambda_t) z_t = p - \sum \lambda_t z_t$, we arrive at
               \begin{eqnarray*}
 C_{0 \cdots 0\,1 \cdots 1\, A_1 \cdots A_p A_{p+1}} 
                 &=&  \frac{  (-1)^{N -M_0 -1}  N}{M_2! \cdots M_{N-1}!  } ~
                   \sum_{\sigma \in S_p} ~   \sum_{ Z \in {\cal P}_p} ~ \frac{ (p; k_1,k_2, \ldots , k_p)^*}{p!}    
           \sum_{\lambda_1, \ldots , \lambda_j =0,1 }(-N)^{ \lambda_1 + \cdots + \lambda_j}  \\
            && \times  ~H \left[~M_1 - \sum_{t=1}^j \lambda_t X_t^{(\sigma)}~\right]  
             \left( \begin{array} {c} N-M_0-1 - \sum _{s=1}^j \lambda_s (X_s^{(\sigma)}  + z_s)  \\~\\ M_1-  \sum \lambda_s X_s^{(\sigma)}
             \end{array} \right)
              \prod_{s=1}^j    \left( \begin{array} {c} X_s^{(\sigma)} + z_s -1\\~\\ z_s-1 \end{array} \right)^{\lambda_s}.    
     \end{eqnarray*}
For the terms for which $\lambda_1 = \cdots = \lambda_j =0$, the sum over $Z \in {\cal P}_p$, and then trivially the sum over $ \sigma \in  S_p$,  can be done using the identity        
        \begin{eqnarray*}
        \sum_{ Z \in {\cal P}_p} ~  (p; k_1,k_2, \ldots , k_p)^* \equiv
       \sum_{j=0}^p ~  \sum_{ k_1, \ldots, k_p} ~\delta_{k_1 +2k_2 + \cdots + pk_p, p}~\delta_{k_1 +k_2 + \cdots + k_p, j}~ (p; k_1,k_2, \ldots , k_p)^* =p!;
        \end{eqnarray*}
(this is eq.(E4) in Appendix E).  Then
         \begin{eqnarray}
  C_{0 \cdots 0\,1 \cdots 1\, A_1 \cdots A_p A_{p+1}} 
                 &=&  \frac{  (-1)^{N -M_0 -1}  N}{M_2! \cdots M_{N-1}!  } \left\{ ~ \frac{(N-M_0-1)!}{M_1!}   
                 + \sum_{ Z \in {\cal P}_p} ~ \sum_{\sigma \in S_p}    \frac{ (p; k_1,k_2, \ldots , k_p)^*}{p!}  ~
                   \right. \nonumber \\
            && \times  ~  \sum_{\mu=1}^j (-N)^{\mu} \sum_{\lambda_1, \ldots , \lambda_j =0,1 }\delta_{\mu,  \lambda_1 + \cdots + \lambda_j}  
               ~ H \left[~M_1 - \sum_{t=1}^j \lambda_t X_t^{(\sigma)}~\right]      \\
            &&  \times \left. \left( \begin{array} {c} N-M_0-1 - \sum _{s=1}^j \lambda_s (X_s^{(\sigma)}  + z_s)  \\~\\
             M_1-  \sum \lambda_s X_s^{(\sigma)}
             \end{array} \right)
              \prod_{s=1}^j    \left( \begin{array} {c} X_s^{(\sigma)} + z_s -1\\~\\ z_s-1 \end{array} \right)^{\lambda_s} \right\} .  \nonumber 
       \end{eqnarray}      
 
We now note that in the sum over $ S_p$, the set of $X_s^{(\sigma)}$'s is invariant for all $\sigma$ that  map a partition onto itself,  $\sigma(\Theta)= \Theta$,  since in this case $\sigma$ either just reorders the elements in each $\theta \in \Theta$ or interchanges the order of the $\theta$'s, and neither of these operations change the set $\{ X_s^{(\sigma)}|_{s=1,\ldots , j}\}$.  The sum over $S_p$ and ${\cal P}_p$ can then be replaced by a sum over the set of partitions of the index set  $\{12 \cdots p\}$ by  multiplying each term by $z_1!z_2! \cdots z_j!~k_1!k_2! \cdots k_p!$: 
   \begin{eqnarray*}
   \sum_{Z \in {\cal P}_p} ~\sum_{\sigma \in S_p} &&\rightarrow \sum_{ \Theta \in {\cal P} \{12 \cdots p\}}; \\
  \frac{ (p; k_1,k_2, \ldots , k_p)^*}{p!} && \rightarrow   \frac{ (p; k_1,k_2, \ldots , k_p)^*}{(p; z_1, \ldots , z_j)} \times k_1! \cdots k_p!
  = \frac{z_1! \cdots z_j!}{1^{k_1}2^{k_2} \cdots p^{k_p}} =  (z_1-1)! \cdots (z_j-1)!.
                 \end{eqnarray*}
With this substitution, 
\begin{eqnarray*}
  X_s^{(\sigma)}\rightarrow  X(\theta_s) \equiv  - \left( \sum_{~~q= n_1^{(s)}}^{~n_{z_s}^{(s)}} ~ A_q  \right) {\rm ~mod~} N,~~{\rm for~}
   \theta_s = \{ ~n_1^{(s)}, \ldots , n_{z_s}^{(s)}~\},
\end{eqnarray*}
and we get
  \begin{eqnarray}
 C_{0\cdots 0 \, 1 \cdots 1\, A_1 \cdots A_p A_{p+1} } &=&   \frac{  (-1)^{N -M_0 -1} N}{ M_2! \cdots M_{N-1}! } 
          \left\{~ \frac{(N-M_0-1)!}{M_1!  }    +    \sum_{\Theta \in {\cal P} \{12 \cdots p  \} }  (z_1-1)! \cdots (z_j-1)!       \right. \nonumber \\
        && \times ~ \sum_{\mu=1}^j  (-N)^{\mu}   \sum_{\lambda_1, \ldots , \lambda_j =0,1 } \delta_{\mu, \lambda_1 + \cdots + \lambda_j } 
                     ~ H \left[~M_1 - \sum_{t=1}^j \lambda_t X(\theta_t)~\right]    \nonumber \\ 
                &&\left.     \times   \left( \begin{array} {c} N-M_0-1 - \sum_{t=1}^j \lambda_t (X(\theta_t)  + z_t ) \\
            ~\\   M_1 - \sum_{t=1}^j \lambda_t X(\theta_t)  \end{array} \right) 
                \prod_{s=1}^j    \left( \begin{array} {c} X(\theta_s) + z_s -1\\~\\ z_s-1 \end{array} \right)^{\lambda_s} \right\} . 
   \end{eqnarray}    
   
 Expression (10d) above is an alternative expression for the coefficient to that  stated in the theorem.   The final expression is derived from (10d) by the replacement of the sum over ${\cal P} \{12 \cdots p  \}$ by a sum over ${\cal P} \{A_1 \cdots A_p  \}$, (which is 
 a smaller set when the coefficients $A_1, \ldots , A_p$ are not all distinct), and then noting that if, for example, $A_1=A_2$, a partition of the form $(A_1 \cdots )(A_2 \cdots )\cdots (\cdots )$ will contribute the same to the sum as the partition obtained by interchanging $1$ and $2$.  So, replacing the sum over  ${\cal P} \{12 \cdots p  \}$ by a sum over the smaller set ${\cal P} \{A_1 \cdots A_p  \}$ requires that the sum be multiplied by a factor of 2, or, more generally, by the factor $M_{A_1}! \cdots M_{A_p}!$.  However, for partitions of the form $(A_1A_2 \cdots )\cdots (\cdots )$ or of the form $(A_1B\cdots C)(A_2B\cdots C) \cdots (\cdots )$, interchanging the indices $1$ and $2$ does not generate additional partitions if $A_1=A_2$.  To correct for this overcounting, the sum then needs to be divided by the factorials of the $\kappa_{\theta}$ and the $m_a^{(\theta)}$ multiplicities.    $\triangle$\\
 
\begin{corollary}
~~If $a_0 + a_1 + \cdots + a_{N-1}  \equiv 0 {\rm ~mod~} N  $, then $C_{[a]} =0 $ only if $N$ divides  $ \frac{(N-M_0-1)!}{M_1!} $
\end{corollary}
~\\
The proof follows by noting that every term inside the brackets in expression (10d) is an integer, and that all terms other than the first one are
proportional to $N$.

\subsection{Coefficients of the form $C_{0\cdots 0\, 1 \cdots 1\, a\cdots a b}$}
    
   In the special case where $A_1=A_2 = \cdots = A_p =a$  the expression for the coefficients simplifies considerably.  In this case, 
 in eq.(10c),  $X_s^{(\sigma)}$ $\rightarrow X_{z_s} \equiv ( - az_s)$ mod $N$, independent of the permutation  $\sigma$.  The sum over ${\cal S}_p$  can then be done and, after some manipulation, so can the sum over ${\cal P}_p$.\\
    
 Summing over the permutations we have, for $M_b=0,1$,
        \begin{eqnarray*}
  C_{0 \cdots 0\,1 \cdots 1\, a \cdots ab} 
                 &=&  \frac{  (-1)^{N -M_0 -1}  N}{M_a! M_b! } \left\{ ~ \frac{(N-M_0-1)!}{M_1!}   
                  +   \sum_{ Z \in {\cal P}_p} ~   (p; k_1,k_2, \ldots , k_p)^*\right. \\
                 && \times \sum_{\mu=1}^j (-N)^{\mu} ~ \sum_{\lambda_1, \ldots , \lambda_j =0,1 } \delta_{\mu,  \lambda_1 + \cdots + \lambda_j}  
                   ~  H \left[~M_1 - \sum_{t=1}^j  \lambda_t X_{z_t}~\right]  \\
            && \left. \times  \left( \begin{array} {c} N-M_0-1 - \sum _{t=1}^j  \lambda_t (X_{z_t} + z_t)  \\~\\ M_1-  \sum_{t=1}^j  \lambda_t X_{z_t}
             \end{array} \right) 
              \prod_{s=1}^j    \left( \begin{array} {c} X_{z_s}+ z_s -1\\~\\ z_s-1 \end{array} \right)^{\lambda_s} \right\} .             
       \end{eqnarray*}
 
We now replace the sum over $(\lambda_1, \ldots, \lambda_j)$ by a sum over $(\beta_1, \ldots, \beta_p)$, which we define as follows:  Let the part  $z_n$ have multiplicity $k ~(= k_{z_n})$.  Then, (after a possible relabeling of the parts), $z_n =z_{n+1} = \cdots = z_{n+k-1}$.  For fixed $(\lambda_1, \ldots , \lambda_j)$, let $\beta$ equal the number of non-zero $\lambda$'s in the subset $(\lambda_n, \ldots , \lambda_{n+k-1});  ~{\rm i.e.,~} \beta = \sum_{t=n}^{n+k-1} \lambda_t$.      $\beta_{z_n}$ is thus defined for each unique part $z_n$.  This definition can be extended to all integers $q \in [1, \ldots , p]$ by setting $\beta_q = 0$ if $k_q=0$. \\ 

  For $q \in [1, \ldots, p]$,  $q$ and $X_q$ will each appear $\beta_q$ times in the sums over t, and the binomial coefficient $ \left( \begin{array} {c} X_q + q -1 \\ q-1 \end{array} \right)$ will  occur $\beta_q$ times in the product over $s$.  For a given $(\beta_1, \ldots , \beta_p)$, there are
   $ \left( \begin{array} {c} k_1 \\ \beta_1 \end{array} \right) \cdots  \left( \begin{array} {c} k_p \\ \beta_p\end{array} \right)$ different $(\lambda_1, \ldots , \lambda_j)$'s that correspond to it.  Replacing the sums over  $(\lambda_1, \ldots, \lambda_j)$ in the above expression by  sums over $(\beta_1, \ldots, \beta_p)$, 
    
\begin{eqnarray*}
\sum _{t=1}^j \lambda_t  z_t ,~   \sum _{t=1}^j \lambda_t X_{z_t}~  \rightarrow \sum _{t=1}^p t \beta_t,~  \sum _{t=1}^p \beta_t X_t ,
\end{eqnarray*}
and including this factor gives  us
       \begin{eqnarray*}
  C_{0 \cdots 0\,1 \cdots 1\, a \cdots ab} 
                 &=&   (-1)^{N -M_0 -1} ~ \frac{ N}{M_a! M_b! } \left\{ ~ \frac{(N-M_0-1)!}{M_1!}    +   \sum_{ j=0}^p ~\sum_{k_1, \ldots , k_p}  
                  (p; k_1,k_2, \ldots , k_p)^*\right. \\
                 && \times \sum_{\mu=1}^j (-N)^{\mu} ~ \sum_{\beta_1, \ldots , \beta_p  } \delta_{\mu,  \beta_1 + \cdots + \beta_p} 
                  \left( \begin{array} {c} k_1 \\ ~\\ \beta_1 \end{array} \right) \cdots  \left( \begin{array} {c} k_p \\~ \\ \beta_p\end{array} \right) 
                  ~  H \left[~M_1 - \sum_{t=1}^p \beta_t X_{t}~\right]  \\
            && \left. \times  \left( \begin{array} {c} N-M_0-1 - \sum _{t=1}^p  \beta_t (X_{t} + t)  \\~\\ M_1-  \sum_{t=1}^p  \beta_t X_{t}
             \end{array} \right)  \prod_{s=1}^p   \left( \begin{array} {c} X_{s}+ s -1\\~\\ s-1 \end{array} \right)^{\beta_s} \right\} .
       \end{eqnarray*}
 We now interchange the sums over the $k$'s and the $\beta$'s and perform the sum over the $k$'s by using the lemma below:
\begin{lemma} 
For $ p \geq   \beta_1 +2 \beta_2 + \cdots + p \beta_p $,
       \begin{eqnarray*}
 \sum_{j=0}^p~\sum_{k_1, \ldots , k_p} \left( \begin{array} {c} k_1 \\~\\ \beta_1 \end{array} \right) \cdots 
          \left( \begin{array} {c} k_p \\~\\ \beta_p \end{array} \right)    ~(p;k_1,k_2, \ldots , k_p)^*     = 
           \frac{p!}{1^{\beta_1} \beta_1!~\, 2^{\beta_2} \beta_2! \, \cdots \, p^{\beta_p}\beta_p!} .
           \end{eqnarray*}
\end{lemma}
(For the proof, see Appendix E.)  Then our result is
         \begin{eqnarray}
 \left. \begin{array} {c}  C_{0 \cdots 0\,1 \cdots 1\, a \cdots ab} \\~\\  C_{0 \cdots 0\,1 \cdots 1\,ba \cdots a} \end{array} \right\} 
                 &=&  (-1)^{N-M_0 -1}~ \frac{    N   }{M_a+M_b} \left( \begin{array} {c} M_a+M_b \\~\\ M_a \end{array} \right)
                  \left\{    \left( \begin{array} {c}N-M_0 -1   \\~ \\ M_1  \end{array} \right)  \right. 
                 \nonumber \\
                &&~~~~~~~~  + ~ \sum_{\mu=1}^{p} (-N)^{\mu}  \sum_{0 \leq \beta_1, \ldots , \beta_{p} \leq \mu }
                   ~ \delta_{\mu,  \beta_1 + \cdots + \beta_{p}}   H \left[~M_1 - \sum_{t=1}^p \beta_t X_{t}~\right]  \\
            &&~~~~~~ \left. \times  \left( \begin{array} {c} N-M_0-1 - \sum _{t=1}^p  \beta_t (X_{t} + t)  \\~\\ M_1-  \sum_{t=1}^p  \beta_t X_{t}
             \end{array} \right)  \prod_{s=1}^{p} \frac{1}{s^{\beta_s}\beta_s!} 
             \left( \begin{array} {c} X_{s}+ s -1\\~\\ s-1 \end{array} \right)^{\beta_s} \right\} ;  \nonumber\\
              X_s &\equiv & - sa {\rm ~mod~} N ~{~\rm such ~that~~} 0 \leq X_s \leq N-1 ;\nonumber \\
              p&=&N-M_0-M_1-1 .\nonumber
       \end{eqnarray}
       \end{subequations}
   for $M_b=0$ or 1.\\
      
As indicated, this expression is valid for $b<a$ as well as for $b \geq a$.  The reason is that, in the proof of Theorem 3, in summing over the equivalence classes we have the freedom to ``hide'' any particular index in the $\sigma_0$ position,  where it's value has no effect,  (other than on the requirement that it and the other indices satisfy condition (8a)).   We can therefore sum over equivalence classes of the form $[b, \sigma_1, \ldots , \sigma_{N-1}]$  for $b$ either less than or greater than $a$.

 \subsection{Coefficients of the form $C_{0\cdots 0\, 1 \cdots 1\, A_1A_2A_3}$}

 As previously noted, the two coefficients $C_{001335}$ and $ C_{0000111368}$, corresponding to $N=6$ and 10 respectively,  are zero despite satisfying condition (8a).   Generalizing, we consider coefficients of the form $C_{0\cdots 0 \, 1 \cdots 1\, A_1 A_2 A_3 }$, with 
 $M_0,M_1\geq 1$ and $2 \leq A_1 \leq A_2 \leq A_3  \leq N-1$; i.e., coefficients  with non-zero $M_0$ and $M_1$ that contain 3 and only 3 indices greater than 1.    We are interested in the question of when coefficients of this form equal zero even if they satisfy conditions (8).   From (8b,c) we have 
 \begin{eqnarray*}
  M_1 +A_1+A_2+A_3 &\equiv& 0 {\rm ~mod~}N, \\
  M_0 +M_1+3 &=& N,
 \end{eqnarray*}
 and from eq.(10d), 
   \begin{eqnarray}      
 C_{0\cdots 0 \, 1 \cdots 1\, A_1 A_2 A_3 } & \propto &   
        \left\{~\frac{}{} (M_1+2)(M_1+1)  -N  H \left[~M_1 -  X(A_1A_2)~\right]  \left( X(A_1A_2) + 1 \right) \right.  \\
  && \left. -N  \left( ~ H \left[~M_1 -  X(A_1) ~\right]   \left( \begin{array} {c} M_1+1 -  X(A_1)  \end{array} \right) +(A_1 \rightarrow A_2)~\right)
            +N^2 H[~M_1 -X(A_1)-X(A_2)~] ~ \frac{}{} \right\} \nonumber.
  \end{eqnarray}   
  For $C_{0\cdots 0 \, 1 \cdots 1\, A_1 A_2 A_3 }$ to be zero, $N$ must divide $(M_1+2)(M_1+1)$ for $M_1 \leq N-4$.  If $N$ divides $(M_1+2)(M_1+1)$, then it also divides $(M_0+2)(M_0+1)$:
 \begin{eqnarray*}
\frac{(M_1+2)(M_1+1)}{N} = \frac{(N-M_0-1)(N-M_0-2)}{N} = N-2M_0-3 +\frac{(M_0+2)(M_0+1)}{N}.
\end{eqnarray*}
  The $X$ functions in eq.(11) are:
  \begin{eqnarray*}
     X(A_{1,2}) &=& N-A_{1,2},\\
     X(A_1A_2) &\equiv& -(A_1+A_2) \equiv M_1+A_3 = \left\{ \begin{array} {l} M_1+A_3 ~~~~~~~~~{\rm if} ~A_3 < N-M_1; \\ ~\\
            M_1+A_3 -N ~\,~{\rm if} ~A_3 \geq N-M_1. \end{array} \right.     
   \end{eqnarray*}    
  If $A_3$ (and therefore $A_1$ and $A_2$) is less than  $N-M_1$, all of the Heaviside step functions in eq.(11) are zero and the coefficient is 
  therefore non-zero.  So we will assume that $A_3 \geq N-M_1$.  That leaves only 4 possibilities: 
     \begin{eqnarray*} 
 (1)&&     N-M_1 \leq A_3 ~{\rm and~}  A_1 \leq A_2 < N -M_1;\\  
 (2)&&     N-M_1 \leq A_2\leq A_3 ,~  A_1 < N -M_1, ~ {\rm and~}  A_1+A_2 <2N-M_1;\\        
 (3)&&     N-M_1 \leq A_1 \leq A_2  \leq A_3~{\rm and~} A_1+A_2 <2N-M_1;\\           
 (4)&&    N-M_1 \leq A_1 \leq A_2  \leq A_3 ~{\rm and~}  A_1+A_2 \geq 2N-M_1.         
\end{eqnarray*}

Now let $\psi$ be the symmetry operator defined as multiplication of an index set $[a]$ by $(N-1)$ followed by addition of 1.  Under this operator,
\begin{eqnarray*}
M_0 \rightarrow M_0'= M_1,~M_1 \rightarrow M_1' = M_0,~A_1 \rightarrow A_1' = N-A_3+1,~~A_2 \rightarrow A_2' = N-A_2+1,
~A_ 3\rightarrow A_3' = N-A_1+1.
\end{eqnarray*}
 Thus, if $N-M_1 \leq A_1\leq A_2 \leq A_3$, then $A_1'\leq A_2' \leq A_3' < N-M_1'$.  From the discussion above  and from 
 Proposition 2, in this case,
 \begin{eqnarray*}
 C_{0 \cdots 0 1 \cdots 1 A_1A_2A_3} = (-1)^{N-1} C_{\psi*[0 \cdots 0 1 \cdots 1 A_1A_2A_3]} \neq 0,
 \end{eqnarray*}
and we need to consider only cases (1) and (2).  Setting the right side of eq.(11) to zero in these two cases gives us the corollary 
to Theorem 3: \\
  
 \begin{corollary}
 ~~For $N=M_0+M_1+3$ with $M_0, M_1 \geq 1$, the coefficient $C_{0 \cdots 0 \, 1\cdots 1\, A_1A_2A_3}$ satisfies conditions (8) but is zero  if $N$  divides $(M_1+2)(M_1+1)$  and either
\begin{eqnarray*}
&& A_1 \leq A_2 <  N-M_1, ~A_1+A_2 =N +1- \frac{(M_1+2)(M_1+1)}{N} , {\rm ~and~}
A_3= M_0+2+\frac{(M_1+2)(M_1+1)}{N}, \\
or &&\\
&&   N-M_1 \leq A_2\leq A_3 ,~A_2+A_3 = N+1+\frac{(M_0+2)(M_0+1)}{N} , {\rm ~and~} 
A_1 =  M_0 +2- \frac{(M_0+2)(M_0+1)}{N} .\\
\end{eqnarray*}
\end{corollary}
~\\
By Proposition 2, all coefficients that are related to these zero coefficients by an additive- and/or a multiplicative-symmetry operator also equal zero. \\

For integers $a$ and $q$, $N$ in the corollary has the form $N=q(q+1)/a,$ with values $ 6,10,12,14,15,18, \ldots, $ corresponding to sequence A139799 in OEIS \cite{oeis}. \\

 As an example, we take $N=6$.  In this case $N$ divides $(M_1+2)(M_1+1)= 6, 12, 30, 42, \ldots$, but only the first two integers satisfy the condition $N \geq M_1 +4 $.  From the corollary,
\begin{eqnarray*}
 && C_{001245} ,~  C_{011235} ~~(= -C_{\psi*[001245 ]}) ,\\
 &&C_{001335} ,~ C_{011244}~~ (=-C_{\psi*[001335]})
 \end{eqnarray*}
 are zero coefficients.  Including all the coefficients generated from these by the Proposition-2 symmetries, there are a total of 12 zero coefficients. \\
   
 As a second example, $N=10$.  $N$ divides $(M_1+2)(M_1+1)= 20,~30$, and the corollary gives us the zero coefficients
 \begin{eqnarray*}
 && C_{0000111377}, ~   C_{0001111448} ~~(=- C_{\psi*[0000111377]}),\\
 &&  C_{0000111368},~C_{0001111358  } ~~(=- C_{\psi*[ 0000111 368 ]}),\\
 && C_{0000111458}, ~C_{0001111367} ~~( = -C_{\psi*[0000111458  ]}).
 \end{eqnarray*}
 Including the Proposition-2-operator-generated coefficients, there are 120 zero coefficients of this form.

\section{Additional results:  Multiplet structure}

Having a formula to calculate any coefficient in the expansion (3a/b), it is of interest to determine how many calculations one needs to perform to solve for the complete determinant.  The number of solutions to eq.(8b,c) is given by the expression \cite{Brualdi, oeis},
 \begin{eqnarray}
F(N) = \frac{1}{2N}\sum_{d|N } ~\phi \left( \frac{N}{d} \right)  \left( \begin{array} {c}  2d \\ ~\\ d \end{array} \right)
    =  ({\rm A}003239)_{N+1}~{\rm  in~ OEIS}.
\end{eqnarray}
where $\phi$ is the Euler totient function \cite{Totient}.  Then, except for  coefficients that equal zero as a result of satisfying the criteria of Corollary 6 and any other ``accidental zeros'', eq. (12) gives the number of terms in the expansions (3a/b). \\

As discussed above following Proposition 2, coefficients can be grouped into additive multiplets, related by the additive symmetry, and into super-multiplets,  related by the additive symmetry and/or the multiplicative symmetry.   For example, for $N=5$, the coefficients can be
formed into the table: 

\begin{eqnarray*}
5:&& \{ C^*_{50000} \}_5  \\
311:&& \{ C^*_{31001}\}_5, ~\{ C^*_{30110}\}_5  \\
221:&&  \{C^*_{12002}\}_5, ~\{C^*_{10220}\}_5    \\
11111:&&    \{C^*_{11111}\}_{ 1} .
\end{eqnarray*}
The set of coefficients consist of six additive multiplets $\{C^*_{[M]} \}_n$, (where $n$ is the number of elements), and four super-multiplets.  (In this example there is only one super-multiplet for each of the allowed partitions of 5, but that is not generally the case; see the table for the $N=8$ coefficients in Appendix G.)\\

 The nonnegative integers less than $N$ form a group under addition modulo $N$, which we will denote as $({\bf Z}/N {\bf Z})^+$.  
 The group $({\bf Z}/N{\bf Z})^+$ is of order $N$, and is a cyclic group, so all its subgroups are cyclic.  The set of  integers less than and relatively prime to $N$ is the multiplicative group of integers modulo $N$, denoted as $( {\bf Z} /N{\bf Z})^{\times}$.  This group has order $\phi(N)$.  The direct product of these two groups,    
   \begin{eqnarray*}
  H_N \equiv  ({\bf Z}/N{\bf Z})^+ \times ({\bf Z}/N{\bf Z})^{\times},
   \end{eqnarray*}
is thus of order $N\phi(N)$.\\

 If the index set $[M]$ is such that $a * [M] =[M]$ for all elements $a$ in a group, we say that $[M]$ is invariant under the action of that group.  To find the numbers of additive multiplets and of super-multiplets, we will need to count the number of index sets that 
 are invariant under the various subgroups of $ ({\bf Z}/N{\bf Z})^+$ and of $H_N$, respectively.  \\

\subsection {The number of additive multiplets}

From eq.(3b) the determinants are expressible in the form
   \begin{eqnarray*}
    {\rm det ~} [  x_0, x_1, \ldots , x_{N-1}] &=& \sum_{0 \leq M_0 ,  \ldots , M_{N-1} \leq N} 
     C^*_{M_0 \cdots  M_{N-1}} x_0^{M_0 }\cdots x_{N-1}^{M_{N-1} } .
\end{eqnarray*}
For a given coefficient $C^*_{M_0 \cdots M_{N-1}}$, the elements of the $n$-element additive multiplet $\{C^*_{M_0M_1 \cdots M_{N-1}}\}_n$ are obtained by circularly permuting the index set $[M_0M_1 \ldots  M_{N-1} ]$:
     \begin{eqnarray*}
     \{ C^*_{M_0M_1 \cdots M_{N-1}} \}_n &=& \left\{~ \frac{}{} C^*_{M_0M_1\cdots M_{N-1}},~  C^*_{P[M_{0} M_1
     \cdots M_{N-1}]},~ C^*_{P^2[M_0M_1 \cdots M_{N-1}]} \ldots , C^*_{P^{n-1}[M_{0}M_{1}\cdots M_{N-1}]}~\right\}
     \end{eqnarray*}
 where $P$  is the circular permutation operator defined in eq.(9).   We define the symbol $[~ \{C^*_{M_0M_1 \cdots M_{N-1}}\}_n~ ]$  as equal to the value of the coefficient that labels the additive multiplet:
\begin{eqnarray}
[~ \{ C^*_{M_0M_1 \cdots M_{N-1}} \}_n ~] = C^*_{M_0M_1\cdots M_{N-1}}.
\end{eqnarray}
Then, as a sum over these multiplets, the above expansion can be written as 
  \begin{eqnarray}
    {\rm det ~} [  x_0, x_1, \ldots , x_{N-1}] &=& \sum_{{\rm Multiplets}}  [~\{  C^*_{M_0M_1\cdots  M_{N-1}}\}_n  ~]
  \left[~ x_0^{M_0}x_1^{M_1} \cdots x_{N-1}^{M_{N-1}}
    +(-1)^{N-1} x_1^{M_0}x_2^{M_1} \cdots x_0^{M_{N-1}}  \right.\nonumber  \\ 
&& ~~~~~~~~~~~~~~~~~~~~~~~~~~~~~~~~~~~~~ 
\left. + \cdots + (-1)^{(N-1)(n-1)} x_{n-1}^{M_0}x_{n}^{M_1} \cdots x_{n-2}^{M_{N-1}} ~ \right]    .
\end{eqnarray}

\begin{prop}
For a given $N$, the number $g_N(n)$ of n-element additive multiplets $\{C^*_{[M]}\}_n$ is
\begin{eqnarray*} 
         g_N(n) =  \frac{ 1  +(-1)^{N-n} }{4n^2} \sum_{d|n} (-1)^{n+d} \mu\left(  \frac{n}{d} \right) \left(  \begin{array} {c} 2d \\~\\ d \end{array} \right)
\end{eqnarray*}
\end{prop}

{\bf Proof:}\\

  Operating on an index set $[M]$, the additive group $({\bf Z}/N{\bf Z})^+$ generates a set with $N$ not-necessarily-distinct elements:
\begin{eqnarray*}
({\bf Z}/N{\bf Z})^+* [M] = \{ ~P^{-a}[M]~ |~0 \leq a \leq N-1 \}.
\end{eqnarray*}
                                       
   Let $\{ M\}_n = \left\{~ [M], ~P[M], P^2[M], \ldots , P^{n-1}[M]~\right\}$ be an $n$-element additive multiplet.  Then $n$ divides $N$ and, for $d=N/n$, each element in the multiplet is  invariant under the action of any element of the subgroup $S^+_d \subseteq  {\bf Z}/N {\bf Z}^+$.  One of the index sets in $\{M\}$, which,  without loss of generality we can take to be the ``labeling'' set $[M]$,  must then have the form
\begin{eqnarray*}
[M_0,M_1, \ldots , M_{N-1}] = [M_0, M_1, \ldots , M_{n-1},\, M_0, M_1, \ldots , M_{n-1}, \, \cdots , M_0, M_1, \cdots , M_{n-1}],
\end{eqnarray*}
in which the subset $\{M_0, M_1, \ldots , M_{n-1}\}$ is sequentially repeated $d$ times.  We have then:
\begin{subequations}
\begin{eqnarray*}
\sum_{q=0}^{N-1} M_q &=& (N/n) ~\sum_{q=0}^{n-1} M_q,\\
 \sum_{q=0}^{N-1} qM_q &=& \sum_{q=0}^{n-1} M_q [~ q +(q+n)+(q+2n) + \cdots + (q+n(~N/n-1)~)~] \\
 &=&  \sum_{q=0}^{n-1} M_q \left[  ~q \frac{N}{n} + n \frac{(N/n-1)(N/n)}{2} ~ \right] \nonumber \\
 &=& \frac {N}{n} \sum_{q=0}^{n-1} qM_q + \frac{N}{n}~\frac{N-n}{2} \sum_{q=0}^{n-1} M_q. \\
 \end{eqnarray*}
 
 Conditions (8b,c) then require that 
\begin{eqnarray}
&&\sum_{q=0}^{n-1} M_q = n,\\
&& {\rm and} \nonumber \\
&& \sum_{q=0}^{n-1} qM_q +  \frac{(N-n)n}{2}   \equiv 0 {\rm ~mod~}  n  .
 \end{eqnarray} 
 The second of these equations can only be satisfied if $N-n$ is even.  (If $N-n$ is odd, $n$ must then be even in order that the right side to be an integer.  But  then $N$ is odd, which contradicts the requirement that $n$ divide $N$.)  In this case then, 
 \begin{eqnarray}
 \sum_{q=0}^{n-1} qM_q \equiv 0 {\rm ~mod~}n,~~N-n {\rm ~even}.
 \end{eqnarray}
  \end{subequations}
From eq.(12), the number of solutions $[M_0, M_1, \ldots , M_{N-1}]$ that satisfy both eqs.(8a,b) and (15a,c) is therefore
\begin{eqnarray*}
F_N(n)  =  \frac{1+(-1)^{N-n} }{4n} ~  \sum_{d|n } \phi \left( \frac{n}{d} \right) \left( \begin{array} {c}  2d \\ ~\\ d \end{array} \right),
\end{eqnarray*}
and $n$-member multiplets $\{M\}_n$ and $\{C^*_{[M]}\}_n$ exist only if $n$ is a same-parity divisor of $N$.  \\

The number of solutions $F_N(n)$ is a sum over all coefficients that lie in $d$-element multiplets for all $d$ that divides  $n$ and is of the same parity as $n$ and $N$.  Let  $l(d)$ denote the total number of coefficients $C_{[M]}$ for which $d$ is the smallest integer such that $P^d [M] =[M]$.  Then the number of $n$-element additive multiplets, $g_N(n)$,  is $l(n)/n$, and  
\begin{eqnarray*}
F_N(n) = \sum _{d|n} \frac{1+(-1)^{n-d}}{2} ~l(d)  \equiv \sum _{d||n} l(d) ,
\end{eqnarray*}
where the notation $d||n$ means that the sum is over all same-parity divisors of $n$.   It is straightforward to show that 
\begin{eqnarray*} 
 l(n) = \sum_{d || n }~ F_N (d) \, \mu \left(  \frac{n}{d} \right) ,
\end{eqnarray*}
by considering separately the $n$ odd and the $n$ even cases and using M\"{o}bius inversion.  Then
\begin{eqnarray*} 
g_N(n) =\frac{1+(-1)^{N-n}}{4n^2}\sum_{d || n}~ \sum_{q | d }~ \frac{ n}{d} ~ \mu \left( \frac{n}{d} \right) 
         ~\phi \left( \frac{d}{q} \right)  \left( \begin{array} {c}  2q \\ ~\\ q \end{array} \right) .\\
    \end{eqnarray*}
    
It remains to prove the identity 
\begin{eqnarray} 
&& \sum_{d||n}  \sum_{q | d }~ \frac{ n}{d} ~   \mu \left( \frac{n}{d} \right) \phi \left( \frac{d}{q} \right)   f(q)
= \sum_{b|n} (-1)^{n+b} \mu \left( \frac{n}{b} \right) f(b)
\end{eqnarray}
where $f(q)$ is an arbitrary function of $q$.  In evaluating the right-hand side of this equation, we express $n$ as $n=2^k m,$ where $m$ is an odd integer and $k\geq 0$.  The divisors of $n$ are of the form $b= 2^l a$, where $a$ divides $m$ and is therefore odd, and $l \leq k$.   Then\\
\begin{eqnarray}
{\rm RHS} = \sum_{b|n} (-1)^{n+b} \mu \left( \frac{n}{b} \right) f(b)
        = \sum_{a|m}\mu\left( \frac{m}{a} \right)  \left\{ \begin{array} {cll} f(a), &~&  k=0;\\~
            \\  f(2a)  + f(a)  , &~& k=1;\\~
             \\  f(2^k a) - f(2^{k-1}a),&~&  k \geq 2.  \end{array} \right.
\end{eqnarray}

We now consider the expression on the left side of eq.(16), and first take $n$ to be odd, ($k$=0).  Then all of its divisors are also odd, and we can make the replacement $d||n \rightarrow d|n$.   For a given divisor $b$, we set $d=bc ,~ c \leq n/b$, and isolate the $q=b$ term in the sum over $q$:
\begin{eqnarray*} 
{\rm LHS_{{\it n}~odd}} = \sum_{d|n}  \sum_{q | d }~ \frac{ n}{d} ~  \mu \left( \frac{n}{d} \right) \phi \left( \frac{d}{q} \right) ~  f(q)
&=& \cdots + f(b)   \sum_{c | (n/b)} \frac{ n/b}{c} ~ \mu \left( \frac{n/b}{c} \right) \phi (c)+ \cdots  \\
&=& \cdots + f(b) \mu\left(  \frac{n}{b}  \right) + \cdots = \cdots +  \sum_{a=b|m=n}  \mu\left( \frac{m}{a} \right) f(a) + \cdots ,
\end{eqnarray*}
where we've used the identity
\begin{eqnarray*}
\sum_{c | m} \frac{ m}{c} ~ \mu \left( \frac{m}{c} \right)  \phi (c)  = \mu(m),
\end{eqnarray*}
which follows from M\"{o}bius inversion of 
\begin{eqnarray}
\phi(m)  = m \sum_{d|m} \frac{\mu(d)}{d}.
\end{eqnarray}
Eq.(16) is therefore valid for $n$ odd.\\

For even $n$, we again set $n=2^k m$, with $k\geq 1$ and $m$ odd.  Its same-parity divisors are $d=2^la$, with $a|m$ and $ 1 \leq l \leq k$.
Setting $q=2^s p$, with $0 \leq s \leq l$ and $p$ an odd integer that divides $a$, the left-hand side of (16) is
\begin{eqnarray*} 
 {\rm LHS_{ {\it n}~ even}} =  \sum_{a|m} \sum_{l=1}^k  \sum_{p | a } \sum_{s=0}^{l} ~ \frac{ 2^{k-l} m}{a} ~  \mu \left( \frac{2^{k-l}m}{a} \right)
 \phi \left( \frac{2^{l-s}a}{p} \right)   f(2^s p) .
 \end{eqnarray*}
 But since $ \mu (2^K m)  =0$ for any $K >1$,  the only non-zero terms in the sum over $d||n$ are those for which $l=k,~k-1 >0$;
 \begin{eqnarray*}
{\rm LHS_{ {\it n}~ even}} &=&  \sum_{a|m}  \sum_{p | a }
    \left[  \sum_{s=0}^k ~ \frac{ m}{a} ~  \mu \left( \frac{m}{a} \right) \phi \left( \frac{2^{k-s}a}{p} \right) ~  f(2^s p)
  +(1-\delta_{1,k})  \sum_{s=0}^{k-1} ~ \frac{ 2m}{a} ~  \mu \left( \frac{2m}{a} \right) \phi \left( \frac{2^{k-1-s}a}{p} \right) ~  f(2^s p)\right]. \end{eqnarray*}
The two identities
  \begin{eqnarray*}
  \mu(2x)=-\mu(x),~~ \phi(2^K x) = 2^{K-1} \phi(x)~~{\rm for~} K>0,
  \end{eqnarray*}
where $x$ is in both cases an odd  positive integer, can be used  to take the M\"{o}bius and Euler functions outside the parentheses:
  \begin{eqnarray*}
    {\rm LHS_{ {\it n}~even }} = \sum_{a|m}   \frac{ m}{a} ~  \mu \left( \frac{m}{a} \right)  \sum_{p | a }~\phi\left( \frac{a}{p} \right)
  \left[ f(2^k p) +  \sum_{s=0}^{k-1} ~  2^{k-s-1}   f(2^s p)  - (1-\delta_{1,k})\left(2 f(2^{k-1}p)
       + \sum_{s=0}^{k-2} ~ 2^{k-s-1}  f(2^s p) \right) \right] .
       \end{eqnarray*}
 The two sums from $s=0$ to $k-2$ cancel and we have 
       \begin{eqnarray*}
   {\rm LHS_{ {\it n}~even }} = \sum_{a|m}   \frac{ m}{a} ~  \mu \left( \frac{m}{a} \right)  \sum_{p | a }~\phi\left( \frac{a}{p} \right)
           \left[~ f(2^k p)  + (2\delta_{1,k} -1) f(2^{k-1}p)~  \right] .
\end{eqnarray*}
Now applying the ``odd-$n$'' result,
\begin{eqnarray*} 
&& \sum_{d|n}  \sum_{q | d }~ \frac{ n}{d} ~      \mu \left( \frac{n}{d} \right) \phi \left( \frac{d}{q} \right) ~  f(q)
= \sum_{b|n}\mu \left( \frac{n}{b} \right) f(b),
\end{eqnarray*}
to the last expression, it becomes
\begin{eqnarray*} 
 {\rm LHS_{ {\it n}~even }} =  \sum_{a|m}  \mu \left( \frac{m}{a} \right)[f(2^ka) +(2 \delta_{1,k} -1) f(2^{k-1}a)~],
\end{eqnarray*}
in agreement with the 2nd and 3rd lines in eq.(17).   $\triangle$\\
~\\

Note that the expression
\begin{eqnarray*} 
  \frac{ 1 }{2n^2} \sum_{d|n} (-1)^{n+d}~ \mu \left( \frac{n}{d} \right)  \left(  \begin{array} {c} 2d \\~\\ d \end{array} \right)
\end{eqnarray*}
corresponds to the sequence A131868 in OEIS.  So $g_N(n) = ({\rm A}131868)_{n||N}$, and the total number of additive multiplets,
$n_{AM}$, is 
\begin{eqnarray*} 
    \sum_{n||N} g_N(n) = \sum_{n||N} ({\rm A}131868)_n .
 \end{eqnarray*}
\\

 \subsection{The number of super-multiplets}
   
   To find the number of super-multiplets, we consider the direct product:
\begin{eqnarray*}
 H_N = ( {\bf Z}/N {\bf Z})^{+} \times ( {\bf Z}/N {\bf Z})^{\times}. 
\end{eqnarray*}
We will denote an element of $H_N$ either with boldface, ${\bf a}$, or as the pair $(a,b)$.  Group multiplication of two elements of $H_N$ is defined by
\begin{eqnarray*}
 (a,b) \cdot (c,d) = (a+c,bd);~~a,c \in ({\bf Z}/N{\bf Z})^+,~~ b,d \in ({\bf Z}/N{\bf Z})^{\times}.
 \end{eqnarray*} 
In this formalism, the symmetry operations in Proposition 2 are expressed as the action of an element $(a,b)$ on an index set $[M]$, which we will denote as $(a,b)*[M]$: 
\begin{eqnarray*}
 (a,b) *[M]  = [M_{N-a},M_{b^{-1}-a},M_{2b^{-1}-a}, \ldots , M_{(N-1)b^{-1}-a}].
 \end{eqnarray*} 
 In particular, if $b=1$,
 \begin{eqnarray*}
 (a,1) *[M]  = P^{-a}[M].
 \end{eqnarray*}
 
  A direct-product group $G_1 \times G_1$ is cyclic iff $G_1$ and $G_2$ are both cyclic and have orders that are relatively prime to
each other.  By Gauss' Theorem,  $( {\bf Z} /N{\bf Z})^{\times}$ is a cyclic group iff $N$ is of the form $ 2,~4,~p^n$ or $2p^n$ where $p$ is an odd prime and $n \geq 1$ \cite{MMG}.   Consequently $H_N$ is cyclic only for the case where $N$ is a prime number.\\
 
 We will denote an order-$d$ subgroup of $({\bf Z}/N{\bf Z})^+$ as $S_d^+$;  since $({\bf Z}/N{\bf Z})^+$ is a cyclic group, all its subgroups are cyclic, and each $d$ corresponds to only one subgroup.   Likewise, we will denote order-$d$ subgroups of $({\bf Z}/N{\bf Z})^{\times}$ as $S_{d;\alpha}^{\times}$, and those of $H_N$ as $S_{d;\alpha}$, where the $\alpha$ indices are used to distinguish same-$d$ subgroups from one another.  \\
 
   Operating on an index set $[M]$, $H_N$ generates a set with $N \phi(N) $  elements:
\begin{eqnarray*}
H_N* [M] = \{ ~(a,b)*[M]~ |~0 \leq a \leq N-1;~1 \leq b \leq N-1,~{\rm gcd}(b,N)=1 ~  \}.
\end{eqnarray*}
  If $[M]$ is a member of an $m$-element super-multiplet, then $[M]$ is invariant under the action of an order-$(N \phi(N)/m)$ subgroup of $H_N$.  If $[M]$ is invariant under the action of all elements of a subgroup $S_{d;\alpha}$, then ${\bf c}_k *[M] = {\bf c}_j *[M]$ for all elements ${\bf c}_a, ~{\bf c}_b$  of a coset $\{{ \bf c}_k\}$ of $S_{d;\alpha}$.  If the elements of the set $\{~[M], ~ { {\bf c}_1} * [M],~\ldots ,{\bf  c}_{m-1} * [M]~\}$ are all distinct, where $m=N\phi(N)/d$ and  ${\bf c}_k$ is any element of the coset $\{ {\bf c}_k\}$,   then this set is an $m$-element super-multiplet.\\
                                     
  For $S_{d;\alpha} \subseteq H_N$ we define $K(S_{d;\alpha})$ to be the number of index sets $[M]$ that are invariant under $S_{d;\alpha}$ and in addition satisfy eqs.(8b,c):
  \begin{eqnarray*}
  K(S_{d;\alpha}) = \#\left\{~[M] ~|~ S_{d;\alpha}*[M]=[M],~S_{d;\alpha} \subseteq H_N; ~\sum_k kM_k \equiv 0 {\rm ~mod~} N;  ~   \sum_k M_k=N ~\right\} .
  \end{eqnarray*}
  If $[M]$ is invariant under the action of the subgroup $S_{d;\alpha}$, it is also invariant under all subgroups of $S_{d;\alpha}$.  It may also be invariant under larger subgroups of $H_N$ that contain $S_{d;\alpha}$; i.e., it may have a higher degree of symmetry than that required to be invariant under just $S_{d;\alpha}$.  $K(S_{d;\alpha})$ therefore counts these higher-symmetry sets as well.  We define  $L(S_{q;\beta})$ as the number of index sets satisfying (8b,c)  for which $S_{q;\beta}$ is the largest subgroup that they are invariant under.   Then
   \begin{eqnarray*}
  K(S_{d;\alpha}) =  \sum_{S_{q;\beta} \supseteq S_{d;\alpha}}~L(S_{q;\beta}).
    \end{eqnarray*} 
  Let
  \begin{eqnarray*}
  n= \frac{N\phi(N)}{q}.
 \end{eqnarray*}
 Then  $L(S_q)/n$ is the number of $n$-element super-multiplets, and the total number, $n_{SM}$, is 
  \begin{eqnarray}
  n_{SM}= \sum_{S_q \subseteq H_N}  \frac{qL(S_q) } {N\phi(N)} .
      \end{eqnarray}
 We can formally solve this equation by expressing the previous equation in matrix form:
 \begin{eqnarray*}
 {\bf K } = {\bf A}  {\bf L},~~{\bf L } = {\bf A}^{-1} {\bf K} , \\
 n_{SM} = \frac{1}{N\phi(N)}~{\bf q } \cdot {\bf L}
 \end{eqnarray*}
 where the elements of the matrix ${\bf A}$ are either 1 or 0 depending on whether or not a particular subgroup $S_{q;\beta} $ contains 
 the subgroup $S_{d; \alpha}$ or not, and the rest of the notation should be clear from the context.\\
 
 Because of the more complicated structure of the multiplicative group, calculating the number of super-multiplets is more involved than for the additive multiplets.  Expression (19) however simplifies in the cases where $N=p$ or $2p$, where $p$ is an odd prime number.  In 
 Appendix A, we calculate the number of super-multiplets for these two cases.  Our results are:
 
\begin{subequations}
 \begin{eqnarray}
    n_{SM}(N=p) &=& ~\frac{1}{p(p-1)} \left\{   \frac{ 1 }{2p} \left(  \begin{array} {c} 2p \\~\\ p \end{array} \right) + \frac{p^2-1}{p} 
    +  p  \sum_{m|(p-1)}^{m< p-1}\phi \left(   \frac{p-1}{m}  \right) \left( \begin{array} {c} 2m \\~\\ m \end{array} \right) \right\}  ; \\
  n_{SM}(N=2p)  &=& \frac{1}{2p(p-1)} \left\{    \frac{1}{4p}     \left( \begin{array} {c}  4p \\ ~\\ 2p \end{array} \right) 
      + \frac{4p^2+1}{4p} \left( \begin{array} {c} 2p \\ ~\\ p \end{array} \right) + \frac{2(p^2-1)}{p}  \right. \\
    && ~~ \left.   +p \sum_{m|(p-1)}^{m<(p-1)/2}   \phi\left( \frac{p-1}{m} \right) \left[\left(  \frac{p+ 4m +1 }{2m+1 }
    - \frac{1-(-1)^{\frac{p-1}{m} }}{4} \right)  
         \left( \begin{array}{c} 4m  \\~\\  2m \end{array} \right)  
 + \frac{1 - (-1)^{\frac{p-1}{m}} }{4} \left( \begin{array} {c} 2m \\~\\ m \end{array} \right)  \right] \right\} . \nonumber
    \end{eqnarray}  
\end{subequations}
    
Using eq.(12) and another series expression in OEIS,
      \begin{eqnarray*}
 \frac{1}{2}  \sum_{m|n}  ~ \phi \left( \frac{n}{m} \right) \left( \begin{array} {c} 2m \\~\\  m \end{array} \right) = ({\rm A}081875)_n,
     \end{eqnarray*}  
 expression (20a) becomes  
 \begin{eqnarray*}
    n_{SM} (N=p)  = \frac{1}{p(p-1)} \left\{  ({\rm A}003239)_{p+1}  +  2p ({\rm A}081875)_{p-1}  +p-1 
       - p \left( \begin{array} {c} 2p-2\\~\\p-1 \end{array} \right) \right\} .
    \end{eqnarray*}

 \section{Conclusion}
 
  The expression for $C_{a_0 \cdots a_{N-1}}$ in Theorem 3 constitutes in principle a complete solution to the problem of calculating the determinant of an arbitrary circulant matrix.  This expression has some similarities to that of Wyn-jones, as well as some significant differences, but the proofs of the two expressions are along completely different lines.  \\
  
  However, Wyn-jones's method does require that one find all of the null subsets/multsets contained in a generally much larger set than the set of partitions employed in Theorem 3.    And since, in applying Theorem 3, one deals only with partitions of sets that do not include any $a=1$ indices, it may be that it is more efficient than the method of Wyn-jones for coefficients with $M_1 >>1$.\\
  
  In many cases, the number of terms in the expression for a particular coefficient can be decreased by using one or both of the symmetry operations in Proposition 2 to equate the coefficient (up to a sign) to another coefficient.  This can often be done in one of two ways:  By either making the new multiplicities $M'_0$ and/or $M'_1$ as large as possible, thereby making the set of partitions to be summed over smaller, or by setting $M'_1$ to 0, making the constraint  $M_1 - \sum \lambda_s X(\theta_s) \geq 0$ more restrictive.\\
 
    To illustrate this second method, let $C_{a_0\cdots a_{N-1}}$ be such that the integer $n$ is not in the index set.  Then if we add $(N+1-n)$ to $ [a_0,\ldots, a_{N-1}]$, the resulting index set will not contain 1, so  $M'_1 =0$.  We can always do this except for the index set $[0,1,2, \cdots, N-1]$. \\ 
    
    Consider the previous example, $C_{0011113788}$.  The missing indices are $n=2,4,5,6,9$, so $N+1-n=9,7,6,5,2$,  and 
    \begin{eqnarray*}
 C_{0011113788} = -C_{0000267799} = -C_{0455778888} = +C_{3446677779}  = -C_{2335566668}= +C_{0022333359}.
 \end{eqnarray*}
 Of these, $C_{0000267799}$ is the easiest to compute as there are only two partitions,  $ (677)( 29)$  and $(677)(2)(9),$
 that contribute to the sum. \\
 
    This method of setting $M_1$ to zero has some similarity to Wyn-jones' method, in its reliance on the null subsets/multisets.   The contributing partitions of the set $\{267799\}$ would be in that case the trivial partition, $(267799)$, and the partition $(677)(299)$.\\

 For $C_{0011113788}$, the number of non-zero terms in the sum can actually be reduced further by  instead multiplying the index set by $(N-1)$ and then adding 1.  This has the effect of interchanging $M_0$ and $M_1$:
 \begin{eqnarray*}
[0011113788]^{\times 9}  \rightarrow [0022379999] ^{+1} \rightarrow[0000113348]
\end{eqnarray*}
In this case the only contribution to the sum comes from the trivial partition $(334)$, with $X(334)  =0$.  All 
other partitions correspond to values of $X(\theta)$ that are greater than $M'_1=2$ and therefore do not contribute.

 \appendix{}  
    
  \section{Derivation of equations $(20)$}
  
  We will use the notation
\begin{eqnarray*}
\langle g|g^d=1\rangle \equiv \{~1,g,g^2, \ldots , g^{d-1}~\},
\end{eqnarray*} 
or the short-hand notation $\langle g|g^d \rangle$,  to denote a general order-$d$ cyclic group generated by $g$. \\

  Throughout this section, $p$ will denote a prime number greater than 2.\\
  
  \begin{prop} 
 For $N=p$ or $N=2p$ ,
   \begin{eqnarray*}
 n_{SM}   = \frac{1}{N\phi(N)}  \sum_{S_d \subseteq  H_N} \phi \left(  d  \right) K(S_d).
    \end{eqnarray*}
  \end{prop}
The proof of this proposition uses the following lemma:\\
 \begin{lemma}:  Let f and g be mappings from the set of subgroups $\{ S_d\}$  of a cyclic group G into the set of nonnegative integers ${\bf Z}$  such that 
 \begin{eqnarray*}
 f(S_d) = \sum_{S_q \supseteq S_d} g(S_q).
 \end{eqnarray*}        
Then
 \begin{eqnarray*}
 g(S_q) = \sum_{S_d \supseteq S_q} \mu \left( \frac{d}{q} \right) ~f(S_d).
 \end{eqnarray*}  
 \end{lemma}
 In both equations, the sums on the right are over all subgroups that contain the subgroup on the left.\\
 ~\\
{\bf Proof of Lemma 4:}  Let $G$ be of order $k$.  Since $G$ is cyclic, for each integer $d$ that divides $k$ there is one and only one subgroup of order $d$ in $G$.  Then we can make the replacements  $f(S_d)\rightarrow f(d),~g(S_q) \rightarrow g(q)$, where $f(d),~g(d)$ are defined for all integers that divide $k$.  Then the condition in the lemma becomes 
   \begin{eqnarray*}
  f(d) |_{d|k}= \sum_{d|q|k} g(q)
 \end{eqnarray*}   
  where the sum is over all $q$ that divides $k$ and is divisible by $d$.   For a given $d$ and $q$ that divide $k$ let 
   \begin{eqnarray*}
  n= \frac{k}{q},~m = \frac{k}{d}
  \end{eqnarray*}
  and
  \begin{eqnarray*}
  ~{\bar f}(m) \equiv f\left( \frac{k}{m} \right),~{\bar g}(n) \equiv g \left( \frac{k}{n} \right) .
 \end{eqnarray*}
 Then the sum above becomes
  \begin{eqnarray*}
   {\bar f}(m) = \sum_{n|m} {\bar g}(n),~~{\rm for~}m|k.
 \end{eqnarray*}   
  Now let ${\widehat f}(s)$ and ${\widehat g}(s)$ be extensions of ${\bar f}(s) $ and ${\bar g}(s) $ to the full domain of the nonnegative integers which preserve this relation but which are otherwise arbitrary:
  \begin{eqnarray*}
  \left. \begin{array} {l} {\widehat f}(s) = {\bar f}(s) \\~\\ {\widehat g}(s) = {\bar g}(s) \end{array} \right\} ~~{\rm for~} s|k ; \\
  {\widehat f}(s) = \sum_{a|s} {\widehat g}(a)~~{\rm for ~all~} s \in {\bf Z}.
  \end{eqnarray*}
   By M\"{o}bius inversion,
   \begin{eqnarray*}
 {\widehat g} (s)  = \sum_{a|s} \mu(s/a) {\widehat f}(a).
     \end{eqnarray*}    
Setting $s=n=k/q$ and $a =m=k/d$ we get
 \begin{eqnarray*}
    g(q) =  \sum_{q|d}~   \mu(d/q)  f(d)
    \end{eqnarray*}
    and
    \begin{eqnarray*}
   g(S_q) =   \sum_{S_d \supseteq S_q}~   \mu(d/q)  f\left( S_d \right) .
     \end{eqnarray*}    
$\triangle$\\

  {\bf Proof of Proposition 4} 
   Let $D = N \phi(N)$. \\

    If $N$ is a prime number then $H_N$ is cyclic and from Lemma 4 we have
    \begin{eqnarray*}
  L(S_q) =  \sum_{S_d \supseteq S_q}~\mu \left( \frac{d}{q} \right) K(S_d)=  \sum_{q|d|D}~\mu \left( \frac{d}{q} \right) K(S_d),
    \end{eqnarray*} 
so
   \begin{eqnarray*}
  n_{SM}&=& \sum_{q|D}  \frac{q}{D}\sum_{q|d|D}  \mu \left( \frac{d}{q} \right) K(S_d)\\
  &=& \sum_{d|D}  K(S_d)  \sum_{q|d}   \frac{q}{D} ~ \mu \left( \frac{d}{q} \right) \\
  &=& \frac{1}{D} ~\sum_{d|D}  K(S_d)  \phi(d),
           \end{eqnarray*}
 where the last line follows from the eq. (18) identity.   Since each divider $d$ corresponds to one and only one subgroup $S_d$, this sum is equal to the sum in the proposition.  \\

 Now let $N$ be equal to twice an odd prime number.  In this case $H_N$ is not cyclic.  We consider  the element $(1,b) \in H_N$, where $b$ is any element of $({\bf Z}/N{\bf Z})^{\times }$, and let $[M]$ be an index set that is invariant under $(1,b)$.  We have
 \begin{eqnarray*}
  [M_0,M_1, M_2, \ldots , M_{N-1}]  = (1,b) *[M]  = [M_{N-1},M_{b^{-1}-1},M_{2b^{-1}-1}, \ldots , M_{(N-1)b^{-1}-1}],
 \end{eqnarray*} 
and so $ M_0 =M_{N-1}$ and $M_k=M_{kb^{-1}-1}$ for all $ k=1, \ldots , N-1.$   $b$ is an odd integer, so if $k$ is even, $(kb^{-1}-1)$ is odd, and vice versa.  Each index with an even subscript is thus equal to an index with an odd subscript, and we can write the sums in $(8b,c)$ as sums over just the even indices:
\begin{eqnarray*}
N= \sum_{k=0}^{N-1} M_k &=&2 \sum_{q=0}^{(N-2)/2}  M_{2q},\\
 \sum_{k=0}^{N-1} kM_k &=& \sum_{q=0}^{(N-2)/2}(2q+2qb^{-1}-1) M_{2q}=2(1+b^{-1})  \sum_{q=0}^{(N-2)/2} qM_{2q} - \frac{N}{2} .\\
  \end{eqnarray*}
Since by assumption $N$ is twice an odd prime,  $\sum kM_k$  is an odd integer and is not congruent to 0 mod $N$.   Thus, index sets that are invariant under any subgroup $S_d$ that contains an element $(1,b)$ for any $b \in ({\bf Z}/N{\bf Z})^{\times }$, (which of course includes the full group $H_N$), do not satisfy condition $(8b)$ and, as a a consequence, $K(S_d)=0$ for such a subgroup.
Then 
 \begin{eqnarray*}
   \sum_{S_d \subseteq  H_N} \phi \left(  d  \right) K(S_d) =
      \sum_{S_d \subseteq  \widehat {H}_N} \phi \left(  d  \right) K(S_d)
     \end{eqnarray*}
where  the proper subgroup
  \begin{eqnarray*}
 \widehat {H}_N = \{0,2,4,  \ldots , N-2\} \times \{1,3,5,\ldots N/2-2,N/2+2, \ldots , N-1\} \subset H_N
 \end{eqnarray*} 
  {\it is} cyclic.  The proof then follows as above for the $N$ equal to a prime number case.  $\triangle$.\\
 ~\\

 If $[M]$ is invariant under a cyclic subgroup $S_d$ of $H_N$, then the conditions ${\bf a}*[M]=[M]$ for all ${\bf a} \in S_d$ can be replaced by the single equation ${\bf g}*[M]=[M]$, where ${\bf g}$ is a generator of $S_d$.  The product subgroup $S^{+}_{N} \times S^{\times}_{1}$ is cyclic, with generator $(1,1)$, since it can be identified with the cyclic group $S^{+}_N= ({\bf Z}/N{\bf Z})^+$.  Its generator, expressed as a permutation on $[M]$, is
\begin{eqnarray*}
(1,1)* = (\,0123\cdots (N-1)\,)
\end{eqnarray*}
The only solution to $(1,1)*[M]=[M]$ and to $N= \sum_q M_q$, (eq.(8c)), is $[M]=[111\cdots 1]$.  This solution  however satisfies conditions (8b) only for odd $N$.  Consequently,
\begin{eqnarray}
K(  S^{+}_{N} \times S^{\times}_{1}) = \left\{ \begin{array} {c} 1~~{\rm for~odd~} N, \\~\\ 0 ~~{\rm for~even~} N .
\end{array} \right. 
\end{eqnarray}

 In addition, since all index sets are invariant under the identity subgroup $S_1 =S^{+}_1 \times S^{\times }_1 = \{(0,1)\}$, 
\begin{eqnarray}
K(S_1) = \frac{1}{2N}\sum_{d|N } ~\phi \left( \frac{N}{d} \right)  \left( \begin{array} {c}  2d \\ ~\\ d \end{array} \right).
\end{eqnarray}

We will use the lemma below to calculate $n_{SM}(N)$ for $N=p$ and $N=2p$:
\begin{lemma} 
If $[M]$ is invariant under the action of the group element $(a,b),~b\neq 1$, then $P^c [M]$ is invariant under the action of
$(a+(1-b^{-1})c,b)$
\end{lemma}
~\\
Proof:  Acting on the multiplicity index set $[M]=[M_0,M_1,\dots ,M_{N-1}]$, the circular permutation operator $P^c$ permutes it to
\begin{eqnarray*}
P^c[M] =[M_{N-c},M_{1-c}, \ldots, M_{N-1-c}]
\end{eqnarray*}
Then
\begin{eqnarray*}
(a+(1-b^{-1}c,b)*(P^c[M]) &=& (a+(1-b^{-1})c,b)*[M_{N-c},M_{1-c}, \ldots, M_{N-1-c}]\\
&=&  [M_{(N-c)b^{-1}- a-c+b^{-1}c},M_{(1-c)b^{-1}-a-c+b^{-1}c}, \ldots, M_{(N-1-c)b^{-1}-a-c+b^{-1}c}]\\
&=&  [M_{- a-c},M_{b^{-1}-a-c}, \ldots, M_{(N-1)b^{-1}-a-c}]\\
&=&  P^c (a,b)*[M] = P^c[M]~~\triangle
\end{eqnarray*}

 \subsubsection{N =p}
 
 For $N=p$, we have, from Proposition 4,
    \begin{eqnarray*}
 n_{SM}(p)   &=& \frac{1}{p(p-1)}  \sum_{S_k \subseteq  H_p} \phi \left(  k  \right) K(S_k),\\
 H_p &=& \{0,1,2, \ldots , p-1\} \times \{ 1,2,\ldots , p-1\}.
    \end{eqnarray*} 
   
$H_{p}$ is cyclic.  All of its subgroups are therefore cyclic and are equal to the direct product of a subgroup of 
$( {\bf Z} /p{\bf Z})^{+}$ and of $( {\bf Z} /p{\bf Z})^{\times}$.   The subgroups of $( {\bf Z} /p{\bf Z})^{+}$ are $ S^{+}_1 =\{0\}$ 
and $S^{+}_p=\{0,1,2, \ldots, p-1\}$, with generators 0 and 1, respectively.  The subgroups of $H_p$ are then of the form
\begin{eqnarray*}
&&S^{+}_{1} \times S^{\times}_{d} = \{0\} \times S^{\times}_d = \langle (0,g)|(0,g)^{d} \rangle,\\
{\rm and} && \\
&&S^{+}_{p} \times S^{\times}_{d} =\{0,1,2, \ldots, p-1\} \times S^{\times }_d = \langle (1,g)|(1,g)^{pd} \rangle,
\end{eqnarray*}
where $d$ divides $(p-1)$ and $g$ is a generator of $S^{\times }_d \subseteq ({\bf Z}/p{\bf Z})^{\times}$.   Then
     \begin{eqnarray}
 n_{SM}(p)  &=& \frac{1}{p(p-1)} \left\{ K(S^+_1 \times S^{\times}_1 ) + (p-1) K(S^+_p \times S^{\times}_1) 
    + \sum_{d|(p-1)}^{d>1} \phi(d)  [ K(S^+_1 \times S^{\times} _d) + (p-1) K(S^+_p \times S^{\times} _d) ] \right\}, \nonumber \\
    \end{eqnarray} 
where we've used the relation $\phi(pd) = \phi(p)\phi(d) = (p-1)\phi(d)$.  The $d=1$ functions, $K(S^+_1 \times S^{\times}_1 ) $ and $ K(S^+_p \times S^{\times}_1)$, are given by eqs.(A2) and (A1), respectively, and we have 
\begin{eqnarray}
 K(S^+_1 \times S^{\times}_1 ) + (p-1) K(S^+_p \times S^{\times}_1)
     &=&  \frac{1}{2p}\left( \begin{array} {c}  2p \\ ~\\ p \end{array} \right) + \frac{p^2-1}{p} .
\end{eqnarray}

For $d>1$, we show that 
\begin{eqnarray*}
K (  S^+_p \times S^{\times}_d ) = K(S^{+}_1 \times S^{\times}_d )   =\left( \begin{array} {c} 
    \frac{2(p-1)}{d} \\ ~\\ \frac{p-1}{d} \end{array} \right).
\end{eqnarray*}
    To prove the first identity,  we set $ (a,b)=(0,g)$ and $  (a+(1-b^{-1})c,b) =(1,g) $  in Lemma 5 above.  Then $c= g(g-1)^{-1}$.  For each $[M]$ that satisfies $(0,g)*[M] =[M]$, there is a circular permutation $[M']= P^{g(g-1)^{-1}}[M]$  that satisfies $(1,g)*[M']= [M']$, and vice versa.  The number of solutions to both eqs.(8b,c) and to $(1,g)*[M] =[M]$ is therefore the same as the number of solutions to (8b,c) and $(0,g)*[M] =[M]$,  and so $K(S^{+}_p\times S^{\times}_d )  =K(S^{+}_1 \times S^{\times}_d)$.\\
    
    [Note that this argument fails for $d=1$  since in this case $g=1$.]\\ 
    
To prove the second identity, we identify $S^{+}_1 \times S^{\times}_d$ with $S^{\times}_d$ and set $m=(p-1)/d$.  $S^{\times}_d$ has $(m-1)$ cosets that we denote as $\{ c_1\}, \ldots ,\{c_{m-1}\}$, with $c_k \in \{c_k\} \subset [1,2,\ldots , (p-1)]$.    Multiplication mod $p$ by $g$ induces a cyclic permutation $g*$ on $[012 \cdots (p-1)]$:
  \begin{eqnarray*}
  g* = (0) (\,1\,g \cdots g^{d-1}\,) (\,c_1\, (gc_1) \cdots( g^{d-1}c_1)\,) \cdots (\,c_{m-1}\, (gc_{m-1}) \cdots (g^{d-1}c_{m-1})\,).
  \end{eqnarray*}  
Then  $[M]$ is invariant under $g*$ if and only if
\begin{eqnarray}
M_1 &=& M_g =\cdots = M_{g^{d-1}} \nonumber, \\
M_{c_1} &=& M_{gc_1} = \cdots =M_{g^{d-1}c_1},\\
&&~~~ \vdots \nonumber \\
M_{c_{m-1}} &=& M_{gc_{m-1}} = \cdots = M_{g^{d-1}c_{m-1}}\nonumber.
\end{eqnarray}
Therefore,
\begin{eqnarray*}
 &&\sum_{q=0}^{N-1} M_q = M_0 + d(M_1 +M_{c_1} + \cdots + M_{c_{m-1}} ),\\
&& {\rm and~}  \\
&&\sum_{q=0}^{N-1}  qM_q = (~M_1 + c_1 M_{c_1} + \cdots + c_{m-1}M_{c_{m-1}}~) (1 + g+ \cdots +g^{d-1}).
 \end{eqnarray*}
 We have
 \begin{eqnarray*}
 1 + g+ \cdots +g^{d-1} = \frac{g^{d}-1}{g-1} = \frac {kp}{g-1}
 \end{eqnarray*} 
for some integer $k$, where the last equality follows from $g^d \equiv 1$ mod $p$.   $(g-1)$ then divides $kp$; since $p$ is prime, it must divide $k$.  Therefore, for $d>1$, $1 +g + \cdots + g^{d-1} \equiv 0 $ mod $p$ and  condition (8b) is satisfied for any index set $[M]$ that satisfies equations (A5).\\

Condition (8c) and the first equation above require that 
 \begin{eqnarray*}
  M_1 +M_{c_1} + \cdots + M_{c_{m-1}}  = \frac{N-M_0}{d} = m - \frac{M_0 -1}{d} ,
\end{eqnarray*}
so we must have that 
\begin{eqnarray*}
M_0 = 1 + d \alpha~\geq 1,~~ M_1 +M_{c_1} + \cdots + M_{c_{m-1}}   = m - \alpha~ \geq 0 ,
\end{eqnarray*}
for some integer $\alpha \in [0,m]$;  ($M_0=0$ is excluded since $d>1$).   For a given value of $\alpha$,  the set $ \{ M_1+1, \ldots , M_{c_{m-1}} +1 \}$ is a  composition of $(2m-\alpha)$.   The number of compositions of an integer $n$ into $k$ parts is $\left( \begin{array} {c} n-1 \\ k-1 \end{array} \right)$,  so the number of solutions to the set of equations in (A5) which also satisfy (8c) is 
\begin{eqnarray}
K(S^+_1 \times S^{\times }_d) =  \sum_{\alpha=0}^m \left( \begin{array} {c} 2m-\alpha -1 \\~\\  m-1 \end{array} \right) 
    &=& \sum_{\beta=m-1}^{2m-1} \left( \begin{array} {c} \beta \\~\\  m-1 \end{array} \right)
    = \left( \begin{array} {c} 2m \\ ~ \\ m \end{array} \right).
    \end{eqnarray}
Putting equations (A4) and (A6) into (A3), we get eq.(20a).\\

  \subsubsection{N =2p}

In this case, 
        \begin{eqnarray*}
 n_{SM}(2p)   &=& \frac{1}{2p(p-1)}  \sum_{S_k \subseteq  {\widehat H}_{2p}} \phi \left(  k  \right) K(S_k),\\
  \widehat {H}_{2p} &=& \{0,2,4,  \ldots , 2p-2\} \times \{1,3,5,\ldots p-2,p+2, \ldots , 2p-1\} .
    \end{eqnarray*}
   ${\widehat H}_{2p}$ is cyclic and its subgroups are of the form
\begin{eqnarray*}
&&S_1^+\times S^{\times}_d  =\{0\} \times S^{\times}_d = \langle (0,g)|(0,g)^p \rangle, \\
&&S^+_{p } \times S^{\times}_d = \{0,2,4,\ldots , 2p-2\}\times S^{\times}_d = \langle (2,g)|(2,g)^{pd} \rangle,
\end{eqnarray*}
where  $S^{\times}_d = \{1,g, \ldots ,g^{d-1} \} \subseteq ({\bf Z}/(2p){\bf Z})^{\times}$ is a $d$-element cyclic subgroup generated by the odd integer $g$.   As indicated above, $S^+_1 \times S^{\times }_d$ and $S^+_p \times S^{\times }_d $ are both cyclic, with generators $(0,g)$ and $(2,g)$, respectively.  Then
  \begin{eqnarray}
 n_{SM}(2p)   &=& \frac{1}{2p(p-1)} \left\{ K( S^+_1 \times S^{\times}_1) + (p-1) K(S^+_p \times S^{\times}_1)
  +\sum_{d|(p-1)}^{d>1}   \phi(d) \left[ K(S^+_1 \times S^{\times}_d) + (p-1) K(S^+_p \times S^{\times}_d) \right] \right\}. \nonumber \\
    \end{eqnarray}
    
 For $d=1$ we have 
 \begin{eqnarray*}
K(S^+_1 \times S^{\times}_1) &=& \frac{1}{4p}\sum_{q=1,2,p,2p } ~\phi \left( \frac{2p}{q} \right) 
     \left( \begin{array} {c}  2q \\ ~\\ q \end{array} \right)
  = \frac{2(p-1)}{p}+  \frac{1}{4p} \left[    \left( \begin{array} {c}  2p \\ ~\\ p \end{array} \right)   + \left( \begin{array} {c}  4p \\ ~\\ 2p \end{array} \right)   \right]\\
\end{eqnarray*}
The cyclic subgroup $S_p^+ \times S_1^{\times}$ has the generator $(2,1)$.  As a permutation on $[M]$, 
\begin{eqnarray*}
(2,1)* = (024\cdots 2p-2)(135 \cdots 2p-1).
\end{eqnarray*}
 The only solutions to $(2,1)*[M] =[M]$ which satisfy $(8b,c) $ are 
$ [M]= [2020\cdots 20]$ and $[0202 \cdots 02]$, and so $ K(S^{+}_p \times S^{\times}_1) =  2.$\\

The identity $K(S^{+}_p \times S^{\times}_d ) =  K(S^{+}_1 \times S^{\times}_d )$ for $d>1$ follows by a similar argument as that in the previous section, using Lemma 5, with $c =2g(g-1)^{-1}$.   It remains then to calculate $K(S^+_1 \times S^{\times}_d)$ for $d>1$, which we
do in essentially the same fashion as was done for the $N=p$ case.  We identify $S^+_1 \times S^{\times}_d$ with $S^{\times}_d$ and consider the permutation $(0,g)*= g*$ acting on $[012\ldots N-1]$:
\begin{eqnarray*}
&&g* = (0) (1g\cdots g^{d-1}) (c_1~ g\cdot c_1 \cdots g^{d-1}\cdot c_1) \cdots (c_{m-1} ~g \cdot c_{m-1} \cdots
g^{d-1}\cdot c_{m-1}) \\
&&~~~~~~~~~(2~ 2g \cdots 2g^{d-1}) (2 c_1~ 2g c_1 \cdots 2g^{d-1} c_1) \cdots (2c_{m-1}~ 2g  c_{m-1} \cdots
2g^{d-1} c_{m-1}) (p).
\end{eqnarray*}
Again the $c_k$'s label the cosets of $S^{\times} _d$ and $m$ is defined as $(p-1)/d$.  If $[M]$ is  invariant under $S^{\times}_d$ then $g*[M] =[M]$, and conditions (8b,c) become
\begin{subequations}
\begin{eqnarray}
&& 2p=\sum_q M_q = M_0 +M_p + d(M_1+ M_{c_1} + \cdots  +M_{c_{m-1}}
+M_2 +M_{2c_1}+ \cdots + M_{2c_{m-1}} )\\
&& {\rm or} \nonumber \\
 && M_1+ M_{c_1} + \cdots  +M_{c_{m-1}}+M_2 +M_{2c_1}+ \cdots + M_{2c_{m-1}}= \frac{2p-M_0-M_p}{d}
= 2m - \frac{M_0+M_p-2}{d},~\\
&&{\rm and} \nonumber \\
&&0 \equiv  \sum_q qM_q = pM_p +(1+g+ \cdots + g^{d-1})(~M_1 +c_1M_{c_1} +\cdots + c_{m-1} M_{c_{m-1}} \nonumber\\
&& ~~~~~~~~~~~~~~~~~~~~~~~~~~~~~~~~
~~~~~~~~~~~~~~~~~~~~~~~~~~~ + 2M_2 +2c_1 M_{2c_1} + \cdots + 2c_{m-1}M_{2c_{m-1}}~).
\end{eqnarray}
\end{subequations}
 Since $g^d \equiv 1 $ mod $2p$, $g^d = q(2p)+1$ for some integer $q$.  Then, for $d>1$,
 \begin{eqnarray*}
 Tr(S^{\times}_d)= 1 + g+ \cdots +g^{d-1} = \frac{g^{d}-1}{g-1} = \frac {2qp}{g-1} .
 \end{eqnarray*} 
The trace of $S^{\times}_d$ is a sum of $d$ odd integers, which is odd or even depending on where $d$ is odd or even.   Since $p$ is odd, the parity of the integer $2q/(g-1)$ is the same as the parity of $d$.  As a result,
\begin{eqnarray*}
&&Tr (S^{\times}_d) \equiv \left\{ \begin{array} {l} 0 {\rm ~mod~} 2p~~{\rm if~} d ~{\rm is~even};\\~\\
                                                                    p {\rm ~mod~} 2p~~{\rm if~} d >1~{\rm is ~odd}.
                                                                    \end{array} \right.                                         
\end{eqnarray*}

We consider the cases $d$ equal to an even integer and $d$ equal to an odd integer greater than 1 separately:\\
~\\
(I)   $d=2k$.   Eq. (A8c) then requires that $M_p$ is even.  From (A8a/b) $M_0$ must also be even.  Setting  $M_0= 2m_0,~M_p=2m_p$, we have from (A8b)
\begin{eqnarray*}
 &&M_1+ M_{c_1} + \cdots  +M_{c_{m-1}}+M_2 +M_{2c_1}+ \cdots + M_{2c_{m-1}}= 2m - \alpha
 \end{eqnarray*}
 where $\alpha$ is an integer less than or equal to  $2m$ defined as
 \begin{eqnarray*}
 && \alpha =\frac{M_0+M_p-2}{d}=  \frac {m_0 +m_p -1}{k};\\
  or~&& m_0+m_p = k \alpha +1, ~\geq 0.
 \end{eqnarray*}
The set $\{ M_1+1, M_{c_1} +1, \ldots , M_{2c_{m-1}}+1\} $ is then a composition of $(4m-\alpha)$.  For each $\alpha$ there are $( k\alpha +2) $ combinations of $M_0$ and $M_p$ such that $M_0/2+M_p/2= k\alpha +1$.  For $k>1$, $\alpha$ can take on values from 0 to $2m$, while for $k=1$ the range of $\alpha$ extends downwards to -1. The total number of solutions $\{M_0, M_1, \ldots , M_{2c_{m-1}}\}$ is then

 \begin{eqnarray*}
K(S^+_1 \times S^{\times}_d )|_{d ~{\rm even}} &=& \sum_{\alpha=0}^{2m} ( k \alpha +2) \left( \begin{array}{c}  4m-\alpha-1\\~\\ 2m-1 \end{array} \right)  
+\delta_{d,2}  \left( \begin{array}{c}  4m \\~\\ 2m-1 \end{array} \right) \\
 &=& \left\{ \begin{array} {l} ~~~~~~~~~~~~~~ \left( \begin{array}{c}  2p \\~\\ p \end{array} \right),~~~~~~~~~~~~~~~~~~~~~~d=2, \\~
\\    2p \left( \begin{array}{c} 4m \\~\\ 2m \end{array} \right)  
- (p-1)  \left( \begin{array}{c} 4m+1 \\~\\ 2m+1 \end{array} \right) ,~~d>2.    \end{array} \right.
\end{eqnarray*}

~\\
(II)  $d=2k+1,~k>0$.  In this case,
\begin{eqnarray*}
 \sum qM_q \equiv  p( M_p   +M_1+  c_1M_{c_1} +\cdots + c_{m-1} M_{c_{m-1}}~)~{\rm mod~} 2p,
\end{eqnarray*}
and condition (8b) requires that the sum $ M_p+ M_1 +c_1M_{c_1} +\cdots + c_{m-1} M_{c_{m-1}}$  is even.  But since all the $c_k$'s are odd integers, this condition becomes 
\begin{subequations}
\begin{eqnarray}
  M_p+ M_1 +M_{c_1} +\cdots +  M_{c_{m-1}} =2n_1
  \end{eqnarray}
  for some nonnegative integer $n_1$.  Then from (A8a)
\begin{eqnarray*}
&&  2p = M_0 +M_p +d( 2n_1 -M_p +M_2 +M_{2c_1}+ \cdots + M_{2c_{m-1}})\\
&&~~~= M_0 +2dn_1 -(d-1)M_p +d(M_2 +M_{2c_1}+ \cdots + M_{2c_{m-1}})
\end{eqnarray*}
so that the sum  $M_0+ M_2 + \cdots +M_{2c_{m-1}} $ is also even:
\begin{eqnarray}
 M_0+ M_2 + \cdots +M_{2c_{m-1}} = 2n_2,~~n_2 \geq 0.
\end{eqnarray}
We define $Z_1$ and $Z_2$ as
\begin{eqnarray*}
Z_1= 2n_1 -M_p,~~Z_2= 2n_2-M_0.
\end{eqnarray*}
Then
\begin{eqnarray}
 2p=M_0+M_p +d(Z_1+Z_2),
\end{eqnarray}
\end{subequations}
and we re-express  (A9a,b,c) as
\begin{subequations}
\begin{eqnarray}
 M_1 +M_{c_1} +\cdots +  M_{c_{m-1}} &=&Z_1,\\
 M_2 + M_{2c_1} +\cdots +M_{2c_{m-1}} &=&  Z_2,\\
 M_0 +M_p  &=& 2p -d( Z_1+Z_2),
 \end{eqnarray}
 \end{subequations}
 which formally decouples the three sets of indices.  The first two equations require that 
 $0 \leq Z_1,Z_2 \leq 2m$.  The third equation requires that $0 \leq Z_1+Z_2 \leq 2m$, and that 
 $ Z_1+Z_2$ has the same parity as $M_0+M_p$.  Equations (A9) however have the additional requirement
  that the pairs $(Z_1,M_p)$ and $(Z_2, M_0)$ individually have the same parity, so we must impose 
  this by hand on the solutions to equations (A10).  (If one set has the same parity, then so does the other.)\\

For a given $(Z_1,Z_2)$,  the number of solutions to (A10c) without this restriction is 
    \begin{eqnarray*}
    \left(  \begin{array} {c} 2p+1-d(Z_1 +Z_2) \\ ~\\ 1 \end{array} \right)   \equiv W.
   \end{eqnarray*}   
If $Z_1+Z_2$ is odd, then $W$ is even and there are $W/2$ compositions $(M_0,M_p)$ in which $M_p$ is even,
and $W/2$ compositions in which it is odd.  So to count only the compositions in which $Z_1$ and $M_p$ have the
same parity, we divide $W$ by 2:  $W \rightarrow W/2$, for $Z_1+Z_2$ odd.\\

If on the other hand $Z_1+Z_2$ is even, then $W$ is odd, and there are $(W\pm1)/2$ such compositions for $Z_1$
even/odd, respectively.\\

  So, to impose the requirement that $(Z_1,M_p)$,  (or equivalently $(Z_2,M_0)$), have the same
parity we make the replacement
\begin{eqnarray*}
W \rightarrow \frac{1}{2}~\left[ W + \frac{(-1)^{Z_1} +(-1)^{Z_2} }{2} \right].
\end{eqnarray*}
The number of solutions to (A10) with this requirement is then
   \begin{eqnarray*}
 K(S^+_1 \times S^{\times}_d)|_{d~{\rm odd}, >1} = \sum_{Z_1=0}^{2m} \sum_{Z_2=0}^{2m-Z_1}\frac{1}{2}
   \left[ 2p+1-d(Z_1 + Z_2) + \frac{(-1)^{Z_1} +(-1)^{Z_2} }{2} \right]
   \left(  \begin{array} {c} Z_1+m -1 \\ ~\\ m-1 \end{array} \right) 
   \left(  \begin{array} {c} Z_2+m -1 \\ ~\\ m-1 \end{array} \right) .
      \end{eqnarray*}  
  From the $(Z_1,Z_2)$ symmetry of the summand we can replace $( ~(-1)^{Z_1} + (-1)^{Z_2}~)/2$ by $(-1)^{Z_1}$.  Evaluating
  this expression we get
     \begin{eqnarray*}
 K(S^+_1 \times S^{\times}_d)|_{d~{\rm odd}, >1}  =   \frac{4p-1}{2} \left(  \begin{array} {c} 4m \\ ~\\ 2m \end{array} \right) 
   - (p-1)  \left(  \begin{array} {c} 4m+1  \\ ~\\ 2m+1 \end{array} \right)   + \sum_{Z_1=0}^{2m}
    \frac{(-1)^{Z_1}}{2}  \left(  \begin{array} {c} Z_1+m -1 \\ ~\\ m-1 \end{array} \right) 
   \left(  \begin{array} {c} 3m-Z_1 \\ ~\\ m \end{array} \right) .
      \end{eqnarray*}  
      To evaluate the remaining sum we use the following lemma:   
   \begin{lemma}
        \begin{eqnarray*}
    \sum_{k=0}^{X}  \left(  \begin{array} {c} k  \\ ~\\ m-1 \end{array} \right)
  \left(  \begin{array} {c} X-k  \\ ~\\ m \end{array} \right)    (-1)^k = (-1)^{m-1} \left(  \begin{array} {c} \lceil X/2 \rceil \\ ~\\ m \end{array} \right)\\
\end{eqnarray*}  
\end{lemma}
where $\lceil ~\rceil$ denotes the ceiling function.  To prove this identity we write it in the form  
     \begin{eqnarray*}
   (A)  &&   \sum_{k=0}^{2n}  \left(  \begin{array} {c} k  \\ ~\\ m-1 \end{array} \right)
  \left(  \begin{array} {c} 2n-k  \\ ~\\ m \end{array} \right)    (-1)^k = (-1)^{m-1} \left(  \begin{array} {c} n  \\ ~\\ m \end{array} \right) ,\\
(B)   &&   \sum_{k=0}^{2n+1}  \left(  \begin{array} {c} k  \\ ~\\ m-1 \end{array} \right)
  \left(  \begin{array} {c} 2n+1-k  \\ ~\\ m \end{array} \right)    (-1)^k = (-1)^{m-1} \left(  \begin{array} {c} n+1  \\ ~\\ m \end{array} \right).  \\
\end{eqnarray*}  
These two identities can be proven by using the identity
\begin{eqnarray*}
  \frac{x^m}{(1-x)^{m+1}} =  \sum_{q=0}^{\infty} \left(  \begin{array} {c} q  \\ ~\\ m \end{array} \right) x^q
  \end{eqnarray*}
  and expressing the product
  \begin{eqnarray*}
     \frac{~~(-x)^{m-1}}{(1+x)^m}~ \frac{x^m}{(1-x)^{m+1}}
        \end{eqnarray*}  
first as the product of two sums,
  \begin{eqnarray}
   \frac{~~(-x)^{m-1}}{(1+x)^m}~ \frac{x^m}{(1-x)^{m+1}}
 = \left(\sum_{k=0}^{\infty} \left(  \begin{array} {c} k  \\ ~\\ m-1 \end{array} \right) (-1)^k x^k \right)
 \left( \sum_{q=0}^{\infty} \left(  \begin{array} {c} q  \\ ~\\ m \end{array} \right) x^q\right),
         \end{eqnarray}   
then as a single sum,
  \begin{eqnarray}
  \frac{~~(-x)^{m-1}}{(1+x)^m}~ \frac{x^m}{(1-x)^{m+1}}   = \frac{(-1)^{m-1}x^{2m-1}(1+x)} {(1-x^2)^{m+1}}
   = (-1)^{m-1} \sum_{p=0}^{\infty} \left(  \begin{array} {c} p  \\ ~\\ m \end{array} \right)( x^{2p-1} +x^{2p}),
           \end{eqnarray}  
and then equating the even-exponent and the odd-exponent terms on the right-hand side of (A11) to those on the right-hand side of (A12). \\

 Now setting $k=Z_1+m-1$ and $X=4m-1$  we have from this lemma
     \begin{eqnarray*}
 \sum_{Z_1=0}^{2m}  \left(  \begin{array} {c} Z_1+m-1  \\ ~\\ m-1 \end{array} \right)
  \left(  \begin{array} {c} 3m-Z_1  \\ ~\\ m \end{array} \right)    (-1)^{Z_1+m-1} = (-1)^{m-1} \left(  \begin{array} {c} 2m  \\ ~\\ m \end{array} \right).
   \end{eqnarray*}   
  And so
     \begin{eqnarray*}
  K(S^+_1 \times S^{\times}_d)|_{d {\rm ~odd},>1}   = \frac{4p-1}{2}  \left(  \begin{array} {c} 4m \\ ~\\ 2m \end{array} \right) 
   - (p-1) \left(  \begin{array} {c} 4m+1  \\ ~\\ 2m+1 \end{array} \right)
   +  \frac{1}{2} \left(  \begin{array} {c} 2m  \\ ~\\ m \end{array} \right) .
   \end{eqnarray*}  
Putting all these results together, we get eq.(20b).   

\section{}

We have, from the discussion in the proof of Theorem 3, that 
   \begin{eqnarray}
A^{(k)}_{0\cdots 0\,1 \cdots 1\, A_1 \cdots  A_{p+1}} \propto \left.   \sum_{ q_1, \ldots,  q_p }
                                                 Q (k-A_1q_1 -\cdots - A_pq_p; M_1,N-1~|~ q_1, \ldots , q_p) \right|_{q_i \neq q_l}.
 \end{eqnarray}

 To determine the relation between this partition function and $A^{(k)}_{0 \cdots 0 \, 1 \cdots 1 \, A \cdots A}$, we first
consider the case where all multiplicities $M_0, M_1,\ldots, M_{N-1}$ have no common factor.  Then, as discussed earlier, all the cyclic-permutation equivalence classes $\{ \sigma_0 =A_{p+1}, \, \sigma_1 , \ldots , \sigma_{M-1}\}$ have $N$ elements.  The coefficient
 $A^{(k)}_{0\cdots 0 \,1 \cdots 1 \, A_1 \cdots A_{p+1}} $ is then proportional to $N$ times a sum of $Q$ over all possible $q$'s divided, to avoid 
 overcounting, by the factorials of the multiplicities:
  \begin{eqnarray*}
A^{(k)}_{0\cdots 0\,1 \cdots 1\, A_1 \cdots  A_{p+1}} \propto \left.  \frac{ N\, M_{A_{p+1}}}{ M_2!  \cdots M_{N-1}! }   ~ \sum_{ q_1, \ldots,  q_p }
                                                 Q (k-A_1q_1 -\cdots - A_pq_p; M_1,N-1~|~ q_1, \ldots , q_p) \right|_{q_i \neq q_l} .
 \end{eqnarray*}
(The factor of $M_{A_{p+1}}$ occurs in the numerator because the multiplicity of $A_{p+1}$ within the indices $\{A_1, \ldots , A_p\}$ is only 
$M_{A_{p+1}}-1$.)\\

Now consider a particular permutation $[ \sigma ] = [ A_{p+1}, \sigma_1, \ldots  \sigma_{N-1} ]$  corresponding to a particular set $\{q_1, \ldots , q_{p} \}$ and a particular partition $(p_1, \ldots , p_{M_1})$.   $[\sigma]$ contains, in addition to the $A_{p+1}$ at the $0$ position, an
additional $(M_{A_{p+1}}-1)$ $A_{p+1}$'s at various $q_j$ positions.   The equivalence class containing this $[\sigma]$, i.e., $\{ \sigma\}$, also contains (among others) all permutations of the form $[\sigma'] =  [\sigma'_0 = \sigma_{q_j} = A_{p+1}, \sigma'_1, \ldots,\sigma_{0}=A_{p+1}, \ldots  \sigma'_{N-1} ]$, obtained by cycling the $A_{p+1}$ at the $q_j$ position to the $0$ position, (and the $A_{p+1}$ at the $0$ position to somewhere else).    These permutations will also be counted in the sum over the $q$'s, under a different set of $q$'s and $p$'s.    Since there are $(M_{A_{p+1}}-1)$ such $q$'s, each equivalence class will be counted $M_{A_{p+1}}$ times in the sum above.  We must therefore divide $N$ by an additional $ M_{A_{p+1}}$ to get the number of {\it distinct} permutations:

   \begin{eqnarray*}
A^{(k)}_{0\cdots 0\,1 \cdots 1\, A_1 \cdots  A_{p+1}} &=&\left.  \frac{ N}{ M_2!  \cdots M_{N-1}! }   ~ \sum_{ q_1, \ldots,  q_p }
                                                 Q (k-A_1q_1 -\cdots - A_pq_p; M_1,N-1~|~ q_1, \ldots , q_p) \right|_{q_i \neq q_l}.  
 \end{eqnarray*}
For example, consider the index set $[0011255]  $ and the two equivalence classes below that appear in the sum over $q_1,q_2$:
\begin{eqnarray*}
&&\{5001125\}:~~ \{p_1,p_2\} = \{3,4\};~\{q_1,q_2\} =\{5,6\};\\
&&\{5500112\}:~~ \{p_1,p_2\} = \{4,5\};~\{q_1,q_2\} = \{1,6\} .
\end{eqnarray*}
Each of the permutations $[5001125]$ and $[5500112]$ are in both of these equivalence classes, so we must divide by $2$.\\

Now let the multiplicities have a common factor $d$:  gcd($M_0,M_1,M_2, \ldots, M_{N-1})=d$.   Then there will be some equivalence classes of the form $\{A_{p+1}x \cdots y\,A_{p+1}x \cdots y \cdots  A_{p+1}x \cdots y\}$ for which 
$\sigma_j = \sigma_{j+N/d}$, so that the ``phrase'' $A_{p+1}x \cdots y$ is repeated $d$ times.  These classes have only $N/d$ distinct elements.  The sum in (B1) will then consist of two sums, one over the $N$-element classes and one over the $(N/d)$-element classes, which we will express in short-hand notation as

\begin{eqnarray*}
   A^{(k)}_{0\cdots 0\,1 \cdots 1\, A_1 \cdots  A_{p+1}} =   \frac{ M_{A_{p+1}}}{ M_2!  \cdots M_{N-1}! } \left\{    \frac{N}{M_{A_{p+1}}}~ \sum_I  Q ~
     +~   \frac{N/d}{X}  \sum_{II}   Q     \right\}                                               
\end{eqnarray*}
~\\
where division by some factor $X$ is needed to correct for overcounting due to the duplication of equivalence classes in the sum.  The index $A_{p+1}$ occurs $M_{A_{p+1}}$ times in the entire index set $[A_{p+1}x \cdots y \,A_{p+1} x \cdots y \cdots A_{p+1} x \cdots y]$, and 
$M_{A_{p+1}}/d$ times in the phrase $A_{p+1}x \cdots y$; $M_{A_{p+1}}/d$ is thus the factor $X$ by which each distinct equivalence class is overcounted:  $(N/d)/(M_{A_{p+1}}/d) =N/M_{A_{p+1}}$. \\

As an example of this situation, we consider $[001122333333]$, with $N=12,~M_0=M_1=M_2 =2,  M_3= 6,~d=2$, and the equivalence classes:
\begin{eqnarray*}
&&\{303312 \, 303312\}:~~\{p_1,p_2\} = \{ 4,10 \};~\{q_1,q_2,q_3,q_4,q_5,q_6,q_7\} = \{2,3,5,6,8,9,11\} ;\\
&&\{331230\, 331230\}:~~\{p_1,p_2\} = \{ 2,8 \};~~\{q_1,q_2,q_3,q_4,q_5,q_6,q_7\} = \{1,3,4,6,7,9,10\} ;\\
&&\{312303\,312303\}:~~\{p_1,p_2\} = \{ 1,7 \};~~\{q_1,q_2,q_3,q_4,q_5,q_6,q_7\} = \{2,3,5,6,8,9,11\} .\\
\end{eqnarray*}
Each equivalence class has $N/d =6$ distinct elements, and each class is overcounted by a factor of $3 = M_3/d$.   The overall correction factor is
then  $6/3 =N/M_3$.  (Although classes 1 and 3 have the same $q$ values, their $p$ values are different and will therefore be counted separately by the partition function.)\\

Since the correction factor is the same as for the first sum, we can use the same multiplying factor for both without any further need to consider them separately:

 \section{ Proof of Lemma 1}

For a given set $\{ \kappa_1, \ldots , \kappa_p\}$ let  
 \begin{eqnarray*}
 \beta_j = \kappa_1 + \cdots + \kappa_j,~~1 \leq j \leq p.
 \end{eqnarray*}
 Then $\kappa_1 = \beta_1$ and  $\kappa_j = \beta_j- \beta_{j-1}$ for  $ j \geq 2$.   In the sums over the $\kappa$'s in the statement of the lemma, the restriction that $ \kappa_1 + \cdots + \kappa_p \leq m$ requires that each $\beta_j \leq m$.  And, since each $\kappa_j$ has the range $[0,m]$, each $\beta_{j-1}$ has the range $[0,\beta_j ]$.  We therefore have 
  \begin{eqnarray*}
 && \sum_{0 \leq \kappa_1, \ldots , \kappa_p \leq m}  \omega^{\kappa_{1} q_1  +\kappa_2 q_2 + \cdots +
                 \kappa_p q_p} H[ m- ( \kappa_1 + \cdots + \kappa_p)] \\
                 && ~~~~~~~~~~~~~~~~~~~~~~~~~~~~~~~~~~~~=  \sum_{\beta_{p}=0}^{m}   \sum_{\beta_{p-1}=0}^{\beta_{p}}  
                 \cdots  \sum_{\beta_2=0}^{\beta_3} \sum_{\beta_1=0}^{\beta_2}  \omega^{(\beta_{p}-\beta_{p-1})q_{p}}
                 \cdots  \omega^{(\beta_2-\beta_1)q_2} \omega^{\beta_1 q_1} ,
 \end{eqnarray*}
and the statement of Lemma I is equivalent to  
 \begin{eqnarray*}
\sum_{m=0}^{N-1} (-y)^m  ~  \prod_{s=1}^{p} \frac{1} { 1+y\omega^{q_s} }
                 &=&  \sum_{m=0}^{N-1-p} (-y)^m   ~  \sum_{\beta_{p}=0}^{m}   \sum_{\beta_{p-1}=0}^{\beta_{p}}  
                 \cdots  \sum_{\beta_2=0}^{\beta_3} \sum_{\beta_1=0}^{\beta_2}  \omega^{(\beta_{p}-\beta_{p-1})q_{p}}
                 \cdots  \omega^{(\beta_2-\beta_1)q_2} \omega^{\beta_1 q_1}  .
 \end{eqnarray*}
The proof of Lemma I, by induction on $p$, is more direct when stated in this form.  \\

To prove the lemma we first set $p$ equal to 1 on the right side,  multiply that expression by $(1+y \omega^q)$, and expand: 
 \begin{eqnarray}
(1+y \omega^q) \sum_{m=0}^{N-2} (-y)^m \sum_{\beta =0}^m \omega^{\beta q}  
     &=& \sum_{m=0}^{N-2} (-y)^m 
         \sum_{\beta =0}^m \omega^{\beta q} - \sum_{m=0}^{N-2}(-y)^{m+1}\sum_{\beta =0}^m \omega^{(\beta +1)q} .
         \end{eqnarray}
The first sum on the right in eq.(C1) is
\begin{eqnarray}
                   \sum_{m=0}^{N-2}(-y)^{m} +\sum_{m=1}^{N-2} (-y)^m \sum_{\beta =1}^m \omega^{\beta q} ,
    \end{eqnarray}
while the 2nd sum, after relabeling $\beta$ and then $m$, is
   \begin{eqnarray}
- \sum_{m=1}^{N-1}(-y)^{m}\sum_{\beta =1}^m \omega^{\beta q}.
         \end{eqnarray}
The 2nd sum in expression (C2) will therefore cancel with all but the $m=N-1$ term in (C3), and we have
  \begin{eqnarray}
(1+y \omega^q) \sum_{m=0}^{N-2} (-y)^m \sum_{\beta =0}^m \omega^{\beta q}  
     &=&    \sum_{m=0}^{N-2}(-y)^{m} - (-y)^{N-1}\sum_{\beta =1}^{N-1} \omega^{\beta q}
     = \sum_{m=0}^{N-1} (-y)^m
     \end{eqnarray}
where we've used eq.(7) to evaluate the sum over $\beta$ as equal to -1.  The lemma is therefore valid for $p=1$.\\

We now assume the lemma is valid for  the product of  $(p-1)$ factors,  $\prod_{s=1}^{p-1} (1+y\omega^{q_s})^{-1},$ and express the product of $p$ factors as the difference of two
$(p-1)$ products:
\begin{eqnarray*}
&&\sum_{m=0}^{N-1} (-y)^m  \prod_{s=1}^p \frac{1} { 1+y\omega^{q_s} }\\
     &&~~~~~~ =  \frac{- 1}{y(\omega^{q_2} -\omega^{q_3})}  \sum_{m=0}^{N-1} (-y)^m    \frac{1}{1+y\omega^{q_1}}   
     \left\{        \frac{1}{1+y\omega^{q_2}} - \frac{1}{1+y\omega^{q_3}} \right\} \prod_{s=4}^p \frac{1}{1+y\omega^{q_s}} \\
     &&~~~~~~ = \frac{1}{\omega^{q_2} -\omega^{q_3}}   \sum_{m=1}^{N-p} (-y)^{m-1}\left\{   \sum_{\beta_{p}=0}^{m}    
                 \cdots \sum_{\beta_2=0}^{\beta_4}  \sum_{\beta_1=0}^{\beta_2}  \omega^{(\beta_{p}-\beta_{p-1})q_{p}}
                 \cdots  \omega^{(\beta_4-\beta_2)q_4}\omega^{(\beta_2-\beta_1)q_2} \omega^{\beta_1q_1}     \right.   \\               
     &&~~~~~~~~~~~~~~~~~~~~~~~~ - \left.   \sum_{\beta_{p}=0}^{m}   
                 \cdots \sum_{\beta_3=0}^{\beta_4} \sum_{\beta_1=0}^{\beta_3}  \omega^{(\beta_{p}-\beta_{p-1})q_{p}}
                 \cdots  \omega^{(\beta_4 -\beta_3)q_4}\omega^{(\beta_3-\beta_1)q_3} \omega^{\beta_1q_1}    \right\}.  
                 \end{eqnarray*}
In the first term inside the brackets we make the replacement $\beta_2 \rightarrow \beta_3$ and factor out the 
$q_p, \ldots , q_4$ sums, so we have

\begin{eqnarray*}
     &&\frac{1}{\omega^{q_2} -\omega^{q_3}}    \left\{ \sum_{\beta_{p}=0}^{m}    
                 \cdots \sum_{\beta_2=0}^{\beta_4}  \sum_{\beta_1=0}^{\beta_2}  \omega^{(\beta_{p}-\beta_{p-1})q_{p}}
                 \cdots  \omega^{(\beta_4-\beta_2)q_4}\omega^{(\beta_2-\beta_1)q_2} \omega^{\beta_1q_1}  \right\}    \\               
     &&~~~~~~~~~~~~~~~~~~~~~~~~ - \left. \sum_{\beta_{p}=0}^{m}   
                 \cdots \sum_{\beta_3=0}^{\beta_4} \sum_{\beta_1=0}^{\beta_3}  \omega^{(\beta_{p}-\beta_{p-1})q_{p}}
                 \cdots  \omega^{(\beta_4 -\beta_3)q_4}\omega^{(\beta_3-\beta_1)q_3} \omega^{\beta_1q_1} \right\} \\
      &&~~~~~~~~~~~~~ =  \sum_{\beta_{p}=0}^{m}    
                 \cdots \sum_{\beta_3=0}^{\beta_4}  \omega^{(\beta_{p}-\beta_{p-1})q_{p}}
                 \cdots  \omega^{(\beta_4-\beta_3)q_4} \sum_{\beta_1=0}^{\beta_3}  
            ~ \frac{  \omega^{(\beta_3-\beta_1)q_2}  - \omega^{(\beta_3-\beta_1)q_3} }{\omega^{q_2} -\omega^{q_3}}   ~  \omega^{\beta_1q_1} .
                 \end{eqnarray*}
 The validity of the lemma for $p$ factors then follows from the identity
  \begin{eqnarray*}
\frac{a^n -b^n}{a-b} = \sum_{p=0}^{n-1}~ a^{n-1-p}\, b^p,
\end{eqnarray*}
and a relabeling of the summation indices.   $\triangle$\\

\section{Proof of Lemma 2}
To remind ourselves, for a partition $Z=(z_1,z_2, \ldots , z_j)$ of $p$, the vector ${\bf q}_{_Z}$ is defined as 
\begin{eqnarray*}
{\bf q}_{_Z}  &=& (q_1, \ldots , q_1, q_2, \ldots , q_2, ~\ldots \ldots ~, q_j, \ldots , q_j) .\\
 && ~~~\,   \longleftrightarrow ~~~~~~ \longleftrightarrow~~~~~~~~~~~~~~~~~~\longleftrightarrow  \\
 &&~~~~~~z_1~~~~~\,~~~~~z_2~~~~~~~~~\cdots ~~~~~~~~~z_j
\end{eqnarray*}
and the statement of the lemma is that
\begin{eqnarray*}
     \sum_{1 \leq q_1< \cdots <  q_p \leq M}   h({\bf q})   
     &=& \frac{1}{p!} \sum_{Z \in {\cal P}_p}  (-1)^{p+j } ~  (p; k_1,k_2, \ldots , k_p)^*
      ~ \sum_{1 \leq q_1, \ldots , q_j \leq M} h( {\bf q}_{_Z})\\
\end{eqnarray*}
where  $h$ is a completely symmetric function of its arguments.\\

This proof is also by induction.   We separate off the $j=p$ term in the sum and write the lemma in the form
\begin{eqnarray}
  \sum_{ q_1, \ldots , q_p } h( q_1, \ldots , q_p)
         &=&  p!  \sum_{ q_1< \cdots <  q_p }   h(q_1, \ldots ,q_p) +  \sum_{j=1}^{p-1}~ \sum_{k_1,\ldots , k_p} (-1)^{p+j -1} ~  (p; k_1,k_2, \ldots , k_p)^* ~ \sum_{q_1, \ldots , q_j }  h( {\bf q}_{_Z})
   \end{eqnarray}
with $j=k_1 + \cdots +k_p$, and the summation limits $1 \leq q_n \leq M$ and $0 \leq k_n \leq p$ for all $q_n$ and $k_n$ understood.\\

 For $p=1$ the equality is trivial.  For $p=2$, this expression reduces to 
\begin{eqnarray}
    \sum_{q_1, q_2}  h( q_1,q_2) = 2 \sum_{q_1<q_2} h(q_1,q_2)  + \sum_{q_1}h(q_1,q_1),
   \end{eqnarray}
which is clearly valid.  We now assume that the lemma is valid for summation over the $(p-1)$-dimensional lattice and show that it is then valid for the $p$-dimensional one.  For clarity in the notation, we will denote a partition of $(p-1)$ by $Y = (y_1,y_2, \ldots, y_j )$, with multiplicities $[l_1, \ldots , l_{p-1}, l_p =0]$.  We then have for the summation over $p$ dimensions,
\begin{eqnarray}
\sum_{q_1, \ldots , q_p} h(q_1, \ldots , q_p) &=& \sum_{q_p}~ \sum_{q_1, \ldots, q_{p-1}} h(q_1, \ldots , q_{p-1}, q_p) \nonumber\\
    &=& \sum_{q_p} \left\{  (p-1)!  \sum_{q_1< \cdots <  q_{p-1} }   h(q_1, \ldots ,q_{p-1} ,q_p) \right. \\ 
   && \left.   + \sum_{j=1}^{p-2} ~\sum_{l_1, \ldots , l_{p-1} } (-1)^{p+j } ~  (p-1; l_1,l_2, \ldots , l_{p-1},0)^* ~ \sum_{ q_1, \ldots , q_j } 
    h( {\bf q}_{_Y},q_p)  \right\}  \nonumber
 \end{eqnarray}
 by the induction assumption.  We now separate the sum over $q_p$ applied to the first term inside the brackets into sums over the $p$ regions
 \begin{eqnarray*}
    ( q_p< q_1 <q_2 < \cdots < q_{p-1}) ~~(q_1 < q_p < q_2 < \cdots < q_{p-1}) ~~ \cdots~~  (  q_1< q_2< \cdots  <q_{p-1} < q_p)
 \end{eqnarray*}
  and the $(p-1)$ boundaries 
  \begin{eqnarray*}
 (q_p=q_1 < q_2 < \cdots < q_{p-1})~~ \cdots~~( q_1< q_2< \cdots < q_{p-1} = q_p) .
 \end{eqnarray*}
   By the symmetry of the function $h$, the $p$-dimensional sum over each subregion is the same, and we can write (making the change $p \rightarrow j+1, q_p \rightarrow q_{j+1}$ in the last sum),
 \begin{eqnarray}
   \sum_{q_1, \ldots , q_p} h(q_1, \ldots , q_p)  
        &=& p!   \sum_{q_1< \cdots <  q_{p-1} < q_p }   h(q_1, \ldots ,q_p)\\
        &&  +(p-1)!  \sum_{q_1< \cdots < q_{p-1}} \{ h(q_1, \ldots ,q_{p-1}, q_{1}) 
         + h(q_1, \ldots ,q_{p-1}, q_{2})  + \cdots + h(q_1, \ldots ,q_{p-1}, q_{p-1}) \}  \nonumber \\ 
        &&    + \sum_{j=1}^{p-2} ~\sum_{l_1, \ldots , l_{p-1} }  (-1)^{p+j } ~  (p-1; l_1,l_2, \ldots , l_{p-1},0)^*
    ~ \sum_{ q_1, \ldots , q_j, q_{j+1}}  h( {\bf q}_{_Y},q_{j+1})  .  \nonumber
    \end{eqnarray}
 Now although the functions $h(q_1,q_2,\ldots , q_{p-1}, q_1), \ldots h(q_1, q_2, \ldots , q_{p-1},q_{p-1})$ in the second term are not completely symmetric functions on the $(p-1)$-dimensional lattice, the sum of these functions {\it is }symmetric, and we can apply the lemma in its $(p-1)$ form to it:
 \begin{eqnarray*}
        &&  (p-1)!  \sum_{q_1< \cdots < q_{p-1}} \{ h(q_1, \ldots ,q_{p-1}, q_{1}) 
         + h(q_1, \ldots ,q_{p-1}, q_{2})  + \cdots + h(q_1, \ldots ,q_{p-1}, q_{p-1}) \}   \\ 
          &&~~~~~~~ =    \sum_{j=1}^{p-1} ~\sum_{l_1, \ldots , l_{p-1} }  (-1)^{p-1+j } ~  (p-1; l_1,l_2, \ldots , l_{p-1},0)^* 
        \sum_{ q_1, \ldots , q_j} \{ h {\bf q}_{_Y},q_{(1)}) + h({\bf q}_{_Y},q_{(2)}) + \cdots + h({\bf q}_{_Y},q_{(p-1)}) ~   \} .
\end{eqnarray*}
The notation $q_{(k)}$ here means, for a partition $Y$,  $q_{(k)} $ is equal to the value of the $k$-th coordinate
 in the associated lattice point ${\bf q}_{_Y}$.  E.g., for $Y=(3,2,2)$ and for $k=2$, we have
\begin{eqnarray*}
 && h(q_1,q_1,q_1, q_2,q_2, q_3,q_3, q_{(2)}) = h(q_1,q_1,q_1,q_2,q_2,q_3,q_3,q_1).\\
 &&~~~~~~  \uparrow~~~~~~~~~~~~~~~~~~~~~~\uparrow
\end{eqnarray*}
Therefore,
   \begin{eqnarray}
     \sum_{q_1, \ldots , q_p} h(q_1, \ldots , q_p)    
      &=& p!   \sum_{q_1< \cdots  < q_p }   h(q_1, \ldots ,q_p)  \nonumber    \\ 
      && + \sum_{j=1}^{p-1} ~\sum_{l_1, \ldots , l_{p-1} }  (-1)^{p-1+j } ~  (p-1; l_1,l_2, \ldots , l_{p-1},0)^*  \\
      &&  ~~~~~~~~~~ \times  \sum_{ q_1, \ldots , q_j} \{ h( {\bf q}_{_Y},q_{(1)}) + h({\bf q}_{_Y},q_{(2)}) + \cdots + h({\bf q}_{_Y},q_{(p-1)} )~   \} 
      \nonumber \\
     &&  + \sum_{j=1}^{p-2} ~\sum_{l_1, \ldots , l_{p-1} }(-1)^{p+j } ~  (p-1; l_1,l_2, \ldots , l_{p-1},0)^* ~
      \sum_{ q_1, \ldots , q_{j+1}}  h( {\bf q}_{_Y},q_{j+1}) . \nonumber  
\end{eqnarray}

The term on the left and the 1st term on the right are in the ``correct" form for the proposition.  The 2nd and 3rd terms on the right contain sums over the partitions of $(p-1)$ that, to complete the proof, we want to express as sums over the partitions of $p$.  Given a partition $(y_1, \ldots , y_j)$ of $(p-1)$, there are two ways a partition of $p$ can be generated from it:   one can add 1 to an individual part,  $y_m \rightarrow (y_m+1)$,  or one can add a 1 to the total sum $y_1+ \cdots + y_j \rightarrow y_1 + \cdots + y_j +1$.  As shown below, the 2nd term converts to a sum over the partitions of $p$ by this first process, and the 3rd term by the second process.\\

We consider the 2nd term first:
   \begin{eqnarray*}
   &&\sum_{j=1}^{p-1} ~\sum_{l_1, \ldots , l_{p-1} }  (-1)^{p-1+j } ~  (p-1; l_1, \ldots , l_{p-1},0)^* 
    \sum_{ q_1, \ldots , q_j} \sum_{s=1}^{p-1} ~ h( {\bf q}_{_Y},q_{(s)})  .
    \end{eqnarray*}
For a partition $Y=(y_1, \ldots , y_j)$ and for a particular $q_{(s)}$ such that
\begin{eqnarray}
y_1+ \cdots + y_{d-1} < s \leq y_1 + \cdots + y_{d-1} + y_d
\end{eqnarray}
for some $d$, the lattice point $({\bf q}_{_Y}, q_{(s)})$ is equivalent to the partition $Z \equiv Z_{y_d}= (y_1, \ldots, y_d+1, \ldots ,y_j)$ of $p$,
since $h({\bf q}_{_Y},q_{(s)}) = h(q_1, \ldots , q_1, q_2, \ldots , q_2, \ldots , q_j, \ldots , q_j, q_d)$, and, by the symmetry of $h$, we can move the $q_d$ in the final position to a position among the other $q_d$'s, so there are $(y_d+1)$ $q_d$'s.  In terms of multiplicities,  if
\begin{eqnarray*}
Y =[l_1, \ldots, l_{y_d}, l_{y_d+1}, \ldots , l_{p-1},0],
\end{eqnarray*}
then
\begin{eqnarray*}
Z_{y_d}= [ l_1, \ldots , l_{y_d}-1, l_{y_d+1}+1, \ldots , l_{p-1}, 0] .
\end{eqnarray*}
In the sum over $s$, the number of $s$'s that satisfy condition (D6) above is $y_d l_{y_d}$, so the partition $Z_{y_d}$ will appear $y_d l_{y_d} $ times in this sum, and we have,
   \begin{eqnarray}
   &&\sum_{j=1}^{p-1} ~\sum_{l_1, \ldots , l_{p-1} }  (-1)^{p-1+j } ~  (p-1; l_1, \ldots , l_{p-1},0)^* 
    \sum_{ q_1, \ldots , q_j} \sum_{s=1}^{p-1} ~  h( {\bf q}_{_Y},q_{(s)}) \nonumber  \\
    &&~~~~=   \sum_{j=1}^{p-1} ~\sum_{l_1, \ldots , l_{p-1} }(-1)^{p-1+j} (p-1; l_1, \ldots, l_{p-1},0)^*
     \sum_{n=1}^{p-1}~n l_n~ \sum_{q_1, \ldots, q_j} h({\bf q}_{_{Z_n}}) 
     \end{eqnarray}
where the partition $Z_n$ has multiplicities $[k_1, \ldots , k_p]$ given by $k_m= l_m - \delta_{n,m},~k_{m+1} = l_{m+1} + \delta_{n,m}$ for $m=1, \ldots , p-1$.  We now multiply and divide on the right by the corresponding multinomial coefficient $(p; k_1, \ldots , k_p)^*_n$:
   \begin{eqnarray}
   &&\sum_{j=1}^{p-1} ~\sum_{l_1, \ldots , l_{p-1} }  (-1)^{p-1+j } ~  (p-1; l_1, \ldots , l_{p-1},0)^* 
    \sum_{ q_1, \ldots , q_j} \sum_{s=1}^{p-1} ~  h( {\bf q}_{_Y},q_{(s)}) \nonumber \\
      &&~~~~=   \sum_{j=1}^{p-1} ~\sum_{l_1, \ldots , l_{p-1} } (-1)^{p-1+j}~ \sum_{n=1}^{p-1}
       \frac{(p-1; l_1, \ldots, l_{p-1},0)^*}{(p;k_1, \ldots , k_{p})^*_n}  ~(p; k_1, \ldots, k_{p})^*_n~n l_n~
        \sum_{q_1, \ldots, q_j} h({\bf q}_{_{X_n}}) .
      \end{eqnarray}
  Since
      \begin{eqnarray}
   \frac{(p-1; l_1, \ldots, l_{p-1},0)^*}{(p;k_1, \ldots , k_{p})^*_n}
             =  \frac{1}{p}~\frac{n^{l_n-1}(l_n-1)!~(n+1)^{l_{n+1}+1}(l_{n+1}+1)!}{n^{l_n}l_n! ~(n+1)^{l_{n+1}} l_{n+1}!} 
             = \frac{1}{p}~\frac{(n+1)(l_{n+1}+1)}{n l_n}   
     \end{eqnarray}
 we have    
       \begin{eqnarray}
   &&\sum_{j=1}^{p-1} ~\sum_{l_1, \ldots , l_{p-1} }  (-1)^{p-1+j } ~  (p-1; l_1, \ldots , l_{p-1},0)^* 
    \sum_{ q_1, \ldots , q_j} \sum_{s=1}^{p-1} ~  h( {\bf q}_Y,q_{(s)}) \nonumber \\
         &&~~~~=   \sum_{j=1}^{p-1} ~\sum_{l_1, \ldots , l_{p-1} } (-1)^{p-1+j}~ \sum_{n=1}^{p-1}
       \frac{(n+1)k_{n+1}}{p}  ~(p; k_1, \ldots, k_p)^*_n~ \sum_{q_1, \ldots, q_j} h({\bf q}_{_{X_n}}) .
      \end{eqnarray}
In this equation, for a given partition $Y$, we sum over all partitions $Z_n$.  We now interchange the sums and, for a given partition $Z$, we sum over the corresponding $Y$ partitions.   This amounts to dropping the subscripts $`n'$ and replacing the sums over the $l$'s by sums over $k$'s.  So, making the replacement $n' = n+1$ and then dropping the prime,
      \begin{eqnarray}
   &&\sum_{j=1}^{p-1} ~\sum_{l_1, \ldots , l_{p-1} }  (-1)^{p-1+j } ~  (p-1; l_1, \ldots , l_{p-1},0)^* 
    \sum_{ q_1, \ldots , q_j} \sum_{s=1}^{p-1} ~  h( {\bf q}_Y,q_{(s)}) \nonumber \\
         &&~~~~=   \sum_{j=1}^{p-1} ~\sum_{k_1, \ldots , k_{p} } (-1)^{p-1+j}~   \sum_{n=2}^{p} \frac{ nk_{n}}{p}~(p; k_1, \ldots, k_p)^*~ 
         \sum_{q_1, \ldots, q_j} h({\bf q}_{_{Z}})  \nonumber\\
         &&~~~~=   \sum_{j=1}^{p-1} ~\sum_{k_1, \ldots , k_{p} } (-1)^{p-1+j} ~   \frac{ p-k_1}{p}~(p; k_1, \ldots, k_p)^*~ 
         \sum_{q_1, \ldots, q_j} h({\bf q}_{_{Z}}) 
      \end{eqnarray}
where we've used $k_1+2k_2 + \cdots + p k_p =p$ in the last line.
~\\

Turning now to the 3rd term, 
  \begin{eqnarray*}
        \sum_{j=1}^{p-2} ~\sum_{l_1, \ldots , l_{p-1} }(-1)^{p+j } ~  (p-1; l_1, \ldots , l_{p-1},0)^* ~
      \sum_{ q_1, \ldots , q_{j+1}}  h( {\bf q}_{_Y},q_{j+1})   ,
\end{eqnarray*}
the lattice point $({\bf q}_{_Y},q_{j+1})$ is equal to the lattice point ${\bf q}_{_Z}$ for the partition $Z = (y_1, \ldots, y_j, 1)
=[l_1+1,l_2, \ldots , l_{p-1},0] \in {\cal P}(p)$ with $j+1$ parts.  Setting $j'=j+1$, we can express this term  as a sum over the partitions of $p$
from  $j' =2$ to $p-1$:
\begin{eqnarray}
&&   \sum_{j=1}^{p-2} ~\sum_{l_1, \ldots , l_{p-1}}  (-1)^{p+j } ~  (p-1; l_1, \ldots , l_{p-1},0)^* 
~ \sum_{ q_1, \ldots , q_{j+1}}  h( {\bf q}_{_Y},q_{j+1}) \nonumber   \\
  &&~~~~~~~~~~~~ = \sum_{j'=2}^{p-1} ~\sum_{l_1, \ldots , l_{p-1}}  (-1)^{p+j'-1 } ~\frac{(p-1; l_1, \ldots , l_{p-1},0)^*}{(p; l_1+1, \ldots , l_{p-1},0)^*}  ~ 
   (p;l_1+1, \ldots , l_{p-1},0)^* ~\sum_{ q_1, \ldots ,  q_{j'}}  h( {\bf q}_{_Z}) \nonumber \\
   &&~~~~~~~~~~~~ =    \sum_{j'=2}^{p-1} ~\sum_{k_1, \ldots , k_{p-1} } (-1)^{p+j'-1 }~  \frac{k_1}{p}  ~  (p;k_1, \ldots , k_{p-1},0)^* 
   ~\sum_{ q_1, \ldots ,  q_{j'}}  h( {\bf q}_{_Z})
   \end{eqnarray}
using
\begin{eqnarray*}
  \frac{(p-1; l_1,l_2, \ldots , l_{p-1})^*}{(p; l_1+1, \ldots , l_{p-1},0)^*}
  =\frac{1}{p} ~ \frac{1^{l_1+1}(l_1+1)!\, 2^{l_2}l_2!\cdots (p-1)^{l_{p-1}} l_{p-1}! }{1^{l_1}l_1!\, 2^{l_2}l_2! \cdots (p-1)^{l_{p-1}} l_{p-1}!}
  = \frac{l_1+1}{p}= \frac {k_1}{p} .
 \end{eqnarray*}
 But since $k_1=0$ if $j'=1$, (corresponding to the partition $p=p$ and $k_p=1$), the summation lower limit can be extended to $j'=1$ and 
 a ``sum" over $k_p$ added, along with the replacement $(p; k_1, \ldots , k_{p-1},0)^*  \rightarrow (p; k_1, \ldots , k_{p-1},k_p)^*$. \\
 
 Putting  (D11) and (D12) into (D5), the $k_1/p$ terms cancel and we get
 \begin{eqnarray*}
  \sum_{ q_1, \ldots , q_p } h( q_1, \ldots , q_p)
         &=&  p!  \sum_{ q_1< \cdots <  q_p }   h(q_1, \ldots ,q_p)  
         +  \sum_{j=1}^{p-1}~ \sum_{k_1,\ldots , k_p} (-1)^{p+j -1} ~  (p; k_1, \ldots , k_p)^* ~ \sum_{q_1, \ldots , q_j }  h ({\bf q}_{_Z}).~~~~\triangle
   \end{eqnarray*}

\section{Proof of Lemma 3}
 The multinomial coefficient $(p;k_1, \ldots , k_p)^*$ has the generating function \cite{Abram}
 \begin{eqnarray}
 &&  j!~\sum_{p=j}^{\infty} ~\frac{t^p}{p!} ~\sum_{k_1, \ldots, k_p} ~(p;k_1,k_2, \ldots , k_p)^*x_1^{k_1}x_2^{k_2} \cdots x_p^{k_p}
            =   \left( \sum_{i=1}^{\infty} ~\frac{x_i}{i}~t^i \right)^j \\
 &&~~ \nonumber
  \end{eqnarray}
 where the sums over the $k$'s are subject to the two constraints 
 \begin{eqnarray*}
 k_1+k_2 + \cdots + k_p &=& j, \\
 k_1+2k_2 + \cdots + p k_p &=& p.
 \end{eqnarray*}
 Dividing both sides by $j!$,  summing over $j$, and interchanging the $j$ and the $p$ summations on the left, this becomes 
 \begin{eqnarray}
 \sum_{p=0}^{\infty} ~\frac{t^p}{p!}~ \sum_{j=0}^{p} ~\sum_{k_1, \ldots , k_p} ~(p;k_1,k_2, \ldots , k_p)^*x_1^{k_1}x_2^{k_2} \cdots x_p^{k_p} 
                     = \sum_{j=0}^{\infty} \frac{1}{j!}~ \left( \sum_{i=1}^{\infty} ~\frac{x_i}{i}~t^i \right)^{j }
                     = \exp \left( \sum_{i=1}^{\infty} ~\frac{x_i}{i}~t^i \right) .
\end{eqnarray}
Now setting all $x_i$ to 1, we get
\begin{eqnarray}
  \sum_{p=0}^{\infty}~\frac{t^p}{p!} \sum_{j=0}^p ~\sum_{k_1, \ldots ,k_p} ~(p;k_1,k_2, \ldots , k_p)^* = \exp \left(  \sum_{i=1}^{\infty} ~\frac{t^i}{i} \right)&=& \exp \left( - \ln (1-t) \right)= \frac{1} {1-t} = \sum_{n=0}^{\infty} t^{n}.
\end{eqnarray}
 Comparing the two,  
 \begin{eqnarray}
 \sum_{j=0}^p ~\sum_{k_1, \ldots , k_p} ~(p;k_1,k_2, \ldots , k_p)^* = p!.
 \end{eqnarray}

Now, for a set on nonnegative integers, $\{ \beta_1, \ldots , \beta_q~|~ \beta_i \geq 0 \} $, we  apply the differential operator 
 \begin{eqnarray*}
 \frac{ \partial^{\beta_1 +\beta_2 + \cdots + \beta_q}} { \partial^{\beta_1} x_1\cdots \partial^{\beta_q} x_{q} }
  \end{eqnarray*}
to both sides of (E2) and again set all $x_i$ equal to 1.  On the right we get
 \begin{eqnarray}
&& \left. \frac{ \partial^{\beta_1 +\beta_2 + \cdots + \beta_q}} { \partial^{\beta_1} x_1\cdots \partial^{\beta_q} x_{q} }
~ \exp \left(  \sum_{i=1}^{\infty} ~\frac{x_i}{i}~t^i \right)\right|_{x_i =1}~
= \frac{t^{\beta_1 +2 \beta_2 + \cdots + q \beta_q}}{1^{\beta_1} \cdots q^{\beta_q}}  \sum_{n=0}^{\infty} ~ t^n  
=  \frac{1}{1^{\beta_1} \cdots q^{\beta_q} } \sum_{n=n_0}^{\infty} ~ t^n,\\
&& {\rm with ~} n_0= \beta_1 +2 \beta_2 + \cdots + q\beta_q , \nonumber
\end{eqnarray}
and on the left,
\begin{eqnarray}
&& \left. \frac{ \partial^{\beta_1 +\beta_2 + \cdots + \beta_q}} { \partial^{\beta_1} x_1\cdots \partial^{\beta_q} x_{q} } ~
\sum_{p=0}^{\infty} ~\frac{t^p}{p!} \sum_{j=0}^p~\sum_{k_1, \ldots , k_p} ~(p;k_1,k_2, \ldots , k_p)^* x_1^{k_1}x_2^{k_2}  \cdots  x_p^{k_p} \right|_{x_i=1}
\nonumber \\
&&~~~~~~~~~~~~~~~~~~ = \sum_{p=n_0}^{\infty} ~\frac{t^p}{p!} \sum_{j=0}^p~\sum_{k_1,\ldots , k_p} (p;k_1,k_2, \ldots , k_p)^*
~k_1(k_1-1) \cdots (k_1-\beta_1+1) \nonumber \\ 
&&~~~~~~~~~~~~~~~~~~~~~~~~~~~~~~~~~~~~~~~~~~~~~~~~~~ \times  k_2(k_2-1)\cdots (k_2-\beta_2+1) 
                                \cdots k_q(k_q-1) \cdots (k_q-\beta_q+1) ~ \nonumber\\
&&~~~~~~~~~~~~~~~~~~ = \sum_{p=n_0}^{\infty} ~\frac{t^p}{p!} \sum_{j=0}^p~\sum_{k_1, \ldots , k_p}  
        \beta_1! \left( \begin{array} {c} k_1 \\~\\ \beta_1 \end{array} \right) \cdots  \beta_q! \left( \begin{array} {c} k_q \\~\\ \beta_q \end{array} \right)  
         ~(p;k_1,k_2, \ldots , k_p)^* .
\end{eqnarray}
Therefore
\begin{eqnarray}
 && \sum_{j=0}^p~\sum_{k_1, \ldots , k_p} \left( \begin{array} {c} k_1 \\~\\ \beta_1 \end{array} \right) \cdots 
          \left( \begin{array} {c} k_p \\~\\ \beta_p \end{array} \right)    ~(p;k_1,k_2, \ldots , k_p)^*= 
           \frac{p!}{1^{\beta_1} \beta_1!~\, 2^{\beta_2} \beta_2! \, \cdots \, p^{\beta_p}\beta_p!}
  \end{eqnarray}
  for $ p \geq \beta_1 +2 \beta_2 + \cdots + p \beta_p. ~~\triangle  $

\section{Sample calculations}

As examples, we will use Theorem 3 to calculate the $N=7$ coefficient $C_{0123456}$ and the $N=8$ coefficient $C_{00224466}$.

\subsection{$C_{0123456}$}
   Since the $A_k$'s are all distinct in this example, the set of partitions of the set $\{A_1,A_2,A_3,A_4\}$, i.e., ${\cal P}\{A_1A_2A_3A_4\}$, and the set of partitions of the set of indices of the $A_k$'s, ${\cal P}\{1234\}$, have the same number of elements.  We will therefore use ${\cal P}\{1234\}$ and the somewhat simpler expression, (10d), to find $C_{0123456}$:  
  \begin{eqnarray*}
C_{0123456} &=&  - 7 \times  \left\{~ 5! +   \sum_{\Theta \in {\cal P} \{1234  \} }  
                      ~  \prod_{k=1}^j ~ (z_k-1)! ~ \sum_{\mu=1}^j  (-N)^{\mu}  \sum_{\lambda_1, \ldots , \lambda_j =0,1 }
                       \delta_{\mu, \lambda_1 + \cdots + \lambda_j } ~ H \left[ 1 - \sum_{t=1}^j \lambda_t X(\theta_t) \right]  \right.   \\ 
    &&\left.  ~~~  ~~~~~  \times   \left( \begin{array} {c} 5 - \sum_{t=1}^j \lambda_t (X(\theta_t) +z_t ) \\
            ~\\   1 - \sum_{t=1}^j \lambda_t X(\theta_t)  \end{array} \right) 
                \prod_{s=1}^j    \left( \begin{array} {c} X(\theta_s) + z_s -1\\~\\ z_s-1 \end{array} \right)^{\lambda_s} \right\}.\\
               && ~\\
              {\cal P}\{1234\} &=& \left\{ \frac{}{} (1234),~(123)(4),~ (124)(3),~(134)(2),~(234)(1),~(12)(34),~(13)(24),~(14)(23), \right. \\
  && \left. \frac{}{}  (12)(3)(4),~(13)(2)(4),~  (14)(2)(3),~(23)(1)(4),~(24)(1)(3),~(34)(1)(2),~(1)(2)(3)(4) \right\}  .
     \end{eqnarray*}   
   The value of, for example, $X(\theta = 1234)$ is $X(1234) \equiv -( A_1+A_2+A_3+A_4)  \mod 7 \equiv -(2+3+4+5)= 0 $.
   The rest of the $X(\theta)$'s are:
   \begin{eqnarray*}
  \begin{array} {cccccc} X(123) =5, & X(124)=4, & X(134) = 3, &X(234)= 2,\\
                            X(12) = 2, & X(34)=5, & X(13)=1, & X(24)= 6, & X(14)=0, & X(23)=0,\\
                            X(1)=5, & X(2)=4, & X(3)=3, & X(4)=2. \end{array}
  \end{eqnarray*}
   The subset of $ {\cal P}\{1234\}$ of partitions that correspond to non-zero terms, for which  $1 - \sum \lambda_t X(\theta_t) \geq 0 $ for at least some values of the $\lambda_t$'s, is then:
    \begin{eqnarray*}
   \left\{ \frac{}{} (1234),~(13)(24),~(14)(23),  ~(13)(2)(4),~  (14)(2)(3),~(23)(1)(4) \right\}.
\end{eqnarray*}
~\\~\\~\\~\\~\\~\\~\\~\\~\\~\\~\\~\\~\\~\\~\\~\\~\\~\\~\\~\\
So
  \begin{eqnarray*}
C_{0123456} 
                &=& -7  \left\{ ~  5!   +   (-7)^{ \lambda_1 } 3!\left.
                \left( \begin{array} {c} 5 - 4\lambda_1     \\~\\ 1    \end{array} \right) 
                    \left( \begin{array} {c} 3\\~\\ 3 \end{array} \right)^{\lambda_1} 
                     \right|^{\lambda_1= 1}_{\Theta = (1234) }\right.  \\
              &&   +  \left. (-7)^{ \lambda_1 + \lambda_2} 1!1!
                \left( \begin{array} {c} 5 - 3\lambda_1 - 8\lambda_2    \\~\\ 1 -  \lambda_1 - 6\lambda_2   \end{array} \right) 
                    \left( \begin{array} {c} 2\\~\\ 1 \end{array} \right)^{\lambda_1} 
                     \left( \begin{array} {c} 7\\~\\ 1 \end{array} \right)^{\lambda_2} \right|^{(\lambda_1,\lambda_2)= (1,0)}_{\Theta = (13)(24) } \\
             && +  \left.  (-7)^{ \lambda_1 + \lambda_2} 1!1!
                \left( \begin{array} {c} 5 - 2\lambda_1 - 2\lambda_2    \\~\\ 1 -0\lambda_1-0 \lambda_2   \end{array} \right) 
                    \left( \begin{array} {c} 1\\~\\ 1 \end{array} \right)^{\lambda_1} 
                     \left( \begin{array} {c} 1\\~\\ 1 \end{array} \right)^{\lambda_2}
                        \right|^{(\lambda_1,\lambda_2)= (1,0),(0,1),(1,1)}_{\Theta = (14)(23)}~\\
                     ~\\
                && \left. +   (-7)^{ \lambda_1 + \lambda_2+\lambda_3} 1!0!0!
                \left( \begin{array} {c} 5 - 3\lambda_1 - 5\lambda_2 -3\lambda_3     \\~\\ 1-\lambda_1  -4\lambda_2 -2\lambda_3  \end{array} \right) 
                    \left( \begin{array} {c} 2\\~\\ 1 \end{array} \right)^{\lambda_1} 
                     \left( \begin{array} {c} 4\\~\\ 0 \end{array} \right)^{\lambda_2}
                  \left( \begin{array} {c} 2\\~\\ 0 \end{array} \right)^{\lambda_3}
                  \right| ^{(\lambda_1,\lambda_2,\lambda_3)= (1,0,0)}_{\Theta = (13)(2)(4) }    \\       
                  ~\\    
            && +  \left.(-7)^{ \lambda_1 + \lambda_2+\lambda_3} 1!0!0!
                \left( \begin{array} {c} 5 - 2\lambda_1 - 5\lambda_2 -4\lambda_3    \\~\\ 1 -0 \lambda_1 -4\lambda_2 -3\lambda_3  \end{array} \right) 
                    \left( \begin{array} {c} 1\\~\\ 1 \end{array} \right)^{\lambda_1} 
                     \left( \begin{array} {c} 4\\~\\ 0 \end{array} \right)^{\lambda_2}
                  \left( \begin{array} {c} 3\\~\\ 0 \end{array} \right)^{\lambda_3} 
                  \right|^{(\lambda_1,\lambda_2,\lambda_3)= (1,0,0)}_{\Theta = (14)(2) (3) }  \\                      
               && \left. \left. +  (-7)^{ \lambda_1 + \lambda_2 +\lambda_3} 1!0!0!
                \left( \begin{array} {c} 5 - 2\lambda_1 - 6\lambda_2 -3\lambda_3    \\~\\ 1-0 \lambda_1 -5 \lambda_2 -2 \lambda_3 
                   \end{array} \right) 
                    \left( \begin{array} {c} 1\\~\\ 1 \end{array} \right)^{\lambda_1} 
                     \left( \begin{array} {c} 5\\~\\ 0 \end{array} \right)^{\lambda_2} 
                   \left( \begin{array} {c} 2\\~\\ 0 \end{array} \right)^{\lambda_3} 
                   \right|^{(\lambda_1,\lambda_2, \lambda_3))= (1,0)}_{\Theta = (23)(1)(4) }~ \right\} 
                   \end{eqnarray*}
 where for clarity we've indicated the particular partition that each term corresponds to and the allowed 
 values of the $\lambda_t$'s.  We have then
                   \begin{eqnarray*}
          C_{0123456} = &=& -7  \left\{   5!  -7 \cdot 3!  \left( \begin{array} {c} 1   \\~\\ 1   \end{array} \right) \left( \begin{array} {c} 3   \\~\\ 3   \end{array} \right)
      -7  \left( \begin{array} {c} 2     \\~\\ 0    \end{array} \right) 
                    \left( \begin{array} {c} 2\\~\\ 1 \end{array} \right) \right.  -7  \left( \begin{array} {c} 3  \\~\\ 1   \end{array} \right) 
                    \left( \begin{array} {c} 1   \\~\\ 1   \end{array} \right)
            -7   \left( \begin{array} {c} 3    \\~\\ 1   \end{array} \right)  \left( \begin{array} {c} 1   \\~\\ 1   \end{array} \right) \\
              && \left. ~~~~~~ +( -7) ^2\left( \begin{array} {c} 1   \\~\\ 1   \end{array} \right)  
             \left( \begin{array} {c} 1   \\~\\ 1   \end{array} \right)   \left( \begin{array} {c} 1   \\~\\ 1   \end{array} \right)
             -7 \left(\begin{array}{c} 2 \\~ \\ 0 \end{array} \right)  \left( \begin{array} {c} 2    \\~\\ 1    \end{array} \right)  
             -7   \left( \begin{array} {c} 3   \\~\\ 1   \end{array} \right)  \left( \begin{array} {c} 1   \\~\\ 1   \end{array} \right)
             -7 \left( \begin{array} {c} 3    \\~\\ 1    \end{array} \right)  \left( \begin{array} {c} 1   \\~\\ 1   \end{array} \right) \right\}
                 \\
                  &=&  -105.
                   \end{eqnarray*}    
                                 
 \subsection{$C_{00224466}$}
  In this example $M_1=0$ and the non-zero terms are characterized by $\sum \lambda_t X(\theta_t) = 0$.  So only partitions
  $\Theta = (\theta_1) \cdots (\theta_j)$ that contain at least one ``null'' part $\theta_k$ such that $X(\theta_k) =0$, will contribute; for a
  non-null part $\theta_q$,  $\lambda_q$ must equal zero for the term to contribute to the coefficient.  When $M_1=0$, the two binomial coefficients in the formula,
      \begin{eqnarray*}
         \left( \begin{array} {c} N-M_0-1 -\sum \lambda_t (X(\theta_t) + z_t) \\~\\M_1 - \sum \lambda_t X(\theta_t)   \end{array} \right) 
         ~~{\rm and~~}            \left( \begin{array} {c} X(\theta_s) + z_s -1\\~\\ z_s-1 \end{array} \right)^{\lambda_s},
        \end{eqnarray*}
  can both be replaced by 1.  The first binomial coefficient equals 1 since, for non-zero terms, it has zero in it's lower entry; in the second coefficient, either $X(\theta_s)=0$ or $\lambda_s=0$.  So we have
      \begin{eqnarray*}
C_{00224466} &=&  -  8 \left\{~ \frac{5!}{2!2!2!} +   \frac{1}{2} \sum_{\Theta \in {\cal P} \{22446 \} }  
                       \prod_{\theta }     \frac{1}{\kappa_{\theta}!}  ~  \prod_{k=1}^j ~\frac{   (z_k-1)! } {m_{2}^{(\theta_k)}!
                        m_4^{(\theta_k)}! m_{6}^{(\theta_k)}!  }  \right. \\
   && ~~~~ \left.  \times  ~ \sum_{\mu=1}^j  (-8)^{\mu}  \sum_{\lambda_1, \ldots , \lambda_j =0,1 } 
   \delta_{\mu, \lambda_1 + \cdots + \lambda_j }   ~ H \left[ - \sum_{t=1}^j \lambda_t X(\theta_t) \right]  \right\}  .
   \end{eqnarray*}
   The   null subsets of \{2,2,4,4,6\} are (2446), (224), (44), and (26).  The 10 contributing partitions are then:
   \begin{eqnarray*}
(2446)(2):&&~~(\lambda_1,\lambda_2)= (1,0)\\
(224)(46):&&~~(\lambda_1,\lambda_2)= (1,0)\\
(226)(44):&&~~(\lambda_1,\lambda_2)= (0,1)\ \\
(244)(26):&&~~(\lambda_1,\lambda_2)= (0,1)\\
 (224)(4)(6):&&~~(\lambda_1,\lambda_2, \lambda_3)= (1,0,0)  \\
(44)(26)(2):&&~~(\lambda_1,\lambda_2,\lambda_3)= (1,0,0),(0,1,0), (1,1,0)\\
(44)(22)(6):&&~~ (\lambda_1,\lambda_2,\lambda_3)=(1,0,0)\\
(26)(24)(4):&&~~ (\lambda_1,\lambda_2,\lambda_3)=(1,0,0)\\
(44)(2)(2)(6):&&~~  (\lambda_1,\lambda_2,\lambda_3,\lambda_4)=(1,0,0,0)\\
 (26)(2)(4)(4):&&~~ (\lambda_1,\lambda_2,\lambda_3,\lambda_4)=(1,0,0,0)
 \end{eqnarray*}
    (To simplify the expression below, we will omit the factor $\prod_{\theta }     \frac{1}{\kappa_{\theta}!} $ when all the $\kappa$'s equal 0 or 1.)  Summing over these partitions, we get
   \begin{eqnarray*}
   C_{00224466} &=&  -  8 \left\{~ 15 -   \frac{8}{2}  \left[  \frac{(4-1)!(~1-1)!}{1!2!1!~1!0!0!}  +   \frac{(3-1)!~(2-1)!}{2!1!0!~0!1!1!}
     +   \frac{(3-1)!~(2-1)!}{2!0!1!~0!2!0!}  +   \frac{(3-1)!~(2-1)!}{1!2!0!!~1!0!1!}  \right. \right. \\
     &&~~~~~~~~~~~~~~~~~~~~~~~~~~~~~~~~~~~~~~~~ +  \frac{(3-1)!~(1-1)!~(1-1)!}{2!1!0!~0!1!0~0!0!1!}  \\
     &&  +  2 \times \frac{(2-1)!~(2-1)!(1-1)!}{0!2!0!~1!0!1!~1!0!0!} 
    + \frac{(2-1)!~(2-1)!(1-1)!}{0!2!0!~2!0!0!~0!0!1!} +   \frac{(2-1)!~(2-1)!(1-1)!}{1!0!1!~1!1!0!~0!1!0!} \\
    &&  \left.  + \frac{1}{1!2!1!}  ~ \frac{(2-1)!~(1-1)!~(1-1)!~(1-1)!}{0!2!0! ~1!0!0!~1!0!0!~0!0!1! }  + 
     \frac{1}{1!1!2!} ~ \frac{(2-1)!~(1-1)!~(1-1)!~(1-1)!}{1!0!1! ~1!0!0!~0!1!0!~0!1!0! }  \right] \\
     && \left.  ~~~~~~~~~~~~~~~~~~~~~~~~~~~~~~~~~~~~~~~~ + \frac{64}{2} \frac{(2-1)!(2-1)!(1-1)!}{0!2!0!~1!0!1!~1!0!0!} ~ \right\}  \\
       &=&  -  120  +96  +32 +16 +  32 +32  +32 +8 +32 + 8 + 16- 128 =56. 
             \end{eqnarray*}
                  
 For this coefficient, the W-J method is actually more efficient.  The null subsets of \{2,2,4,4,6,6\} are 
(224466), (2266), (2446), (466), (224), (44), and (26).  There are then  5 null multisets to be summed over:
(224466), (2266)(44), (2446)(26), (224)(466), (26)(26)(44).  We have then, using Wyn-jones' formula,
 \begin{eqnarray*}
 C_{00224466} &=& (-8) \frac{(6-1)!}{2!2!2!} +(-8)^2 \frac{(4-1)!(2-1)!}{2!2!~2!} +(-8)^2\frac{(4-1)!(2-1)!}{2!} +(-8)^2 \frac{(3-1)!(3-1)!}{2!~2!}\\
 && +(-8)^3 \frac{1}{2!}~\frac{(2-1)!(2-1)!(2-1)!}{2!}\\
 &=& -120 +48 + 192 + 64 -128 = 56.
 \end{eqnarray*}

   \section{ The coefficients for the $N=6,7,8$ determinants  }   
      
 In this section we give our results for the coefficients $ [~ \{ C^*_{M_0 \cdots  M_{N-1}} \} ~] $ for $N$=6, 7, 8, using Theorem 3.
 The results are expressed in the form of a table.  In the first column are listed the allowed partitions of $N$; the second and third columns list, respectively, the form of the coefficient $C_{a_0 \cdots a_{N-1}}$  and all the additive multiplets associated with this partition; and the forth column lists the values of $[~ \{ C^*_{M_0 \cdots  M_{N-1}} \} ~] $.\\
 
 The subscript $n$ on $\{ C^*_{M_0 \cdots M_{N-1}}\}_n $ is used to indicate the number of elements in a particular additive multiplet.  To simplify the notation however it is understood that if there is no subscript then the multiplet has $N$ elements.\\
 
 Our results are:
 
 \begin{eqnarray*}
 N=6:&& (12 ~{\rm super-multiplets}).~~~~~~~~~~~~~[~ \{ C^*_{M_0 \cdots  M_5} \} ~]  \\
6&& C_{000000}   ~~~\,~~~   ~\{C^*_{600000}\}~~~~~~~~~~~~~~~~~~~~~~~~1 \\
42&& C_{0000aa}   ~~~\,~~~   ~ \{ C^*_{400200} \} ~ ~~~~~~~~~~~~~~~~~~~ -3 \\
411&& C_{0000ab} ~~ \left\{ \begin{array} {l} ~ ~\{C^*_{410001}\}~   \\~ ~ \{C^*_{401010}\}~   \end{array} \right. 
                                                    ~~~~~~~\,~~~~~  \left\{ \begin{array} {r} ~ -6 \\ ~-6  \end{array} \right. \\
33&& C_{000aaa}   ~~~\,~~~ ~\{C^*_{303000}\}~~~~~~~~~~~~~~\,~~~~~~~~~2 \\ 
321&& C_{000aab} ~~~\,~~~~\{C^*_{320010}\} ,~\{C^*_{301002}\} ~~~~~~~~~6 \\
3111&& C_{000abc} ~~~\,~~~~\{ C^*_{311100}\},~\{C^*_{300111}\}~ ~~~~~~~ 12  \\
222&& C_{00aabb} ~~  \left\{ \begin{array} {l}~~ \{C^*_{222000 }\}~ \\
                                                                               ~~ \{ C^*_{202020}\}_{ n=2}   \end{array} \right. 
                                                                                  ~~~~~~~~ ~  \left\{ \begin{array} {r}~ -9 \\ ~9   \end{array} \right. \\
2211&& C_{00aabc} ~~   \left\{ \begin{array} {l } ~~\{C^*_{210120}\}~ \\ 
                                                                                ~~  \{ C^*_{210201}\}~  \end{array} \right.
                                                                        ~~~~~~~~~~~\,~ \left\{ \begin{array} {r} \,\, -18 \\ ~ \,0   \end{array} \right.   \\
21111&& C_{00abcd} ~~~\,~~~~   \{C^*_{211011}\}~~~~~~~~~~~~~~~~~~~~~~~  ~0 
\end{eqnarray*}

\begin{eqnarray*}
  N=7 :&& (12 ~{\rm  super-multiplets}) ~~~~~~~\,~~~~~~~~~~~~~~~~~~~
   ~~~~~~~~~~~~~~~~~~~~~~~~~~~~~~~~~~~~~~~~~~~~~~~~~[~ \{ C^*_{M_0 \cdots  M_6} \} ~] \\
  7&& C_{0000000} ~~~~~  \{ C^*_{7000000} \} ~~~~~~~~~~~
  ~~~~~~~~~~~~~~~~~~~~~~~~~~~~~~~~~~~~~~~~~~~~~~~~~~~~~~~~\,~~~~~~~~~~~~~~~~~~1 \\
511&& C_{00000ab}~~~~~ \{ C^*_{5100001}\}, ~\{ C^*_{5010010}\} ,~ \{C^*_{5001100}\}   
~~~~~~~~~~~~~~~~~~~~~~~~~~~~~~~~~~~~~~~~~~~~~~~~~~  -7 \\
421&& C_{0000aab}~~~~~ \{C^*_{4200010}\}, ~\{C^*_{4021000}\}  ,~  \{  C^*_{4102000}\}, 
                                                          ~\{C^*_{4000201}\},~\{C^*_{4000120\}} ,   ~\{C^*_{4010002}\}  ~~\,~~~7   \\
4111&& C_{0000abc}~~~~~ \{C^*_{4110100}\},~\{C^*_{4001011}\} 
~~~~~~~~~~~~~~~~~~~~~~~~~~~~~~~~~~~~~~~~~~~~~~~~~~~~~~~~~~~~~\,~~~~~~~  14  \\
331&& C_{000aaab} ~~~~ ~ \{C^*_{3300100} \}, ~\{ C^*_{3130000}\}, ~\{C^*_{3003010}\} 
    ~~~~~~~~~~~~~~~~~~~~~~~~~~~~~~~~~~~~~~~~~~~~~~~~~~   -7\\
322&& C_{000aabb}~~~~~    \{ C^*_{3200002}\}, ~ \{ C^*_{3020020}\} , ~\{ C^*_{3002200}\} 
~~~~~~~~~~~~~~~~~~~~~~~~~~~~~~\,~~~~~~~~~~~~~~~~~~~~~~  14  \\
3211&& C_{000aabc} ~~~~~ \{C^*_{3211000}\}, ~\{C^*_{3020101}\}, ~\{ C^*_{3012001}\},
                                                       ~\{C^*_{3100210}\},~ \{ C^*_{3101020}\},~\{C^*_{3000112}\}   ~  -21 \\
31111&& C_{000abcd}~~~~~  \{C^*_{3110011}\},~\{C^*_{3011110}\},~\{ C^*_{3101101}\} 
~~~~~~~~~~~~~~~~~~~~~~~~~~~~~~~~~~~~~~~~~~~~~~~~~  ~ ~~~~7\\
2221&& C_{00aabbc} ~~~~~ \{C^*_{1002022}\}, ~\{C^*_{1220200}\}   
~~~~~~~~~~~~~~~~~~~~~~~~~~~~~~~~~~~~~~~~~~~~~~~~~~~~~~~~~~~~~~\,~~~  ~ -7  \\
22111&& C_{00aabcd}  ~ \left\{ \begin{array} {l}   \{ C^*_{1012210 } \},~ \{C^*_{1201102}\},~\{C^*_{1120021} \} 
                                             \\  \{ C^*_{1102201}\},~\{C^*_{1210012}\}, ~\{C^*_{1021120 }\}  \end{array} \right.
    ~~~~~~~~~~~~~~~~~~~~~~~~~~~~~~~\,~~~~~~~~~~~~~~      \left\{\begin{array} {c} -14 \\ ~~ 35 \end{array}  \right.  \\
1111111&& C_{0abcdef} ~~~~~  \{C^*_{1111111}\}_{n=1}
 ~~~~~~~~~~~~~~~~~~~~~~~~~~~~~~~~~~~~~~\,\,~~~~~~~~~~~~~~~~~~~~~~~~~~~~~~~~~~~~  -105 
 \end{eqnarray*}

~\\~\\~\\~\\~\\~\\~\\~\\~\\~\\~\\~\\~\\~\\~\\~\\

\begin{eqnarray*}
 N=8 :&& (49~ {\rm super-multplets})~~~~~\,~~~~~~~~~~~~~~~~~~~
   ~~~~~~~~~~~~~~~~~~~~~~~~~~~~~~~~~~~[~ \{ C^*_{M_0 \cdots  M_7} \} ~] \\
8&& C_{00000000}  ~~~~~~\{C^*_{80000000}\}
~~~~~~~~~~~~~~~~~~~~~~~~~~~~~\,~~~~~~~~~~~~~~~~~~~~~~~~~~~~~~~~~~ ~~1 \\
62&& C_{000000aa}  ~~~~\, ~\{C^*_{60002000}\}
~~~~~~~~~~~~~~~~~~~~~~~~~~~~~~~~~~~~~~~~~~~~~~~~~~~~~~~~~~~~~~ -4  \\
611&&C_{000000ab}~~\left\{ \begin{array} {l} \{ C^*_{61000001} \}, ~\{ C^*_{60010100}\} \\
                                                                                        \{ C^*_{60100010} \}~
                                                 \end{array} \right. 
         ~~~~~~~~~~~~~~~~~~~~~~~~~~~~~~\,~~~~~~~~ 
          \left\{ \begin{array} {r} -8  \\~ -8 \end{array} \right.\\
521&&C_{00000aab}~~ \left\{ \begin{array} {l} \{C^*_{52000010}\}, ~\{C^*_{50120000}\} ,~ \{C^*_{50000210}\},
                                                                                           ~\{C^*_{50100002}\}  \\
          \{C^*_{50201000}\} ,~\{C^*_{50001020}\}  \end{array} \right. ~~~~   \left\{ \begin{array} {r} 8
          \\ ~~~\,~8  \end{array} \right. \\
5111&&C_{00000abc} ~\,~ \left\{ \begin{array} {l} \{C^*_{51100100}\},~\{C^*_{50010011}\}~\\
                                                    \{ C^*_{51011000}\},~\{C^*_{50001101}\} \end{array} \right. 
~~~~~~~~~~~~~~~~~~~~~~~~~~\,~~~~~~~~~~~   \left\{ \begin{array} {r} ~~16 \\~~~16
 \end{array} \right. \\
44&&C_{0000aaaa} ~~ \left\{ \begin{array} {l} \{ C^*_{40400000} \} \\ 
                                                                                      \{ C^*_{40004000}\}_{n=4} \\  \end{array} \right. 
        ~~~~~~~~~~~~~~~~~~~~~~~~~~~~~~~~~~~~~~~~~~~~~~~~~~~  \left\{ \begin{array} {r}  -2 \\
         ~~~~ 6 \end{array} \right. \\
431&&C_{0000aaab}~~~~~~  \{ C^*_{43000100}\},~\{C^*_{40030001} \},~ \{ C^*_{41000300}\},~\{ C^*_{40010003}\} 
                                                                     ~~~~~~~~~~  -8      \\
422&&C_{0000aabb} ~~ \left\{ \begin{array} {l} \{ C^*_{42020000}\},~  \{ C^*_{40000202}\} \\
                                                                                    \{C^*_{40200020}\} \\ 
                                                                                    \{ C^*_{42000002}\} ,~\{C^*_{40020200}\} \end{array} \right. 
           ~~~~~~~~~~~~~~~~~~~~~~~~~~~~~\, ~~~~~~~~~    \left\{ \begin{array} {r} -12\\ 20  
            \\  20  \end{array} \right.  \\
4211&&C_{0000aabc} ~~  \left\{ \begin{array} {l}  \{C^*_{42101000}\},~\{ C^*_{40021010}\}, ~\{ C^*_{40101200}\},
                                                                   ~\{C^*_{40001012}\}\\ 
                                                  \{ C^*_{41210000} \},~\{ C^*_{41010020} \}, ~\{ C^*_{40200101} \}, ~\{C^*_{40000121} \} \\
                                                         \{  C^*_{41002001} \}, ~ \{ C^*_{40012100} \}  \\
                                                              \{ C^*_{40102010} \} \end{array} \right. 
                                                          ~~~~  \left\{ \begin{array} {r} -24 \\ -24 \\ 8 \\ 8
                                                           \end{array} \right.  \\~ \\
41111&&C_{0000abcd} ~~ \left\{ \begin{array} {l} \{ C^*_{41100011} \},~\{C^*_{40110110} \} \\
                                                                    \{ C^*_{41010101} \} \\
                                                                    \{ C^*_{41001110}\},~\{ C^*_{40111001} \}  \end{array} \right. 
~~~~~~~~~~~~~~~~~~~~~~~~~\,~~~~~~~~~~~~~  \left\{ \begin{array} {r} 16 \\ 16\\-48 \end{array} \right.  \\~\\
332&& C_{000aaabb} ~~  \left\{ \begin{array} {l} \{C^*_{23000003}\},~\{C^*_{20030300 }\}\\ 
                                                                                          \{C^*_{20300030}\}
                                                                                            \end{array} \right. 
                 ~~~~~~~~~~~~~~~~~~~~~~~~~~~~~~~~~~~~~~~  \left\{ \begin{array} {r}  16 \\ -16  \end{array} \right.  \\~\\
3311&& C_{000aaabc} ~~ \left\{ \begin{array} {l}  \{C^*_{33110000}\},~\{C^*_{31030010 }\},
                                                                                      ~\{C^*_{30100301}\},~\{C^*_{30000113}\} \\
                                                                                             \{C^*_{00130310}\} \\
                                                                                             \{C^*_{01030301} \} \\ 
                                                                                             \{C^*_{01300031}\} \end{array} \right.
             ~~~~           \left\{ \begin{array} {r} 32 \\ 32 \\ -32  \\ -32\end{array} \right. \\
3221&& C_{000aabbc} ~~  \left\{ \begin{array} {l}\{C^*_{32100020 } \},~\{C^*_{30220010}\}, ~\{C_{30100220}\},
                                                                                              ~\{C^*_{30200012}\} \\ 
                                                                                        \{ C^*_{32001200}\} ,~\{C^*_{30021002 }\}~]\\
                                                                                         \{C^*_{32002010}\},~\{C^*_{30122000}\},~\{C^*_{30002210}\},
                                                                                          ~\{C^*_{30102002 }\}  \end{array} \right. 
                                                       ~~~~      \left\{ \begin{array} {r} -16   \\  48   \\ 16  \end{array} \right. \\~\\
32111&& C_{000aabcd} ~~  \left\{ \begin{array} {l}
                         \{ C^*_{32100101}\} ,~\{ C^*_{30020111}\},~\{C^*_{31110200}\},~\{C^*_{31010012}\} \\
                         \{C^*_{32010110}\},~\{C_{31120001}\}, ~\{C_{31000211}\}, ~\{C^*_{30110102 } \}\\                       
                            \{C^*_{32011001}\},~\{C^*_{31021100 }\},~\{C^*_{30011201}\},~\{C^*_{31001102}\}  \\
                              \{C^*_{31201001}\},~\{C^*_{30011120}\},~\{C^*_{30211100}\},~\{C_{31001021}\} \\
                          \{  C^*_{31200110} \},~\{C_{30110021} \} \\
                          \{C^*_{31102100}\},~ \{C^*_{30012011}\}  \end{array} \right. 
                         ~~~~          \left\{ \begin{array} {r} -32 \\32  \\-32\\ 32  \\-32 \\32     \end{array} \right. \\   
 311111&& C_{000abcde}  ~~~~~~ \{C^*_{31111010}\},~\{C^*_{30101111}\}
 ~~~~~~~~~~~~~~~~~~~~~~~~~~~~~~~~~~~~~~~~~~~~\,~~  64 \\
2222&& C_{00aabbcc} ~~ \left\{ \begin{array} {l} \{C^*_{20022200}\} \\
                                                                \{ C^*_{02200022}\}\\ 
                                                              \{  C^*_{20202020} \}_{n=2}   \end{array} 
                                                              \right. 
~~~~~~~~~~~~~~~~~~~~~~~~~~~~~~~~~~~~~~~~~~~~~~~~~~~ 
 \left\{ \begin{array} {r} 8 \\ ~~56 \\ 56 \end{array} \right. \\~\\  
 22211&& C_{00aabbcd}~~  \left\{ \begin{array} {l} \{C^*_{00122210} \},~\{C^*_{02102012} \} \\
                                                                                             \{ C^*_{01022201} \}, ~\{C^*_{02012102} \}   \\
                                                                                             \{ C^*_{22120010 }\},~\{C^*_{20100212 }\}\\
                                                                                             \{C^*_{00212120} \},~\{ C^*_{01202021} \}   
 \end{array} \right. ~~~~~~~~~~~~~~~~~~~~~~~~~~~~~~~~~~~~~~ \left\{ \begin{array} {r} -48  \\ 16 \\ -80 \\ -16\end{array} \right. \\                         
221111&& C_{00aabcde} ~~   \left\{ \begin{array} {l} 
                                            \{  C^*_{22111100} \},~\{C^*_{21021011} \},~ \{C^*_{21101201} \}, ~\{C^*_{ 20011112 } \} \\ 
                                            \{C^*_{02110112} \} \\
                                             \{C^*_{21102011 } \} \\
                                            \{C^*_{21012101} \}_{n=4}  \end{array} \right. 
                                                         ~~     \left\{ \begin{array} {r} -32 \\32 \\ -160
                                                         \\ 96 \end{array} \right. \\
2111111&& C_{00abcdef} ~~~~~~ \{  C^*_{21110111}\}
~~~~~~~~~~~~~~~~~~~~~~~~~~~~~~~~~~~~~~~~~~~~~~\,~~~~~~~~~~~~~ -64
\end{eqnarray*}

  To avoid writing out all $N$ of the circular permutations for the $N=7$ and 8 determinants, we will use the short-hand
 notations $``+ \cdots'' $  and $`` \pm \cdots ''$ to mean, respectively, ``plus all cyclic permutations of the preceding terms'' and ``plus all alternating-in-sign cyclic permutations of the preceding terms''.  Then the determinants are:
  
 \begin{eqnarray*}
{\rm det} [A,B,C,D,E,F]   &=& ~A^6 -B^6 +C^6-D^6 +E^6-F^6 \\
                                        && -~~3~ [~A^4D^2 - B^4E^2+ C^4F^2-D^4A^2+E^4B^2-F^4C^2 ~] \\
                                        && -~~6~ [~A^4BF -B^4CA+C^4DB-D^4EC+E^4FD-F^4AE ~] \\
                                        && -~~6~ [~A^4CE -B^4DF+C^4EA-D^4FB+E^4AC-F^4BD ~] \\                                 
                                        && +~~2~ [~ A^3C^3 -B^3D^3+C^3E^3-D^3F^3+E^3A^3-F^3B^3~] \\                                            
                                         && +~~6~ [~A^3B^2E-B^3C^2F+C^3D^2A-D^3E^2B+E^3F^2C-F^3A^2D ~] \\
                                         && +~~6~[~A^3CF^2 -B^3DA^2 + C^3EB^2 -D^3FC^2 + E^3AD^2 -F^3BE^2~ ] \\
                                         && +~12~ [~ A^3BCD-B^3CDE+C^3DEF-D^3EFA+E^3FAB-F^3ABC~] \\  
                                         && +~12~ [~A^3DEF-B^3EFA + C^3FAB - D^3 ABC +E^3 BCD -F^3CDE~]\\ 
                                         && -~~9~ [~A^2B^2C^2 - B^2C^2D^2 + C^2D^2E^2 - D^2E^2F^2 +E^2F^2A^2 -F^2A^2B^2 ~] \\   
                                         &&+ ~~9~ [~A^2C^2E^2 - B^2D^2F^2  ~] \\                                                                       
                                         && -~18~[~A^2BDE^2 - B^2CEF^2 +C^2DFA^2 - D^2EAB^2 +E^2FBC^2 -F^2ACD^2~] ; \\    ~\\                              
  \det [A,B,C,D,E,F,G]   &=& A^7 + B^7 +C^7+D^7+E^7 +F^7+G^7\\
  &&- ~~ 7~[~A^5BG +A^5CF+A^5DE  + \cdots  ~]\\
  && +~~7 ~[~A^4B^2F +A^4C^2D +A^4BD^2+A^4E^2G+A^4EF^2 +A^4CG^2+ \cdots  ~] \\
   && + ~ 14~ [~A^4BCE+A^4DFG + \cdots  ~] \\
    &&- ~~7  ~[~A^3B^3E +A^3BC^3+A^3D^3F+ \cdots  ~ ]\\
   &&+~14 ~[~A^3B^2G^2 +A^3C^2F^2 + A^3D^2E^2+ \cdots  ~] \\
   &&- ~21 ~[~A^3B^2CD + A^3C^2EG+A^3CD^2G+A^3BE^2F +A^3BDF^2 +A^3EFG^2 +  \cdots  ~] \\
   && +~~7 ~[~A^3BCFG +A^3CDEF+ A^3BDEG+ \cdots   ~]\\
   && - ~~7 ~[~A^2B^2CF^2  +A^2B^2D^2G + \cdots   ~]\\
   && -~14 ~[~A^2B^2CEG + A^2BC^2EF+ A^2D^2EFG + \cdots ~] \\
     &&+~ 35 ~[~A^2B^2DEF  + A^2BC^2DG +A^2BCD^2F+ \cdots  ~] \\
   && -105 ~ABCDEFG  ;
   \end{eqnarray*}
   ~\\~\\~\\~\\~\\~\\~\\~\\~\\~\\~\\~\\~\\~\\~\\~\\~\\~\\~\\~\\~\\~\\~\\
   \begin{eqnarray*}
 \det  [A,B, C,D,E,F,G, H] &=& A^8-B^8+C^8-D^8+E^8-F^8+G^8-H^8\\
 && -~~4~ [~A^6E^2 \pm  \cdots  ]\\
&&-~~8 ~[~A^6BH +A^6DF \pm \cdots] -8 ~[~A^6CG \pm \cdots ]\\
&& +~~8 ~[~A^5B^2G +A^5CD^2 +A^5F^2G +A^5CH^2 \pm  \cdots ] +8~ [~A^5C^2E+A^5EG^2 \pm \cdots ]\\
&& +~16 ~[~A^5BCF + A^5DGH \pm \cdots ] +16 ~[~A^5BDE + A^5EFH \pm \cdots ]\\
&& -~~2~[~A^4 C^4 \pm \cdots ] +6 ~[~A^4E^4 -B^4F^4 +C^4G^4- D^4 H^4  ~]\\
&& -~~8~[~A^4B^3F +A^4D^3H + A^4BF^3 +A^4DH^3 \pm \cdots ]\\
&&-~12~[ ~A^4B^2D^2 + A^4F^2H^2 \pm \cdots]+20 ~[ ~A^4C^2G^2 \pm \cdots ]+20~[ ~A^4B^2H^2 +A^4D^2F^2 \pm \cdots ]\\
&&-~24~[~A^4B^2CE +A^4D^2EG +A^4CEF^2 +A^4EGH^2 \pm \cdots     ]\\
&&~~~~~ -24~[~A^4BC^2D +A^4BDG^2 +A^4C^2FG +A^4FG^2H  \pm \cdots     ]\\
&&~~~~ ~+~8~[~A^4BE^2H + A^4DE^2F \pm \cdots  ] + 8~[~A^4CE^2G \pm \cdots ]  \\
&&+~16~[~A^4BCGH +A^4CDFG \pm \cdots] + 16~[ ~A^4BDFH \pm \cdots ] \\
&&~~~~~ - 48~[ ~A^4BEFG +A^4CDEH \pm \cdots  ] \\
&&+~16~[ ~A^3B^2C^3 +A^3D^2G^3 \pm \cdots  ]  -16~[~ A^3C^2E^3\pm \cdots  ] \\
&& +~32~[~A^3B^3CD +A^3BD^3G +A^3CF^3H+A^3FGH^3 \pm \cdots ] \\
&&~~~~~ +32~[~A^3C^3DH \pm \cdots ] -32~[~A^3C^3EG \pm \cdots ]-32~[~A^3BDE^3 \pm \cdots ] \\
&&-~16~[~A^3B^2CG^2 + A^3C^2D^2G +A^3CF^2G^2 +A^3C^2GH^2 \pm \cdots ] \\
&&~~~~~ +48~[~A^3B^2EF^2 + A^3D^2EH^2 \pm \cdots ] \\
&&~~~~~ +16~[~A^3B^2E^2G +A^3CD^2E^2 +A^3E^2F^2G +A^3CE^2H^2 \pm \cdots ] \\
&&-~32~[~A^3B^2CFH +A^3D^2FGH +A^3BCDF^2 +A^3BDGH^2 \pm \cdots] \\
&&~~~~~ +32~[~A^3B^2DFG+A^3BCD^2H+A^3BF^2GH+A^3CDFH^2 \pm \cdots ]\\
&& ~~~~~ -32~[~A^3B^2DEH +A^3BD^2EF + A^3DEF^2H +A^3BEFH^2 \pm \cdots]\\
&& ~~~~~ +32~[~A^3BC^2EH+A^3DEFG^2+A^3C^2DEF+A^3BEG^2H \pm \cdots ]\\
&& ~~~~~ -32~[~A^3BC^2FG +A^3CDG^2H \pm \cdots]\\
&& ~~~~~ +32~[~A^3BCE^2F +A^3DE^2GH  \pm \cdots ]\\
&&+~64 ~[~A^3BCDEG + A^3CEFGH \pm \cdots] \\
&& +~~8~[ ~A^2B^2C^2F^2 \pm  \cdots]  +56~[~ A^2B^2 D^2E^2 \pm \cdots ]
              + 56~[ ~A^2C^2E^2G^2 - B^2D^2F^2 H^2~]\\
&& -~48~[ ~A^2B^2C^2DH +A^2BD^2FG^2 \pm \cdots ] + 16~[ ~A^2B^2C^2EG +A^2CD^2EG^2 \pm \cdots ] \\
&&~~~~~ -80 ~[ ~ A^2B^2CD^2G + A^2CF^2GH^2 \pm \cdots ] -16~[ ~A^2BC^2DE^2+ A^2BDE^2G^2 \pm \cdots ]\\
&& -~32~ [~A^2B^2CDEF +A^2BD^2EGH + A^2BCEF^2H + A^2DEFGH^2 \pm \cdots ] \\
&&~~~~~ +32~ [~A^2BCEFG^2 \pm \cdots  ] - 160~ [ ~A^2BCE^2GH \pm \cdots       ]\\
&&~~~~~ +96~ [~A^2BDE^2FH - B^2CEF^2GA +C^2DFG^2HB -D^2EGH^2AC~ ] \\
&& -~64~ [~A^2BCDFGH \pm  \cdots] .
\end{eqnarray*}

There are various methods that can be used to check these results.  One way is  to use the identities
\begin{eqnarray*} 
 &&  {\rm det} [1,1, \ldots , 1] = 0 ,\\
   && {\rm det} [0,1, \ldots , 1]  = (-1)^{N-1}(N-1) . 
 \end{eqnarray*} 
 [The second identity follows directly from eqs.(2) and (7)].  The first identity does not however constitute a check for even $N$;  
 since the coefficients in additive multiplets alternate in sign, it is satisfied automatically.  As sums over multiplets, these identities are:

 \begin{eqnarray*} 
   \sum~ n [~\{C^*_{M_0M_1 \cdots M_{N-1}}\} ~] &=& 0 ,~~(N {\rm ~odd}),\\
   \sum~ \frac{n}{N}~[~ \{C^*_{M_0M_1 \cdots M_{N-1}} \} ~]
        \sum_{k=0}^{N-1}(-1)^{k(N-1)}  \delta_{M_k,0}   &=&  (-1)^{N-1}(N-1).
 \end{eqnarray*}
  
Another method, when $N$ is divisible by some integer $d$ and when the coefficients for the smaller $(N/d)$-dimensional determinant are known, is to use the relation
\begin{eqnarray*}
\det [x_0,0,\ldots, 0,x_d,0,\ldots ,0,x_{2d},0,\ldots, 0, x_{N-d},0 ,\ldots ,0] = \left( ~  \det [x_0,x_d,x_{2d},\ldots , x_{N-d}] ~\right)^{d},
\end{eqnarray*}
where the only non-zero matrix elements in the determinant on the left are those of the form $x_{nd~{\rm mod~}N}$.  For example, for $N=6$ and $d=2$,
   \begin{eqnarray*}
  \det [ A,0,C,0,E,0] = \left(~ \det [A,C,E] ~\right)^2 =\left( A^3+C^3+E^3-3ACE  \right)^2.
\end{eqnarray*}
Expanding the right side, all coefficients of the form $C^*_{M_0 0 M_2 0 M_4 0}$ can be read off and compared with the calculated values.\\


\begin{thebibliography}{99}
  \bibitem{Wyn} Wyn-Jones, Alun, {\it Circulants}, Chapters 10 and 11.  www.circulants.org/circ/circall.pdf.
  \bibitem{Ore} Ore, Oyster,  ``Some Studies on Cyclic Determinants'', Duke Mathematical Journal, {\bf 18}, 1951, pp. 343-354.
   \bibitem {Moebius} Wikipedia, {\it M{\"o}bius function}. http://en.wikipedia.org/wiki/M{\"o}bius\_function.
       \bibitem{Abram} Abramowitz and Stegun, {\it Handbook of Mathematical Functions}, Section 24.1, pg. 823.
     \bibitem{oeis}``The Online Encyclopedia of Integer Sequences", published electronically at 
  http://oeis.org, 2010.
  \bibitem{Brualdi} Brualdi and Newman, ``An Enumeration Problem for a Congruence Equation'', 
  J. Res. Nat. Bureau Standards, B74 (1970), 37-40.
     \bibitem{Totient} Weisstein, Eric W., ``Totient Function''. http://mathworld.wolfram.com/TotientFunction.html.
   \bibitem{MMG} Weisstein, Eric W., ``Modulo Multiplication Group.''  http://mathworld.wolfram.com/ModuloMultiplicationGroup.html.
 \end{thebibliography}
\end{document}